\documentclass[12pt]{article}
\usepackage{amsfonts}
\usepackage{amsthm}
\usepackage{CJK}
\usepackage{amssymb}
\usepackage{amsmath}    
\usepackage{graphicx,floatrow}   
\usepackage{verbatim}   
\usepackage{color}      
\usepackage{subfigure}  
\usepackage{epsfig}
\usepackage{float}
\usepackage{epstopdf}
\usepackage{cases}
\usepackage{appendix}
\usepackage{bm}
\usepackage{multirow}
\usepackage{algorithm}
\usepackage{algorithmic}
\usepackage{booktabs}

\allowdisplaybreaks

\newtheorem{theorem}{Theorem}[section]
\newtheorem{lemma}[theorem]{Lemma}
\newtheorem{definition}{Definition}[section]
\theoremstyle{definition}
\newtheorem{remark}{Remark}[section]

\numberwithin{equation}{section}
\renewcommand{\theequation}{\arabic{section}.\arabic{equation}}
\renewcommand{\thetable}{\arabic{section}.\arabic{table}}
\renewcommand{\thefigure}{\arabic{section}.\arabic{figure}}

\newcommand{\abs}[1]{\left\vert#1\right\vert}

\begin{document}
\title{Two-phase segmentation for intensity inhomogeneous images by the Allen-Cahn Local Binary Fitting Model}
\date{}
\author{Chaoyu Liu \thanks{Department of Applied Mathematics, The Hong Kong Polytechnic University, Hong Kong, E-mail: polyucy.liu@connect.polyu.hk)}
	\and Zhonghua Qiao \thanks{Department of Applied Mathematics \& Research Institute for Smart Energy, The Hong Kong Polytechnic University, Hong Kong, E-mail:{zhonghua.qiao@polyu.edu.hk})}
	\and Qian Zhang\thanks{Corresponding author. Department of Applied Mathematics, The Hong Kong Polytechnic University, Hong Kong, E-mail:{qian77.zhang@polyu.edu.hk})}}

\maketitle

\begin{abstract}
This paper proposes a new variational model by integrating the Allen-Cahn term with a local binary fitting energy term for segmenting images with intensity inhomogeneity and noise. An inhomogeneous graph Laplacian initialization method (IGLIM) is developed to give the initial contour for two-phase image segmentation problems. To solve the Allen-Cahn equation derived from the variational model, we adopt the exponential time differencing (ETD) method for temporal discretization, and the central finite difference method for spatial discretization. The energy stability of proposed numerical schemes can be proved.
  Experiments on various images demonstrate the necessity and superiority of proper initialization and verify the capability of our model for two-phase segmentation of images with intensity inhomogeneity and noise.
\end{abstract}


\section{Introduction}
\label{sec1}

Image segmentation aims to divide an image domain into
disjoint areas according to a characterization of the image within or in-between the regions. It plays a crucial role in computer vision, pattern recognition, and has many applications in the field of medical image recognition, satellite remote sensing, and visual field monitoring \cite{Abburu,mitiche2010variational,Pham}.

Various approaches have been developed for image segmentation, wherein active contour models are of particular interest. In the past few decades, active contour models have been widely used to detect edges of objects in images \cite{CasellesVicent1997Geodesic,chan2001active,GUO2021108013,Kass1988Snakes,ma2020fast,niu2017robust}. The basic idea is to create an initial contour and then drive it to evolve to the edges of objects according to certain information from the image. Generally, the information relies on the edge \cite{Caselles1993A,CasellesVicent1997Geodesic,Kass1988Snakes} or the region \cite{chan2001active,li2007implicit,mumford1989optimal,zhang2013local} of the given image. Edge-based active contour models take stopping functions that mainly depend on the gradient information, which makes results very sensitive to the initialization, boundary strength and noise \cite{li2011level,min2021inhomogeneous}. In contrast with edge-based active contour models, global region-based active contour models, such as the Mumford-Shah model \cite{mumford1989optimal} and the Chan-Vese (CV) model \cite{chan2001active,vese2002multiphase}, have better performance on images with noise and weak boundaries. Nevertheless, most of global region-based active contour models are not applicable to images with intensity inhomogeneity. Therefore, many local region-based models \cite{li2011level,li2007implicit,li2008minimization,zhang2013local} have been proposed in the last few years. In these models, contributors of the fitting energy at each pixel are mainly from pixels around it, which can effectively reduce the influence of intensity inhomogeneity. For example, Li et al. proposed a local binary fitting (LBF) model \cite{li2007implicit} and a local intensity clustering (LIC) model \cite{li2011level}. Both of them can segment inhomogeneous images effectively. Recently, Zosso et al. \cite{zosso2017image} proposed the CV-XB model, which integrates an artifact indicator function $X$ and a smooth bias field term into the Chan-Vese model and thus can successfully correct the bias and segment images with noise and intensity inhomogeneity. In \cite{min2021inhomogeneous}, Min et al. developed a model by integrating the smooth bias field term into the LIC model. This model utilizes the local constant and global smoothness priors to describe the bias field and thus can give quite exact segmentation results. Moreover, many hybrid models employ both local and global energy fitting terms to achieve more delicate segmentation for various images, see e.g., \cite{akram2017active,song2020active,wang2009active} and references therein.

Although more and more active contour models with better performance have been designed for image segmentation, the initialization and noise effect are still significant issues. A proper selection for the initial contour can increase the probability of successful segmentation and reduce the segmentation time, while an improper initialization can seriously affect the model performance and lead to unacceptable segmentation results. Generally, initial contours are selected by simply taking a threshold value of the original image \cite{jung2007multiphase,li2011multiphase,yang2019image}, or selecting a part of the image 
\cite{chan2001active,li2007implicit,zhang2013local}. However, these initialization methods are not convenient or flexible and likely to give improper initial contours. 
When an improper initialization is applied to non-convex energy functionals, the minimization process can easily get stuck in poor local minima.
To overcome numerical difficulties caused by the non-convexity, convex relaxation has been widely studied for active contour models  \cite{bae2015efficient,bae2011global,bresson2007fast,brown2010convex,Cai_Zeng,chan2006algorithms,pock2009algorithm}. It can transfer the original energy minimization problem to a convex optimization problem and then find the global minima. However, the convex relaxation may lead to the loss of non-convex boundary information, which makes it difficult to preserve the sharpness and neatness of edges \cite{chan2018convex,wu2021two}. 
As a consequence, it is well known that initialization is a vital step for active contour models. In this work, inspired by the nonlocal edge detection initialization method in \cite{qiao2021twophase}, we construct an anisotropic Laplacian operator to provide a proper initial contour for segmenting images with inhomogeneity.

In addition to the initialization problem, the robustness and efficiency of a segmentation model is also a challenge for images with severe intensity inhomogeneity and noise. In \cite{DongWang,2019The}, Wang et al. used a concave functional of characteristic functions of segments to approximate the contour length and then proposed the iterative convolution-thresholding method (ICTM) to minimize modified energy functionals in an efficient and energy stable way. It is applicable to a range of active contour models. But for noise images, the results of ICTM are not so satisfactory if the given active contour models are sensitive to the noise. In \cite{niu2017robust}, Niu et al. proposed a region-based model
via local similarity factor (RLSF), which has excellent performance on images with strong noise. But
it can not have a very delicate result for images with severe intensity inhomogeneity.
 It is worth mentioning that due to the $\Gamma$-convergence theory \cite{jung2007multiphase,modica1987gradient,modica1977esempio}, the Allen-Cahn functional, which is non-convex, has also been utilized to approximate the length term in phase-field models \cite{esedog2006threshold,zosso2017image}.
 With the Allen-Cahn term, some efficient algorithms can be designed based on the MBO method \cite{merriman1992diffusion,merriman1994motion} or the ETD method \cite{qiao2021twophase}. However, it is still an intractable problem for these algorithms to handle images with both intensity inhomogeneity and strong noise.

In this paper, we employ the Allen-Cahn term to approximate the contour length and propose Allen-Cahn local binary fitting (ACLBF) model based on the LBF energy for image segmentation.  {It is widely known that the phase-field model can be approximately attributed to the solution of a diffuse interface problem. As a result, the smooth effect of the Allen-Cahn equation, which stems from its diffusion property, enables our segmentation model to reduce the influence of noise significantly. The LBF energy functional in our model can sharpen the gradient along the boundaries. For the majority of points, the LBF term will be forced to evolve to the correct phase. For some noise points forced to the incorrect phase, the Allen-Cahn term helps to bring them back to correct phase after smoothing. Combining these two terms, we obtain our desired segmentation model. A good balance between the LBF term and Allen-Cahn term can be obtained through adjusting corresponding parameters so that all noise can be removed, at the same time, the accuracy near the boundary can be guaranteed.} Then we introduce an IGLIM based on edge detection to tackle the initialization problem for two-phase image segmentation. With a proposed anisotropic Laplacian operator, our initialization method can generate a more reliable initial contour even for images with intensity inhomogeneity and noise. Both first- and second-order ETD schemes are designed to solve the evolving equation derived in our model. These ETD schemes can be implemented efficiently on a uniform mesh by Discrete Cosine Transform (DCT).  Using ETD schemes, we solve the contour evolution equation directly, which can exploit the advantage of the Allen-Cahn term to reduce the effect of noise. Meanwhile, the energy stability of our schemes can be proved when the stabilizer satisfies a certain condition by using similar techniques developed in \cite{du2019maximum,DuQiaoSIREV}. Experiments show that our methods can achieve the segmentation in an effective and efficient way, even for the images with severe intensity inhomogeneity and strong noise. %
  {This paper utilizes the ETD1 and ETDRK2 schemes for time discretization. For more details about higher-order ETD schemes, one can refer to \cite{ju2015fast2,zhu2016fast}. In addition, some ETD schemes of arbitrary accuracy have recently been proposed. We recommend interested readers to \cite{chen2021energy, Li_Yang_Zhou}.}

The rest of this paper is organized as follows. In section \ref{sec2}, we give a detailed description of the proposed ACLBF model. IGLIM is introduced to generate the initial contour and an alternating minimization method is used to solve the ACLBF model in section \ref{sec3}.  Numerical examples are given in section \ref{sec4} to show the performance of the developed algorithms. Finally, the paper ends with some conclusions in section \ref{sec5}.

\section{Allen-Cahn local binary fitting model}
\label{sec2}
In this section, a phase-field approach to the LBF model is considered to achieve segmentation for images with intensity inhomogeneity and noise.


As shown in \cite{li2007implicit}, the LBF energy functional is given based on the level set method. In the LBF model, each pixel $x\in \Omega$ in an image is equipped with a local binary fitting energy, which is defined as follows:
\begin{equation*}
	\begin{aligned}
		e_x^{LBF} (\phi, f_1(x), f_2(x))
	&= \lambda_1 \int_\Omega H(\phi(y))K_{\sigma}(x - y)|I(y) - f_1(x)|^2dy\\
	&+ \lambda_2\int_\Omega (1-H(\phi(y)))K_{\sigma}(x - y)|I(y) - f_2(x)|^2dy.
	\end{aligned}
\end{equation*}
Here, $\phi(x)$ is the level set function and the active contour is represented by its zero level set. $\lambda_1$ and $\lambda_2$ are non-negative parameters. $H(x)$ is the Heaviside function and $I: \Omega \to R$ is a given gray level image. $f_1(x)$ and $f_2(x)$ are functions used to fit image intensities near the point $x$ and $K_{\sigma}(x)$ is the Gaussian kernel:
\begin{equation*}
	K_\sigma(x) = \frac{1}{(2\pi)\sigma^2} e^{-\frac{|x|^2}{2\sigma^2}},  \ \sigma> 0.
\end{equation*}
Then the energy functional of LBF model is given by
\begin{align}\label{LBF}
E^{LBF}(\phi,f_1,f_2) = \mu e^{LBF}  + \int_{\Omega} \delta(\phi(x))|\nabla \phi(x)| dx + \nu \int_{\Omega} \frac{1}{2}(\abs{\nabla\phi}-1)^2dx,
\end{align}
where the fitting term $e^{LBF}\left(\phi, f_{1}(x), f_{2}(x)\right)$ is the integral of $e_x^{LBF}$ over $\Omega$:
\begin{align*}
	e^{LBF}\left(\phi, f_{1}(x), f_{2}(x)\right)& = \int_{\Omega} e_x^{LBF}\left(\phi, f_{1}(x), f_{2}(x)\right)dx\\
	&= \lambda_1 \int_\Omega\int_\Omega H(\phi(y))K_\sigma(x - y)|I(y) - f_1(x)|^2dy dx\\
	&+ \lambda_2\int_\Omega\int_\Omega (1-H(\phi(y)))K_\sigma(x - y)|I(y) - f_2(x)|^2dy dx.
\end{align*}
The second term of \eqref{LBF} formulated by $\phi(x)$ is an approximation of the contour length and $\delta(x)$ represents the Dirac delta function, which is the derivative of $H(x)$ in the distribution sense. The last term of \eqref{LBF} is a penalty term to constrain $|\nabla \phi| = 1$, which is used for the re-initialization. $\mu$ is a non-negative parameter used to control the contribution of the fitting energy to the energy functional and $\nu$ is a non-negative parameter used to tune the penalty force.

 By introducing the Allen-Cahn functional to the LBF model, we propose the ACLBF model to solve segmentation problems of images with intensity inhomogeneity and noise. The energy functional of the ACLBF model is defined as follows:
\begin{equation}\label{energyfunctional}
	E(u,f_1,f_2) = \mu e^{LBF}(u,f_1,f_2) + \int_{\Omega}\left(\frac{\varepsilon}{2}|\nabla u|^{2} + \frac{1}{\varepsilon} W(u)\right) dx.
\end{equation}
Here, $\varepsilon>0$ is a diffusion parameter and $W(u)$ is given by:
$$
W(u) = \sin^2(\frac{\pi}{2}(u+1)).
$$
As pointed in \cite{jung2007multiphase}, for binary phase-fields $u\in\{-1,1\}$,
$$
\int_{\Omega}\left(\frac{\varepsilon}{2}|\nabla u|^{2} + \frac{1}{\varepsilon} \sin^2(\frac{\pi}{2}(u+1))\right) dx \xrightarrow[\Gamma]{\varepsilon \rightarrow 0^{+}} C \int_\Omega |\nabla u| dx,
$$
 where $C$ is a constant independent on $u$ and $\int_\Omega |\nabla u| dx$ is a common approximation of the contour length \cite{chan2006algorithms}.
 The two terms in \eqref{energyfunctional} correspond to the first two terms in the LBF energy functional \eqref{LBF}. The last penalty term in \eqref{LBF} can be removed because re-initialization is not required for the phase-field model. In the ACLBF model, the segmented edges are represented by the zero level set of $u$.

The advantages of the ACLBF model can be concluded as follows. Firstly, the LBF fitting energy term can help to segment images with intensity inhomogeneity effectively. Secondly, the phase-field approach can reduce the noise effect significantly in comparison with the level set method. Moreover, the re-initialization is not required.
Finally, we can easily design energy stable numerical methods for this phase-field model.


\section{The numerical scheme}
\label{sec3}

In this section, we will give a detailed description of the numerical scheme proposed for the ACLBF model. An alternating minimization method will be used to minimize the energy functional \eqref{energyfunctional}. As mentioned in section \ref{sec1}, how to give the initial data is a vital step for segmentation of images with intensity inhomogeneity.  To tackle the initialization issue, we introduce IGLIM, which can find exact partial edges in images and subsequently select them as initial contours.

\subsection{Inhomogeneous graph Laplacian initialization method}
\label{subsec31}

\cite{qiao2021twophase} shows that it is very effective to solve region-based active contour models for image segmentation when taking edges detected by gradient-based detection methods as the initial contour. Inspired by this, we propose IGLIM to generate the initial contour for our numerical scheme to solve the ACLBF model for inhomogeneous images.
Generally, gradient-based edge detection methods can be grouped into two categories \cite{castleman1998digital}. Some are based on the first-order derivative \cite{cherri1989optical} and the others are based on the second-order derivative \cite{canny1986computational, marr1980theory}. These edge detection methods are early techniques for detecting edges in images and easy to implement. For second-order derivative based edge detection methods, signs of the Laplacian values defined on pixels will change through the edge due to the rapid change of the image intensity. Therefore, the edge should consist of all these zero-cross points of Laplacian. While for images with intensity inhomogeneity, classical discrete Laplacian operators may fail to obtain correct edges of objects in images since  the image intensity may not change so rapidly. Other traditional gradient-based edge detection operators can not give satisfactory results for inhomogeneous images either.  
We will introduce the IGLIM algorithm below,  which can give appropriate initial contours for region-based active contour models. This initialization method contains two steps. Firstly, we introduce an inhomogeneous Laplacian operator by which most edges can be determined even for images with intensity inhomogeneity. Then a denoising method is applied to remove the misclassification caused by noise.

\subsubsection{Inhomogeneous graph Laplacian operator}

Let $\Omega$ be a 2D discrete image domain, and $I$ be an image defined on it with $M_1\times M_2$ pixels. Denote a pixel $x_0 = (i,j)\in \Omega$, $I(x_0)$ by $I_{i,j}$. We define the inhomogeneous graph Laplacian operator $L$ as
\begin{equation}\label{Graph_L}
L(x_0)= \sum_{k=1}^{8}c_kI_{i,j}^k -I_{i,j},
\end{equation}
where $I_{i,j}^k$ is intensity value of the $k$-th neighbour point of $I_{i,j}$, more precisely,
$$
\begin{aligned}
	& I_{i,j}^1 = I_{i-1,j-1}, I_{i,j}^2 = I_{i-1,j}, I_{i,j}^3 = I_{i-1,j+1}, I_{i,j}^4 = I_{i,j+1}, \\
	& I_{i,j}^5 = I_{i+1,j+1}, I_{i,j}^6 = I_{i+1,j},
	I_{i,j}^7 = I_{i+1,j-1}, I_{i,j}^8 = I_{i,j-1},
\end{aligned}
$$
and
$$
c_k=\frac{e^{\lambda(I_{i,j}-I_{i,j}^k)^2}}{\sum_{k=1}^{8}e^{\lambda(I_{i,j}-I_{i,j}^k)^2}},\text{ $\lambda$ is a given non-negative parameter}.
$$\par
In the expression of inhomogeneous Laplacian operator $L$, $c_k$ is between $0$ and $1$. And the larger difference between $I_{i,j}$ and $I_{i,j}^k$ is, the bigger value of $c_k$ would be. Compared with the classical discrete Laplacian, the inhomogeneous graph Laplacian operator enlarges the influence of points with big ``difference" from the central point $I_{i,j}$ (Points with big ``difference" are more likely to be edge points or noise).
\begin{remark}
  If $\lambda = 0$, then the inhomogeneous graph Laplacian \eqref{Graph_L} degenerates to a homogeneous discrete Laplacian operator, i.e.,
  $$
  \begin{aligned}
    L(x_0) &= \sum_{k=1}^{8}\frac{1}{8}I_{i,j}^k -I_{i,j}\\
    & = \frac{1}{8}(I_{i-1,j}+I_{i+1,j}+I_{i,j-1}+I_{i,j+1}-4I_{i,j})\\
    &+ \frac{1}{8}(I_{i-1,j-1}+I_{i+1,j+1}+I_{i+1,j-1}+I_{i-1,j+1}-4I_{i,j}).
  \end{aligned}
  $$
\end{remark}

We need to approximate the zero-cross points of the inhomogeneous Laplacian operator to obtain rough initial edges. Some relevant definitions are stated below.
\begin{definition} Let $k_1, k_2$ be two small non-negative numbers.
\label{bp}
\leftline{}
1. If $L(x_0)\le-k_1$, $L(x_0)$ is defined as negative.\\
2. If $L(x_0)\ge k_2$, $L(x_0)$ is defined as positive.\\
3. The set of all zero-cross points of $L$ is denoted as
\begin{equation*}
  \begin{aligned}
    S = \{ &x|x \text{ has at least one neighbor point } y \text{ such that the sign of }\\
    & L(y) \text{ is different from } L(x) \}.
  \end{aligned}
\end{equation*}
4. S has following two subsets:
   \begin{equation}
\begin{split}
	&\text{Positive Laplacian edge points set:} \ S_p = \{x|x\in S, L(x) \ge k_2 \},\\
	&\text{Negative Laplacian edge points set:} \ S_n = \{x|x\in S, L(x) \le -k_1 \}.
\end{split}
\end{equation}
\end{definition}

The boundary of the object can be divided into an inner boundary and an outer boundary. The difference is that the pixels on the inner boundary belong to the object while those on the outer boundary belong to the background. In our method, the edge points are divided into two groups. One of $S_p$ and $S_n$ consists of pixels on the inner boundary and the other one consists of pixels on the outer boundary. Generally, we choose the inner boundary to be the rough initial contour.
\begin{remark}
  \label{bd_judge}
 One can determine which of $S_p$ and $S_n$ corresponds to the inner boundary according to the intensity information of the object. If the intensity of the object is smaller than that of the background, i.e. the object is darker than the background, then $S_p$ consists of pixels on the inner boundary. Otherwise, $S_n$ corresponds to the inner boundary.
\end{remark}

\subsubsection{A denoising method based on the connectivity of edge points}

Although the inhomogeneous graph Laplacian operator can give rough edges for images with intensity inhomogeneity, noise in images can affect its performance heavily. To solve this problem, a denoising method is proposed here to remove the possible noise pixels in the rough initial contour obtained from the inhomogeneous graph Laplacian operator.
This denoising method is motivated by the fact that an edge should have connectivity, which means that points of edges connect with each other, while the noise doesn't possess this property. As a result, we can remove most noise points from the rough initial contour by judging their connectivity. In our denoising method, we mainly consider the diagonal connectivity of edge points.
The diagonally connected points are defined as follows:
\begin{definition}
Suppose that $S_p$ ($S_n$) consists of pixels on the inner boundary.
When $x = (i,j)\in S_p$ ($S_n$), the neighbor areas are divided into four parts:
\begin{align*}
	&S_1 = \{(i-1,j-1),(i,j-1),(i-1,j)\}, \ S_2 = \{(i-1,j+1),(i,j+1),(i-1,j)\},\\
	&S_3 = \{(i+1,j-1),(i,j-1),(i+1,j)\}, \ S_4 = \{(i+1,j+1),(i,j+1),(i+1,j)\}.
\end{align*}
We call $x\in S_p$ ($S_n$) a diagonally connected point if both $S_1$ and $S_4$ or both $S_2$ and $S_3$ have at least one pixel that also belongs to $S_p$ ($S_n$).
\end{definition}

To eliminate noise in the rough initial contour, we keep all the diagonally connected points in $S_p$ ($S_n$) and remove the other points from $S_p$ ($S_n$). The denoising process needs to be repeated $M$ times where $M$ is a pre-setting small integer.
\begin{remark}
	All the edge points are connected with each other but not all of them are diagonally connected, and thus a few of them will be removed from the rough initial contour points set $S_p$ ($S_n$) after denoising. But the majority of edge points will remain in $S_p$ ($S_n$). Meanwhile, a few noise points can be diagonally connected with each other but most of them are not, and thus most of noise points will be removed from $S_p$ ($S_n$) after denoising. The reason why we consider the diagonal connectivity rather than the common connectivity is that the former one is less likely to appear on noise points, which implies that we can remove more noise points by judging the diagonal connectivity.
\end{remark}

 Now combining the inhomogeneous graph Laplacian operator and the denoising method yields the IGLIM. The algorithm for IGLIM is organized in \ref{IGLIM}:\\
 \begin{algorithm}
   \caption{IGLIM}
   \label{IGLIM}
   \begin{algorithmic}
   \STATE{\textbf{Step 1}: Compute the inhomogeneous Laplacian value of each pixel by \eqref{Graph_L}.}
   \STATE{\textbf{Step 2}: Set $k_1$ and $k_2$. Determine $S_p$ and $S_n$ according to  Definition \ref{bp}.}
   \STATE{\textbf{Step 3}: Take $S_p$ $(S_n)$ as a rough initial contour.}
   \STATE{\textbf{Step 4}: Go through every pixel in $S_p$ $(S_n)$ and judge whether it is diagonally connected. If not, remove it from $S_p$ $(S_n)$.}
   \STATE{\textbf{Step 5}: Set an appropriate integer $M$, and repeat Step 4 for $M$ times.}
   \STATE{\textbf{Step 6}: Output $S_p$ $(S_n)$ as the initial contour.}
   \end{algorithmic}
 \end{algorithm}

 \subsubsection{Initialization of the ACLBF model}
 It can be seen that the $S_p$ or $S_n$ given by IGLIM is a set of curves. However, curves are not the best choice for the initialization of our model, for which a single curve is prone to disappear quickly due to the existence of the diffusive interface term $\int_{\Omega}\frac{\varepsilon}{2}|\nabla u|^2 d x$. Therefore, we extend all points in the initial contour points set $S_p$ $(S_n)$ from one single point to a small region. The way we extend a point in $S_p$ $(S_n)$ is to add into it all its neighbor points $y$ whose inhomogeneous graph Laplacian values are not negative (positive), i.e., $y\notin S_n$ ($S_p$). \\
 Let
 $$
 R_p = \{y| y\notin S_n\text{ and } y \text{ has at least one neighbor point }x\in S_p\},
 $$
 $$
 R_n = \{y| y\notin S_p\text{ and } y \text{ has at least one neighbor point }x\in S_n\},
 $$
 and then the final region extended from $S_p$ $(S_n)$ is $S_p\cup R_p$ ($S_n\cup R_n$).
 Finally, the initial value $u_0$ for the ACLBF model is given by
 \begin{equation}
 	u_{0}(i,j) = \left\{\begin{array}{c l}
 	1, &\text{if } (i,j) \in S_p\cup R_p\text{ } (S_n\cup R_n),\\
 	-1,&\text{otherwise.}
 \end{array}\right. \label{u0}
 \end{equation}

\subsection{Energy minimization}
We will adopt an alternating minimization method to minimize the energy functional of the ACLBF model \eqref{energyfunctional}. The initial contour obtained in \ref{subsec31} will be used to start the alternating minimization iteration.

\subsubsection{An iterative method for energy minimization}
In this part, we will demonstrate the specific procedure of our iterative method for minimizing the energy functional (\ref{energyfunctional}). It consists of two parts. We first fix $u$ and minimize \eqref{energyfunctional} with respect to the functions $f_1(x)$ and $f_2(x)$.
By variation calculus, one can show that the functions $f_1(x)$ and $f_2(x)$ are given by
\begin{equation}\label{f1f2}
	\displaystyle f_{1}(x)=\frac{K_{\sigma}(x) *\left[H\left(u(x)\right) I(x)\right]}{K_{\sigma}(x) * H\left(u(x)\right)} , \quad
	\displaystyle f_{2}(x)=\frac{K_{\sigma}(x) *\left[\left(1-H\left(u(x)\right)\right) I(x)\right]}{K_{\sigma}(x) *\left[1-H\left(u(x)\right)\right]}.
\end{equation}
In calculation, the Heaviside function $H(x)$ is approximated by the following smooth function
\begin{equation}\label{Heviside}
	H_{\varepsilon_1}(x) = \frac{1}{2}[1+\frac{2}{\pi}\arctan(\frac{x}{\varepsilon_1})].
\end{equation}
Correspondingly, the function used to approximate $\delta(x)$ is defined as:
\begin{equation}\label{delta}
	\delta_{\varepsilon_1} (x)= H'_{\varepsilon_1}(x)=\frac{1}{\pi}\frac{\varepsilon_1}{\varepsilon_1^2 + x^2}.
\end{equation}

Next, keeping $f_1$ and $f_2$ fixed, and minimizing the energy functional \eqref{energyfunctional} with respect to $u$, lead to the Allen-Cahn equation:
\begin{equation}
	u_{t}=\varepsilon \Delta u-\frac{1}{\varepsilon} W^{\prime}(u)-\mu\delta_{\varepsilon_1}\left(u\right)\left(\lambda_{1} e_{1}-\lambda_{2} e_{2}\right),\label{evo_eq}
\end{equation}
where
\begin{equation*}
	e_{1}=\displaystyle\int_{\Omega} K_{\sigma}(y-x)\left|I(x)-f_{1}(y)\right|^{2} d y, \quad
	e_{2}=\displaystyle\int_{\Omega} K_{\sigma}(y-x)\left|I(x)-f_{2}(y)\right|^{2} d y.
\end{equation*}
Let $x\in\Omega$, the value of $e_k$ at $x$ is evaluated by:
\begin{align*}
  e_k(x) & = (K_\sigma*f_k^2)(x)-2I(x)(K_\sigma*f_k)(x) + I^2(x)(K_\sigma*\textbf{1}_\Omega)(x), \quad k = 1,2,
\end{align*}
where $\textbf{1}_\Omega$ is the characteristic function of $\Omega$ and $*$ represents the discrete convolution operator.

\subsubsection{Exponential time differencing method}
\label{subsec32}
To solve the Allen-Cahn equation (\ref{evo_eq}) efficiently and accurately, we will use the ETD methods for temporal discretization and the central finite difference method for spatial discretization with the homogeneous Neumann boundary condition.

Let $h$ be the spacing distance between two adjacent pixels. 
By taking the central finite difference discretization of (\ref{evo_eq}) in space, we obtain the following ordinary differential equations (ODE) system:

\begin{equation}\label{ODE}
	U_{t}=-\boldsymbol{L}_h U+N(U),
\end{equation}
where
$$
\boldsymbol{L}_h = S\boldsymbol{I} - \varepsilon D_h, \quad \ N(U) = SU-\frac{1}{\varepsilon} W^{\prime}(U)-\mu\delta_{\varepsilon_1}\left(U\right)\left(\lambda_{1} e_{1}-\lambda_{2} e_{2}\right).
$$
Here, $U = [U_k]_{k=1}^{M_1M_2}\in R^{M_1M_2}$ is the semi-discrete numerical solution after spatial discretization using column-wise ordering and the $k$-th equation of \eqref{ODE} corresponding to $(i,j)-$th point in $\Omega$ has the following relation:
$$
k = i + M_1(j-1), \quad i = 1,2,\cdots,M_1,\quad j = 1,2,\cdots,M_2.
$$
$D_h$ is the 2D discrete Laplacian matrix obtained from the central finite difference discretization of $\Delta$:
$$
  D_{ h } = \frac{1}{h^2}(I_{M_2}\otimes\Lambda_{M_1}+\Lambda_{M_2}\otimes I_{M_1}),
$$
where $I_{M_i}$ is an $M_i\times M_i$ identity matrix and
$$
\Lambda_{ M_i } = \left[ \begin{array} { c c c c c c } - 1 & 1 & & & & 0\\ 1 & - 2 & 1 & & & \\ & & \ddots & \ddots & \ddots & \\ & & & 1 & - 2 & 1 \\ 0 & & & & 1 & - 1 \end{array} \right]_{ M_i \times M_i }, \quad i=1, 2.
$$
$\boldsymbol{I}$ is an $(M_1M_2)\times (M_1M_2)$ identity matrix identity matrix. $S>0$ is a constant called the stabilizer. 

Solving this ODE system, we obtain
\begin{equation*}
	  U\left(t_{n+1}\right) = \mathrm{e}^{-\boldsymbol{L}_h \Delta t} U\left(t_{n}\right)+\int_{0}^{\Delta t} \mathrm{e}^{-\boldsymbol{L}_h(\Delta t-s)} N\left(U\left(t_{n}+s\right)\right) \mathrm{d} s.
\end{equation*}
If we approximate $N\left(U(t_{n}+s)\right)$ by $N(U(t_n))$, then we obtain the first-order ETD (ETD1) scheme:
\begin{equation}
  \label{etd1}
	U^{n+1}=\mathrm{e}^{-\boldsymbol{L}_h \Delta t} U^{n}+\Delta t \phi_{0}(\boldsymbol{L}_h \Delta t) N\left(U^{n}\right),
\end{equation}
where
\begin{equation*}
	\phi_{0}(\boldsymbol{L}_h \Delta t)=\int_{0}^{\Delta t} \mathrm{e}^{-\boldsymbol{L}_h(\Delta t-s)} \mathrm{d} s=(\boldsymbol{L}_h \Delta t)^{-1}\left(\boldsymbol{I}-\mathrm{e}^{-\boldsymbol{L}_h \Delta t}\right). \label{phi_0}
\end{equation*}
If we approximate $N\left(U\left(t_{n}+s\right)\right)$ by a linear approximation
\begin{equation*}
	N\left(U\left(t_{n}+s\right)\right) \approx\left(1-\frac{s}{\Delta t}\right) N\left(U\left(t_{n}\right)\right)+\frac{s}{\Delta t} N\left(\tilde{U}\left(t_{n+1}\right)\right), \quad s \in[0, \Delta t],
\end{equation*}
where $\tilde{U}(t_{n+1})$ is an approximation of $U(t_{n+1})$ obtained by (\ref{etd1}), we obtain the second-order ETD Runge-Kutta (ETDRK2) scheme:
\begin{equation}
  \label{etd2}
	\left\{\begin{aligned}
	\tilde{U}^{n+1} &=\mathrm{e}^{-\boldsymbol{L}_h \Delta t} U^{n}+\Delta t \phi_{0}(\boldsymbol{L}_h \Delta t) N\left(U^{n}\right) \\
	U^{n+1} &=\mathrm{e}^{\boldsymbol{L}_h \Delta t} U^{n}\\
    &+\Delta t\left\{\phi_{0}(\boldsymbol{L}_h \Delta t) N\left(U^{n}\right)
   +\phi_{1}(\boldsymbol{L}_h \Delta t)\left[N\left(\tilde{U}^{n+1}\right)-N\left(U^{n}\right)\right]\right\},
	\end{aligned}\right.
\end{equation}
where
\begin{equation*}
	\phi_{1}(\boldsymbol{L}_h \Delta t)=\int_{0}^{\Delta t} \frac{s}{\Delta t} \mathrm{e}^{-\boldsymbol{L}_h(\Delta t-s)} \mathrm{d} s=(\boldsymbol{L}_h \Delta t)^{-2}\left(\boldsymbol{L}_h \Delta t-\boldsymbol{I}+\mathrm{e}^{-\boldsymbol{L}_h \Delta t}\right).
\end{equation*}
In fact, the ETD schemes can be solved by Discrete Cosine Transform (DCT). For more information see e.g.,\cite{hochbruck2010exponential,ju2015fast}. The time complexity is only $O((M_1M_2)\log (M_1M_2))$ at time step. Consequently, it is remarkably efficient to apply the ETD schemes for solving the evolving equation (\ref{evo_eq}).

The algorithm for minimizing the energy functional is organized in \ref{ACLBF}.
\begin{algorithm}
  \caption{The ETD-based iterative method for ACLBF}
  \label{ACLBF}
  \begin{algorithmic}
  \STATE{\textbf{Step 0}: Generate an initial function $U^0=u_0$ by \eqref{u0} based on IGLIM in \ref{IGLIM}.}
  \STATE{\textbf{Step 1}: Substitute $u=U^n$ into (\ref{f1f2}) and calculate $f_1^n$ and $f_2^n$.}
  \STATE{\textbf{Step 2}: Obtain $U^{n+1}$ from $U^n$ by the ETD1 (or ETDRK2) scheme.}
  \STATE{\textbf{Step 3}: Set $n = n+1$ and repeat Step 1 and Step 2 until the evolution is stationary.}
  \end{algorithmic}
\end{algorithm}

In \ref{ACLBF}, the stopping criterion is that the contour evolution is stationary, which means the contour does not move anymore. In practice, the iteration will be stopped if the contour is identical to that in the last iteration.

\subsubsection{Discrete energy stability}
For given $f_1$ and $f_2$, the Allen-Cahn equation \eqref{evo_eq} holds the energy stability
\begin{align*}
\frac{dE}{dt} \leq 0.
\end{align*}
In the following part, we will show that ETD1 and ETDRK2 schemes can preserve the discrete energy stability. For a rectangular image with $M_1\times M_2$ pixels, order pixels column by column and the $k$-th pixel denoted by $x_k$ corresponding to $(i,j)\in \Omega$ has the following relations:
the discrete energy $E_h (U,f_1,f_2)$ is defined as follows
\begin{align*}
E_h (U,f_1,f_2)
& = \sum_{k=1}^{M_1M_2}\frac{1}{\varepsilon}W(U_k) + \mu\sum_{k=1}^{M_1M_2}(\lambda_1H_{\epsilon_1}(U_{k})e_1(x_k) +\lambda_2(1 - H_{\epsilon_1}(U_{k}))e_2(x_k))\\
& -\frac{\varepsilon}{2}U^TD_hU. 
\end{align*}
\begin{lemma}\label{lemma}
	For any fixed $f_1$ and $f_2$,  when the stabilizer satisfies $S>\frac{G}{2}$ with
\begin{align*}
G := \parallel\tilde{N}'\parallel_\infty, \quad \tilde{N}(U) = SU - N(U),
\end{align*}
we have
	\begin{align*}
    &ETD1: E_h(U^{n+1},f_1,f_2) \le E_h(U^{n},f_1,f_2),\quad \forall \Delta t>0,\\
    &ETDRK2: E_h(U^{n+1},f_1,f_2) \le E_h(U^{n},f_1,f_2) + Ch^{-\frac{1}{2}}(h^2+\Delta t)^2,\ \forall 0<\Delta t<1,
	\end{align*}
where $C\ge0$ is a constant and independent of $\Delta t$ and $h$.
\end{lemma}

\begin{proof}
We only prove the boundedness of $\tilde{N}'$. With the boundedness of $\tilde{N}'$, the remaining part of the proof is similar to the proofs of Theorem 5.1 and Theorem 5.2 in \cite{du2019maximum}.

By a simple calculation, we can obtain
\begin{equation*}
  \begin{aligned}
      \tilde{N}'(U) &= \frac{1}{\varepsilon} W^{\prime\prime}(U) + \mu\delta^{\prime}_{\varepsilon_1}\left(U\right)\left(\lambda_{1} e_{1}-\lambda_{2} e_{2}\right)\\
      & = \frac{1}{\varepsilon} W^{\prime\prime}(U) - \frac{2\mu\varepsilon_1}{\pi}\frac{U}{(\varepsilon_1^2+U^2)^2}\left(\lambda_{1} e_{1}-\lambda_{2} e_{2}\right),
  \end{aligned}
\end{equation*}
where
$$
W^{\prime\prime}(U) = \frac{\pi^2}{2}\cos(\pi(u+1)).
$$
It is easy to check that there exists a constant $G_1 = \frac{\pi^2}{2\varepsilon}+1$ such that
$$
\parallel \frac{1}{\varepsilon}W^{\prime\prime}(U)\parallel_\infty < G_1.
$$
And for any fixed $f_1$ and $f_2$, there exists a positive constant $G_2$, which is independent of $U$ such that
$$
\left\|\frac{2\mu\varepsilon_1}{\pi}\frac{U}{(\varepsilon_1^2+U^2)^2}\left(\lambda_{1} e_{1}-\lambda_{2} e_{2}\right)\right\|_\infty < G_2,\quad \forall U\in(-\infty,\infty).
$$
Therefore, $\parallel\tilde{N}'\parallel_\infty < G_1 + G_2$.
\end{proof}



The theorem below shows that the total discrete energy is stable for both ETD1 and ETDRK2 schemes during the whole iteration process, provided that the stabilizer $S$ is sufficiently large.
\begin{theorem}
	Given $f_1^n$ and $f_2^n$, if $S>\frac{G}{2}$, then we have for ETD1,
		\begin{equation}
		\label{ineq_etd1}	 E_h(U^{n+1},f_1^{n+1},f_2^{n+1}) \le E_h(U^{n},f_1^{n},f_2^{n}),\ \forall \Delta t>0;
        \end{equation}
        and for ETDRK2,
        \begin{equation}
		\label{ineq_etd2} 	 E_h(U^{n+1},f_1^{n+1},f_2^{n+1})
\le E_h(U^{n},f_1^{n},f_2^{n}) + Ch^{-\frac{1}{2}}(h^2+\Delta t)^2,\ \forall 0<\Delta t<1,
		\end{equation}
	where $C\ge0$ is a constant and independent of $\Delta t$ and $h$.
\end{theorem}

\begin{proof}
 Through Lemma \ref{lemma}, we can get
\begin{equation*}
  E_h(U^{n+1},f_1^{n+1},f_2^{n+1}) \le E_h(U^{n},f_1^{n+1},f_2^{n+1}).
\end{equation*}
To obtain the inequality (\ref{ineq_etd1}), it suffices to show that
\begin{equation*}
  E_h(U^{n},f_1^{n+1},f_2^{n+1}) \le E_h(U^{n},f_1^{n},f_2^{n}).
\end{equation*}
Indeed, this is a direct consequence of variation calculus. Thus the inequality (\ref{ineq_etd1}) holds. The discrete energy stability for ETDRK2 (\ref{ineq_etd2}) can be derived in a similar way.
\end{proof}

\section{Experimental results}
\label{sec4}
This section displays experiments to test the energy stability of our algorithm and demonstrate the effectiveness of IGLIM and the ACLBF model for segmenting various images. 
 {We observe that the total iteration number of the ETDRK2 scheme is generally less than that of the ETD1 scheme for each simulation in our experiments and moreover, Qiao et al. find that ETDRK2 has better performance than ETD1 on image segmentation in terms of efficiency and accuracy in \cite{qiao2021twophase}. Therefore, unless otherwise specified, all images are segmented by the ACLBF model with ETDRK2 methods.} 
All numerical experiments are implemented on a laptop with 2.60-GHz CPU, 16GB RAM, and MATLAB R2020a.
\subsection{Energy stability test}
\begin{figure}[t!]
		\subfigure[\scriptsize The original image]{
			\includegraphics[width=0.45\linewidth]{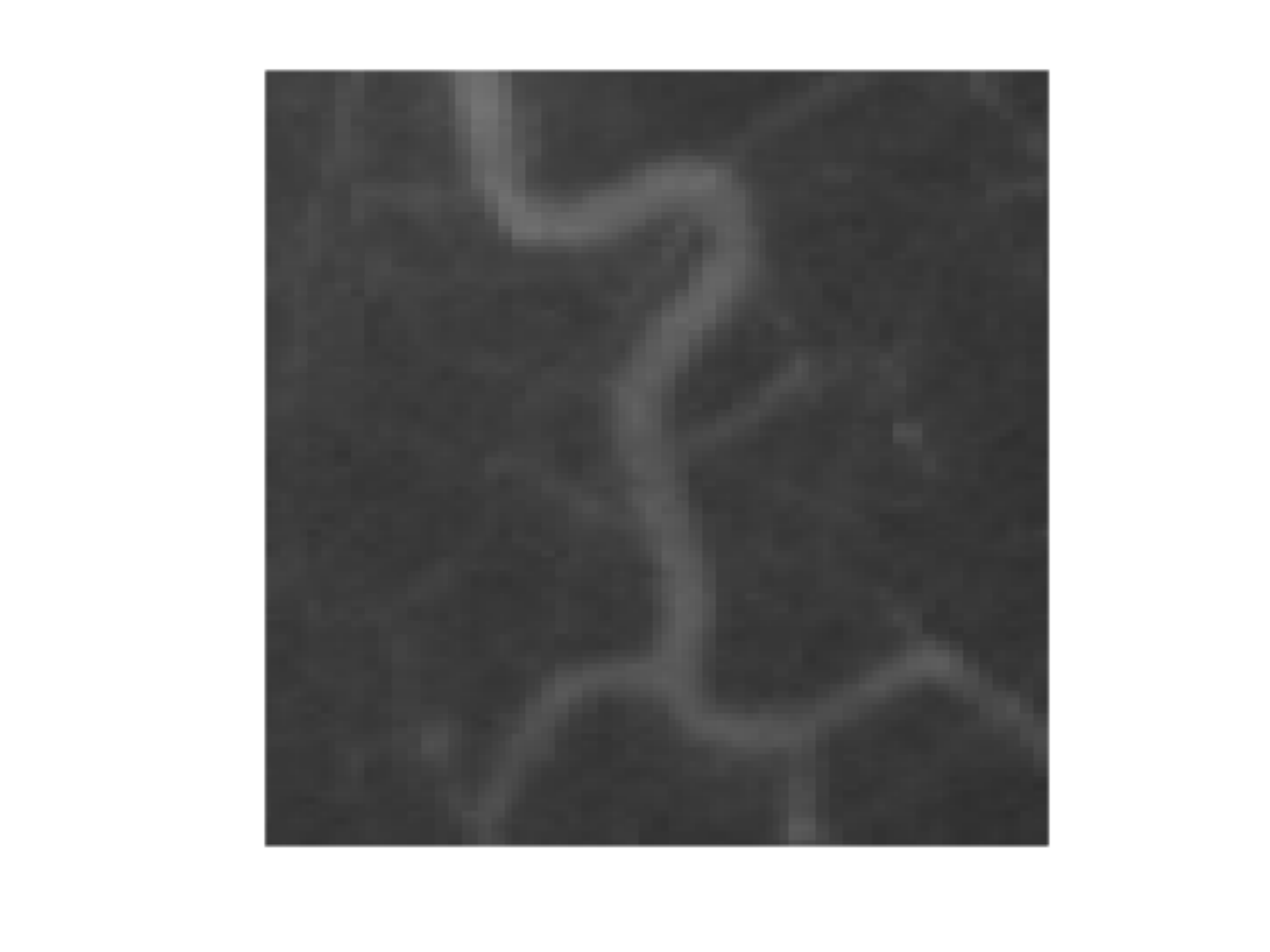}
		}
		\hspace{-2cm}\subfigure[\scriptsize The segmentation result by ETD1 (ETDRK2)]{
			\includegraphics[width=0.45\linewidth]{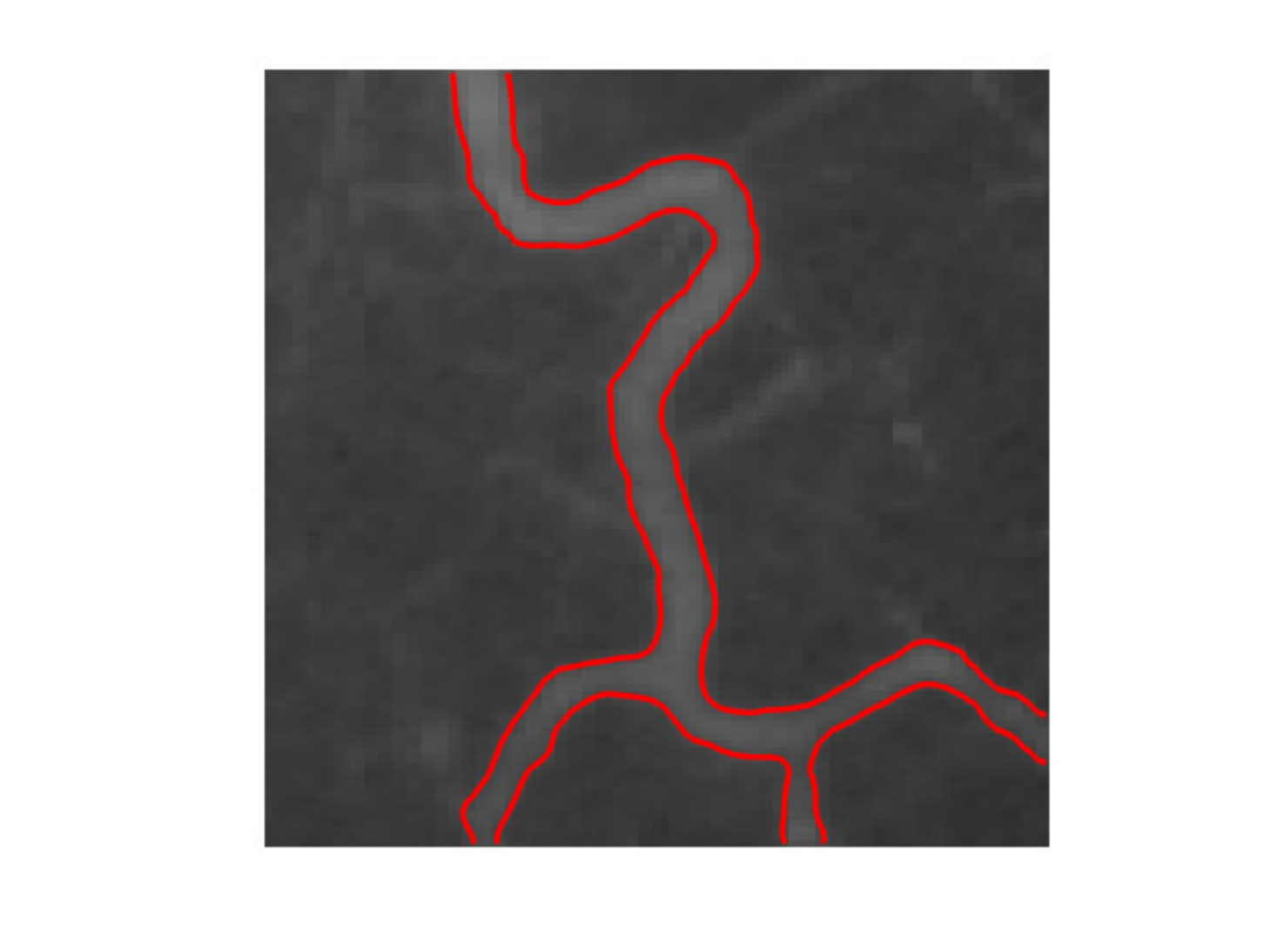}
		}
		\subfigure[\scriptsize Energy diagram of ETD1]{
			\includegraphics[width=0.45\linewidth]{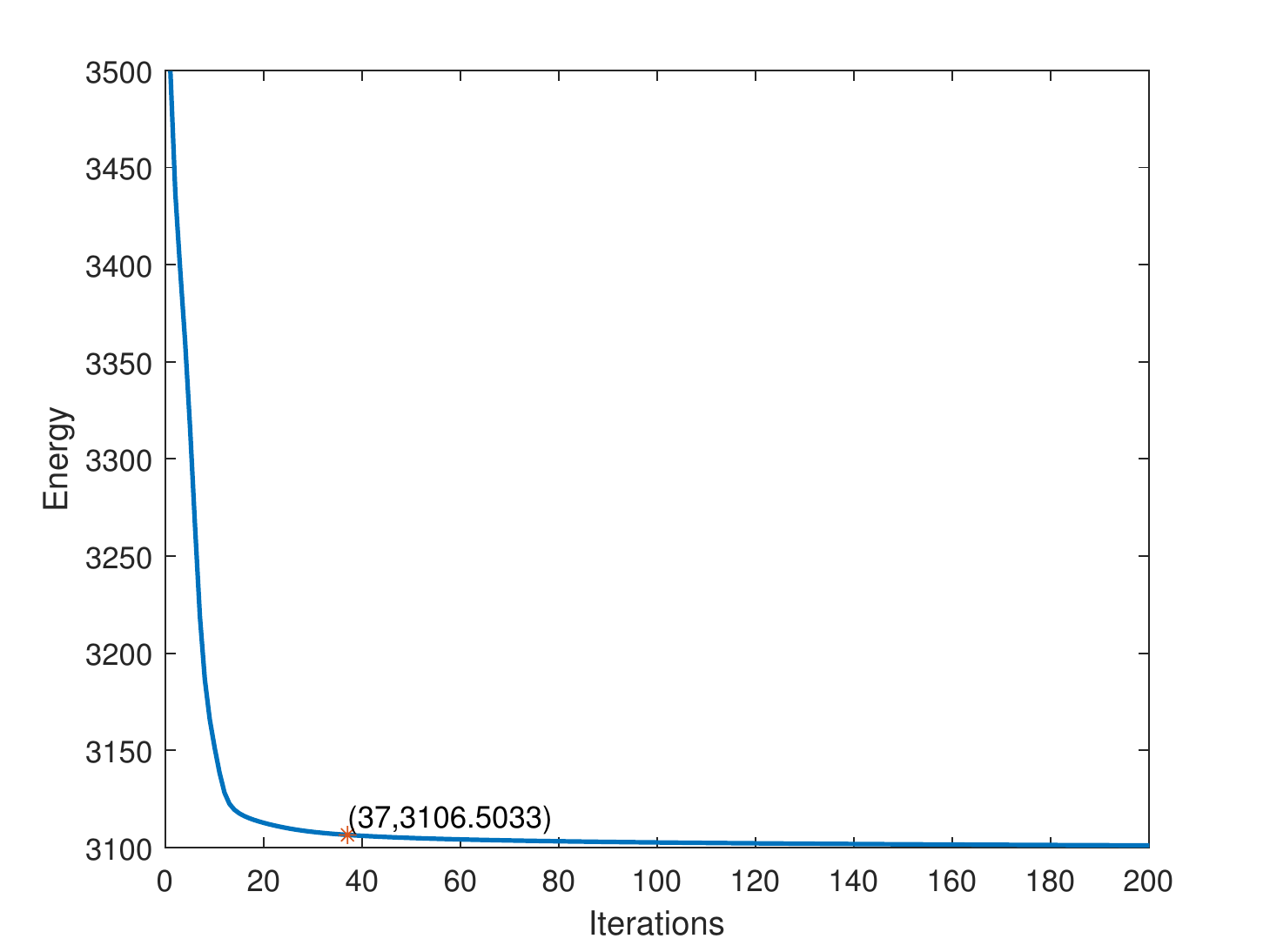}
		}
		\subfigure[\scriptsize Energy diagram of ETDRK2]{
			\includegraphics[width=0.45\linewidth]{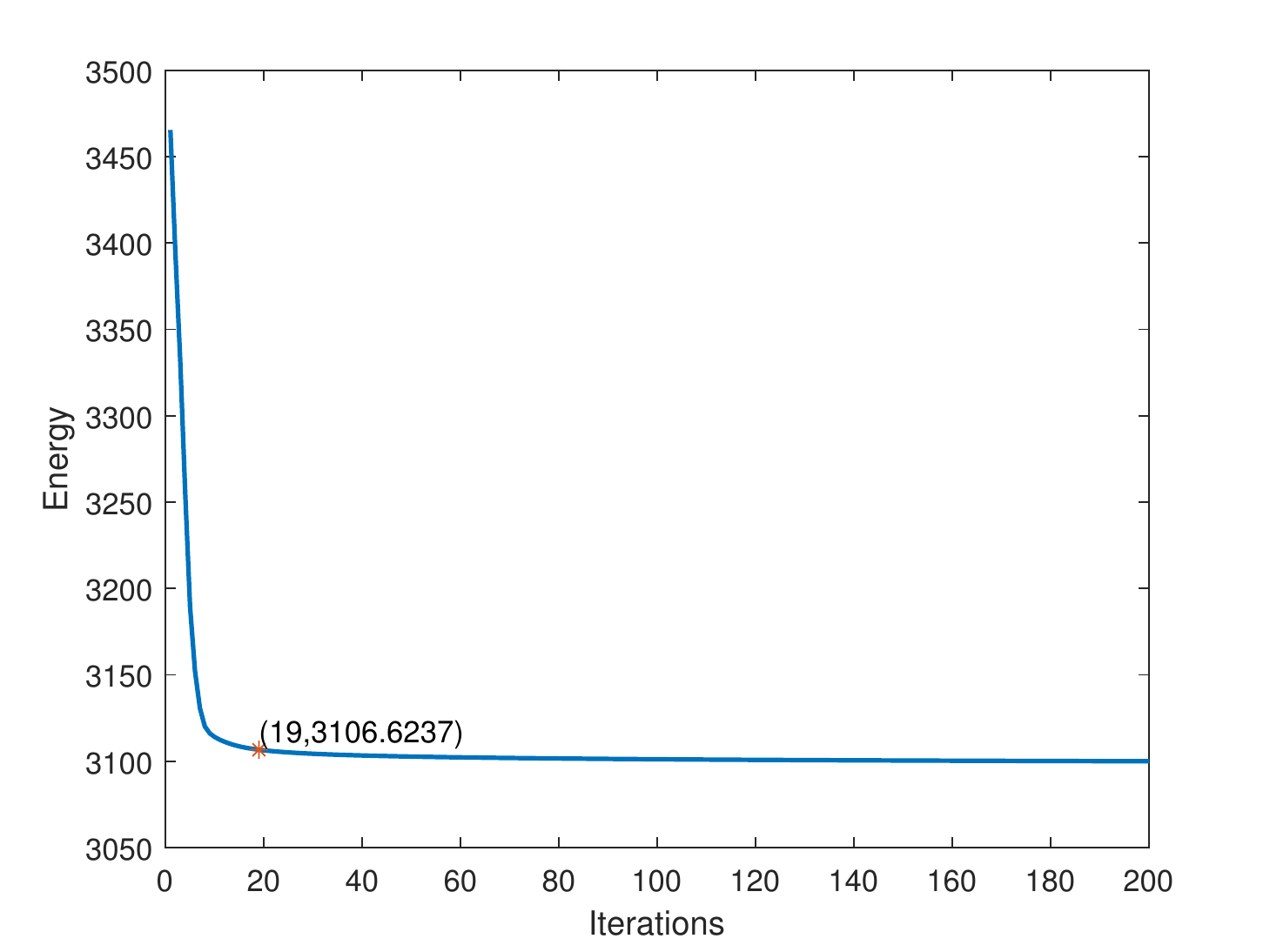}
		}
		\caption{\label{exETD1} Segmentation results and energy diagrams obtained by ETD1 and ETDRK2 schemes.}
\end{figure}

First, we are going to test the energy stability of the algorithms designed for our model. Fig. \ref{exETD1} shows segmentation results and the discrete energy evolution for a vessel image solved by the ACLBF model with ETD1 and ETDRK2 schemes. The segmentation result of the ETD1 scheme, which is almost the same as that of the ETDRK2 scheme, is represented in Fig. \ref{exETD1}(b). Two energy diagrams of ETD1 and ETDRK2 schemes are displayed in Fig. \ref{exETD1}(c) and Fig. \ref{exETD1}(d), respectively. The smallest iteration numbers to obtain segmentation results by ETD1 and ETDRK2 schemes are marked by red points on the energy curves. One can observe that the ACLBF model solved by the ETDRK2 scheme needs 19 iterations  {(CPU time: 0.222315s)} for the segmentation, which is much more efficient than ETD1 (37 iterations, {CPU time: 0.251346s}).
Moreover, the energy diagrams indicate that our algorithms are discrete energy stable. In this energy stability test, parameters for IGLIM are set as: $\lambda = 50,k_1 = k_2 = 0.01, M = 1$ and $S_n$ is chosen as the initial contour, and parameters for the ACLBF model are set as: $\lambda_1=\lambda_2=1, \mu =80,\sigma = 3,h=0.01,\Delta t =0.1, \epsilon =0.5, \epsilon_1 = 0.5$. At each iteration, $S$ is chosen to be $\frac{G}{2} + 1$ to guarantee the energy stability.

\subsection{Initialization comparison}
\begin{figure}[t!]
	\centering
	\subfigure{
	\begin{minipage}[b]{0.18\textwidth}
		\includegraphics[width=1.8cm]{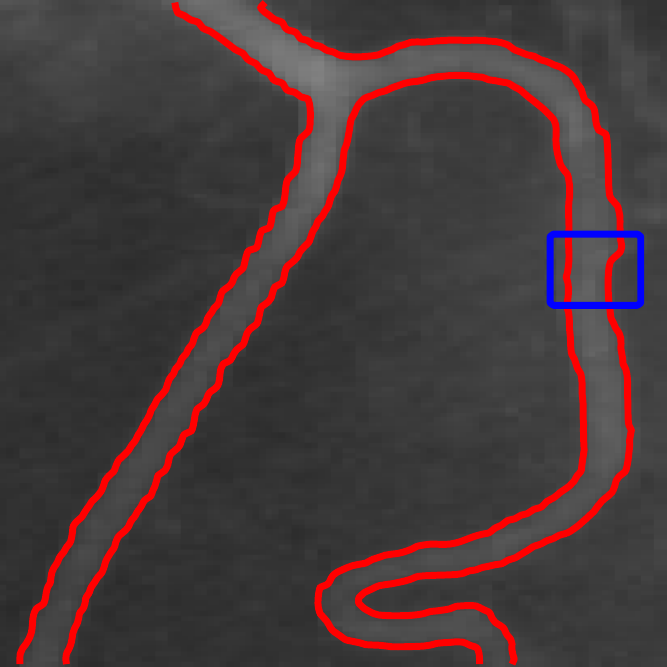}
		\includegraphics[width=1.8cm]{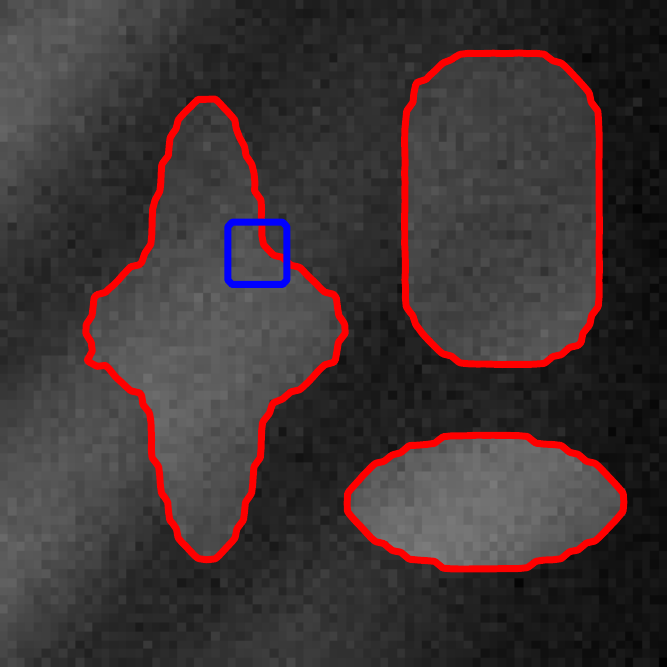}
		\includegraphics[width=1.8cm]{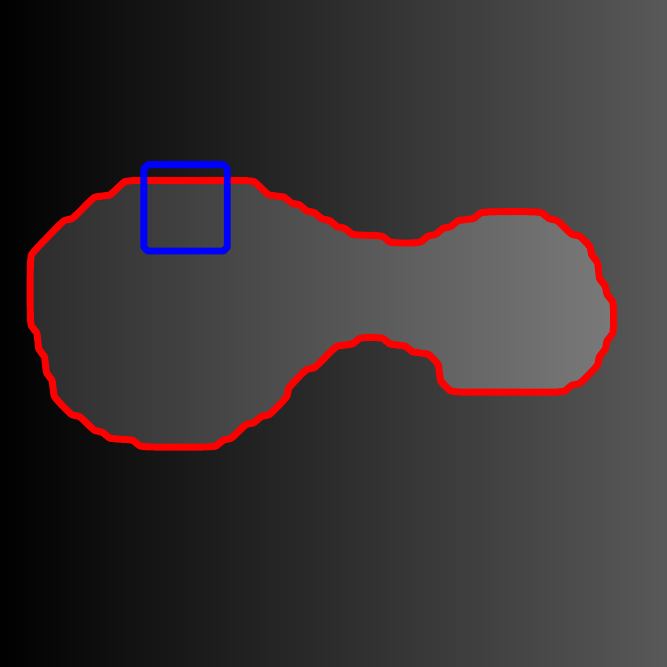}
		\centerline{(a)}
	\end{minipage}
	}
  \hspace{-10.3mm}
	\subfigure{
	\begin{minipage}[b]{0.18\textwidth}
	\includegraphics[width=1.8cm]{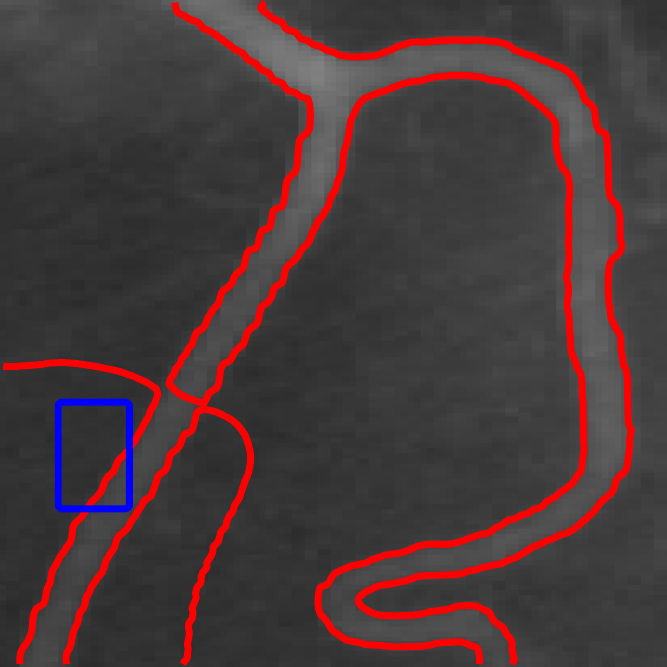}
	\includegraphics[width=1.8cm]{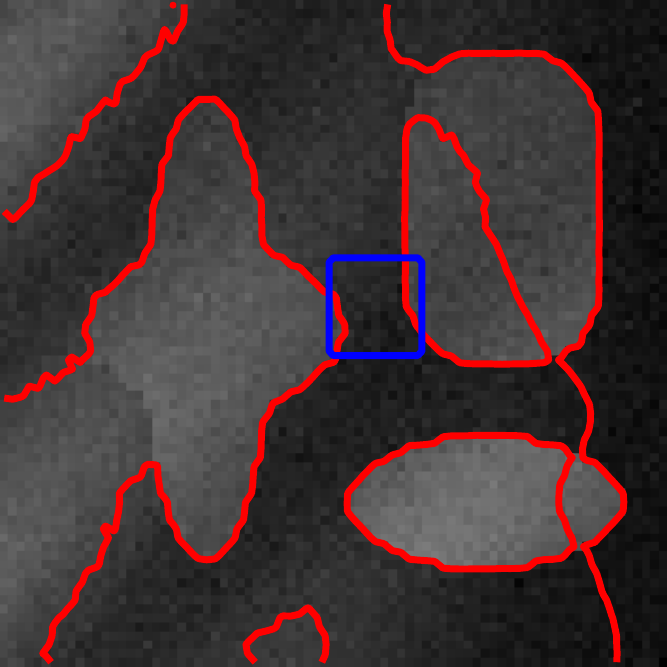}
	\includegraphics[width=1.8cm]{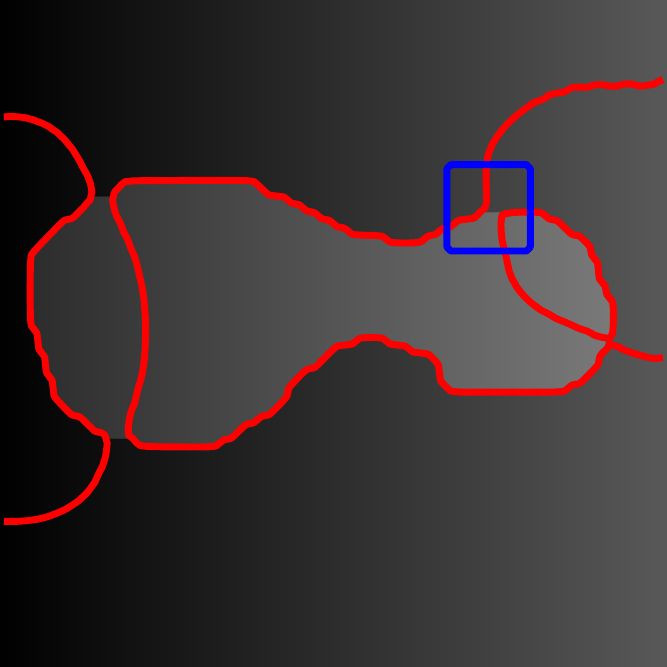}
	\centerline{(b)}
	\end{minipage}
	}
  \hspace{-10.2mm}
	\subfigure{
	\begin{minipage}[b]{0.18\textwidth}
		\includegraphics[width=1.8cm]{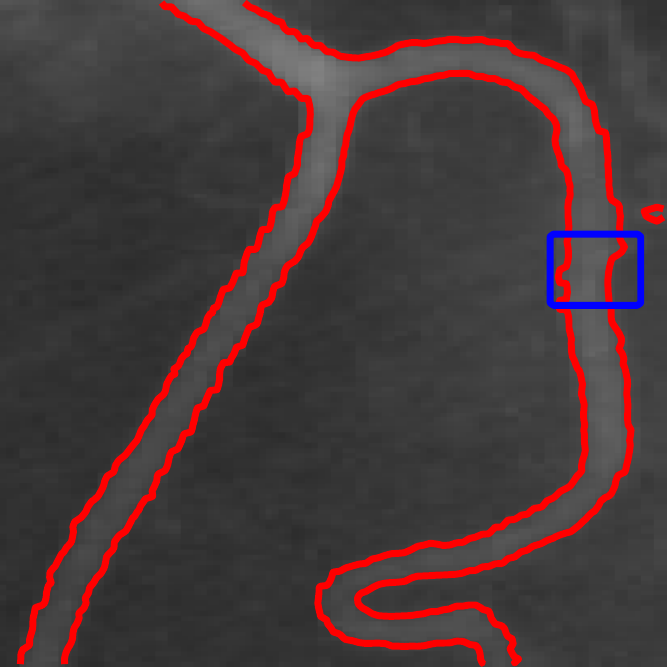}
		\includegraphics[width=1.8cm]{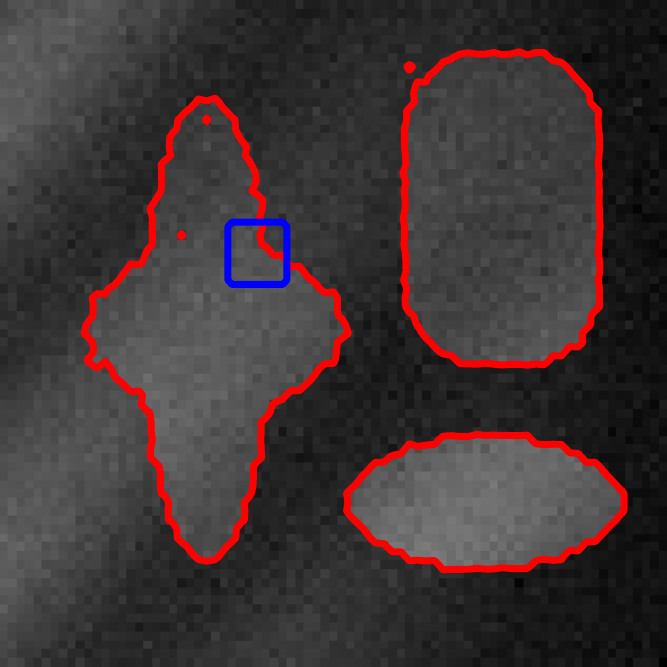}
		\includegraphics[width=1.8cm]{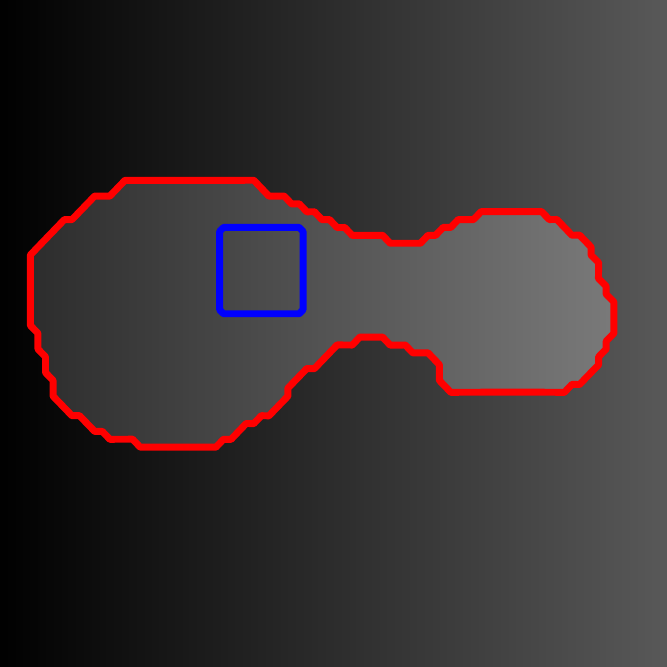}
		\centerline{(c)}
  \end{minipage}
	}
  \hspace{-10.3mm}
	\subfigure{
	\begin{minipage}[b]{0.18\textwidth}
	\includegraphics[width=1.8cm]{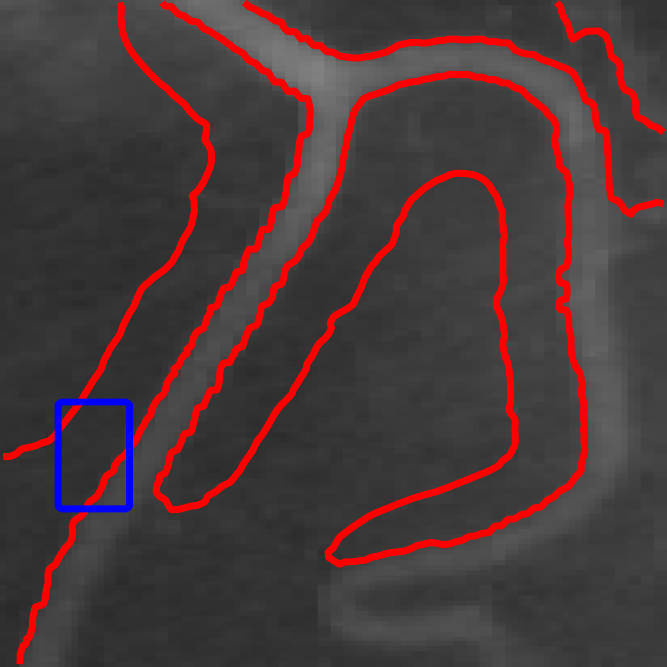}
	\includegraphics[width=1.8cm]{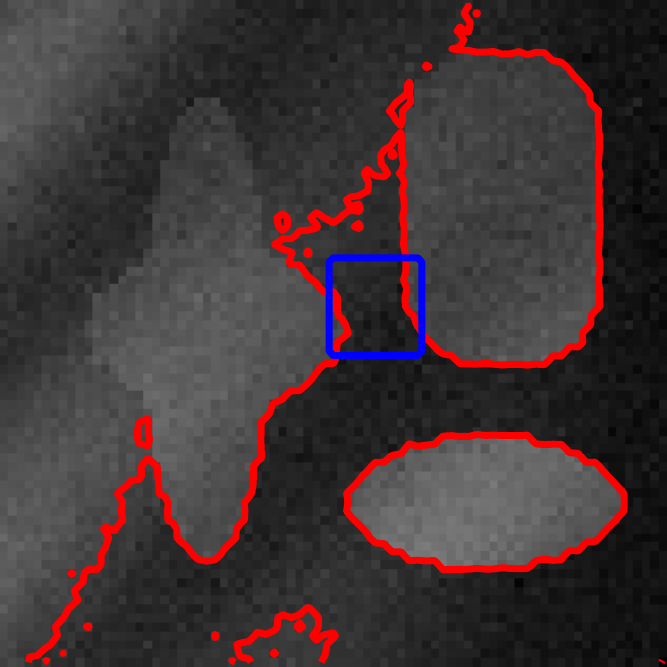}
	\includegraphics[width=1.8cm]{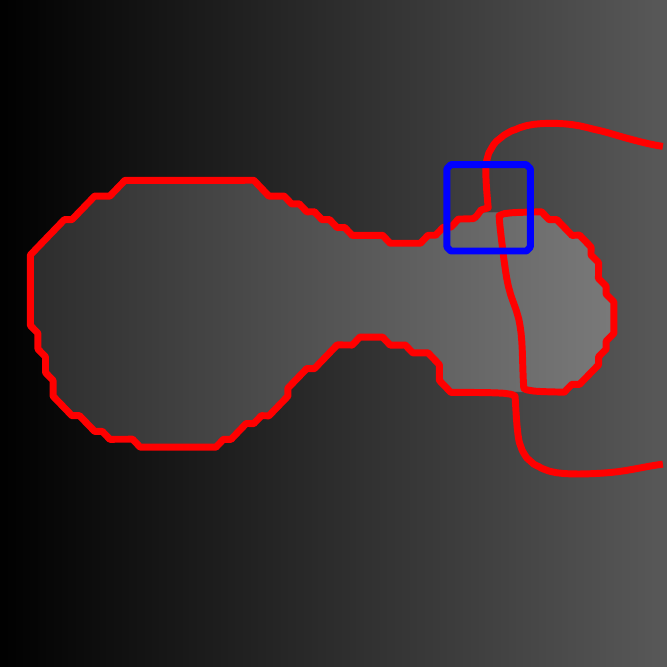}
	\centerline{(d)}
	\end{minipage}
	}
  \hspace{-10.3mm}
	\subfigure{
	\begin{minipage}[b]{0.18\textwidth}
	\includegraphics[width=1.8cm]{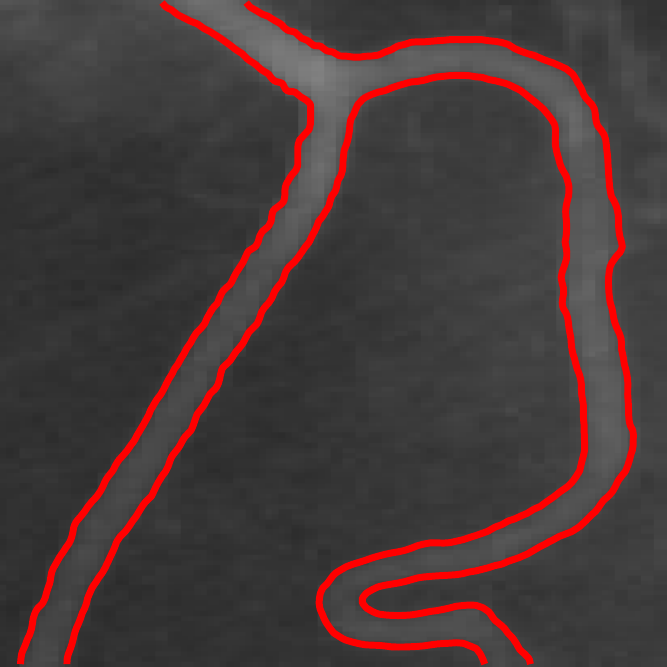}
	\includegraphics[width=1.8cm]{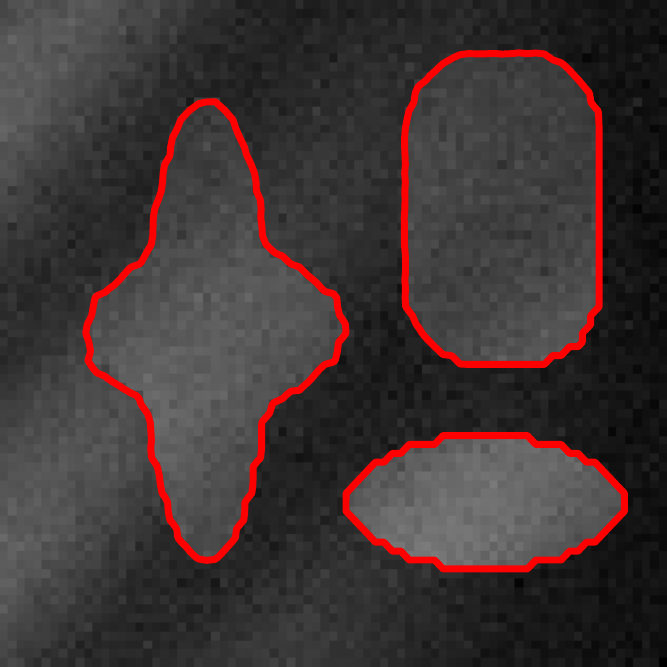}
	\includegraphics[width=1.8cm]{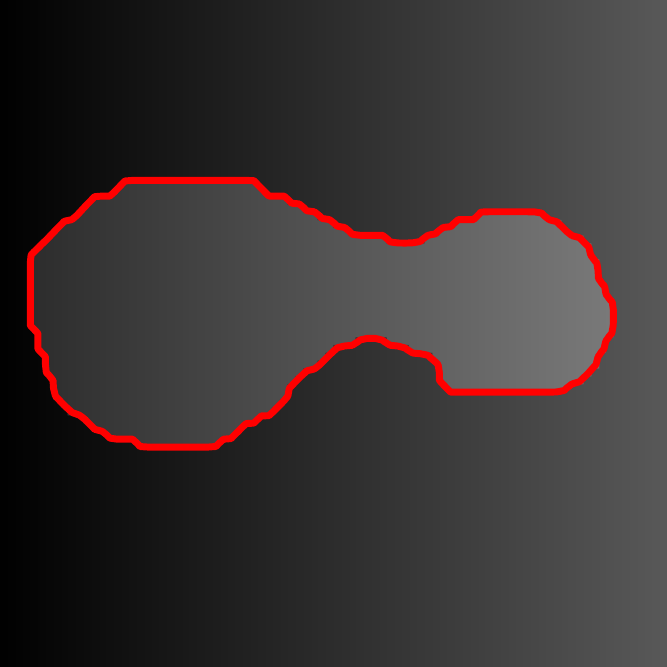}
	\centerline{(e)}
	\end{minipage}
	}
  \hspace{-10.2mm}
	\subfigure{
	\begin{minipage}[b]{0.18\textwidth}
	\includegraphics[width=1.8cm]{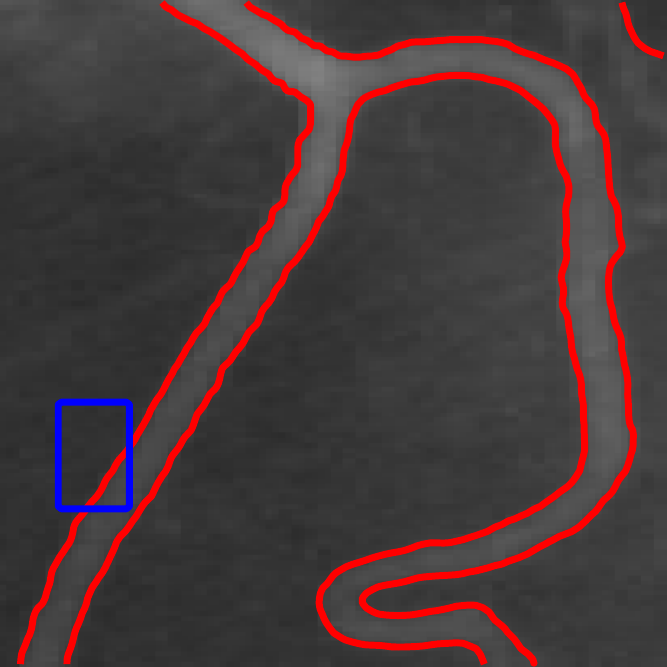}
	\includegraphics[width=1.8cm]{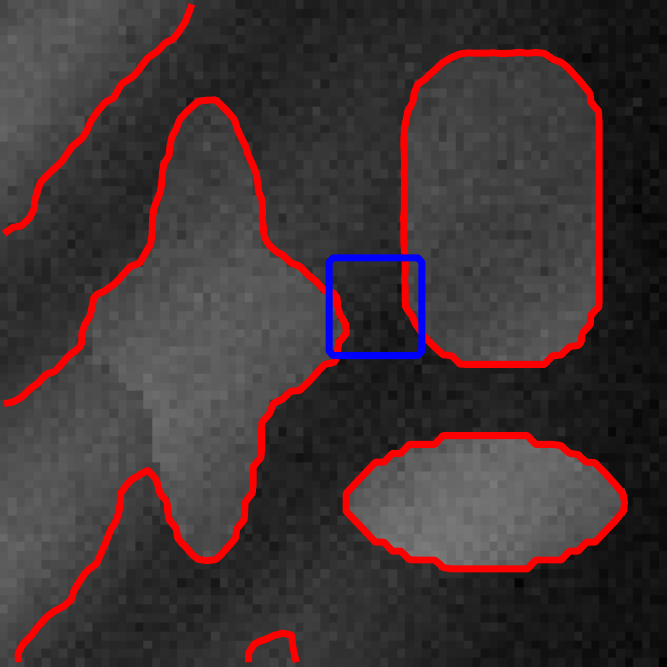}
	\includegraphics[width=1.8cm]{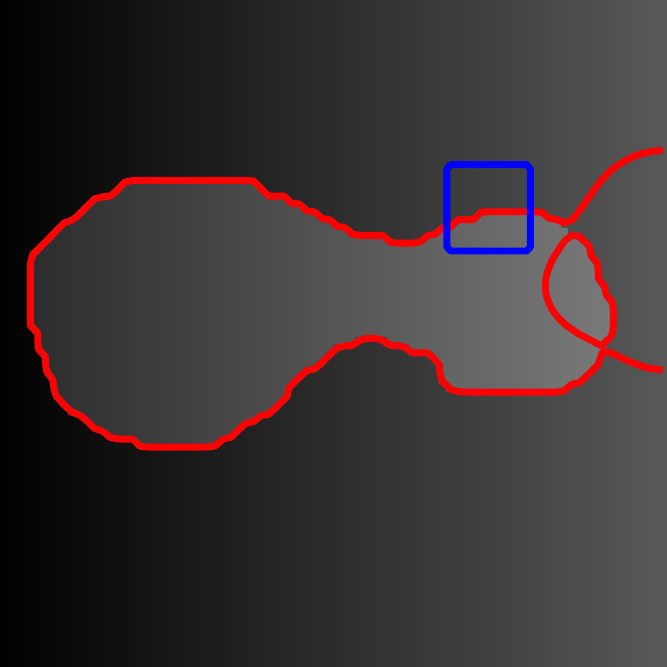}
	\centerline{(f)}
	\end{minipage}
	}
	 \caption{ (a) and (b) the results of LBF model solved by level set method with proper and improper initial contours, respectively; (c) and (d) the segmentation results of LBF model solved by ICTM with proper and improper initial contours, respectively; (e) the result of ACLBF model with IGLIM and (f) the results of ACLBF model with improper initial contours. Initial contours are represented by blue outlines, and the results are highlighted by red ones.} \label{fig1}
\end{figure}

The LBF model can effectively segment inhomogeneous images, but the segmentation result may be seriously affected by initialization.
In this part, we segment three images by solving the LBF model with the level set method, ICTM, and our phase-field approach. Initial contours of the first two methods are given by selecting different parts of the images. The corresponding results are shown in (a)-(d) of Fig. \ref{fig1}. Then we compare segmentation results of the ACLBF model with initialization by IGLIM and by selecting a part of the images, which are exhibited in (e) and (f) of Fig. \ref{fig1}. It can be seen that all these methods are sensitive to initialization when solving the LBF model for segmentation and solving the ACLBF model with initial contours from IGLIM gives satisfactory segmentation results. Table \ref{table1} shows the iteration numbers and the CPU time for the three methods successfully segmenting the three images. One can see that iteration numbers of the ACLBF model solved by the ETDRK2 scheme are smaller than those of the other two methods, while ICTM costs less CPU time if appropriate initial contour is given.

\begin{table}[t!]
\centering
	\begin{tabular}{lllllll}
		\hline
		\multirow{2}{*}{Images} &
		\multicolumn{2}{l}{Level set} & \multicolumn{2}{l}{ICTM} & \multicolumn{2}{l}{ACLBF-ETDRK2}\cr\cmidrule(lr){2-3}\cmidrule(lr){4-5}\cmidrule(lr){6-7}
		&Ite.&Time(s)&Ite.&Time(s)&Ite.&Time(s)\cr
		\hline
		Row1 &124 &0.378146 &51  &0.087089 &35 & 0.314895 \cr
		Row2 &69  &0.164980 &43  &0.053831 &18  & 0.145300 \cr
		Row3 &29  &0.106062 &12  &0.035150 &8  & 0.112060\cr\hline
	\end{tabular}
	\label{table1}\caption{Comparison of the iteration number and running time in Fig. \ref{fig1}.}
\end{table}

\subsection{Experiments on images with intensity inhomogeneity}
Next, we solve the ACLBF model for segmentations of various images, including images with intensity inhomogeneity. At first, two synthetic images with severe intensity inhomogeneity are used to test the performance of the ACLBF model with IGLIM. Next, we apply our model to two medical images, magnetic resonance images of a human brain and an angiogram of a blood vessel. Then we verify the capability of the ACLBF model and IGLIM on two real natural images. The parameters settings are displayed in Table \ref{parameters}.

There are some strategies for the parameters setting. For IGLIM, $\lambda$ is always set as 50. Generally, $k_1$ and $k_2$ are both set as 0.01. But for images with inconspicuous edges, $k_1,k_2$ can be set smaller to capture more edges. Conversely, they can be tuned to be larger numbers to avoid excessive detection if the edges are apparent. The value of $M$ mainly depends on the noise level. One can choose a relatively large $M$ for images with strong noise. According to Remark \ref{bd_judge}, $S_p$ is generally chosen as the initial contour if the intensity of the object is smaller than that of the background. Otherwise, $S_n$ is taken as the initial contour. 

 {For ACLBF, $h,\Delta t$, and $\epsilon$ are parameters in the discrete Allen-Cahn term and we fix $h=0.01, \Delta t = 0.1$ and $\epsilon = 0.5$ in the experiments. The spatial step $h$ controls the smoothness of the segmentation result. Smaller $h$ can help to recognize smooth boundaries and remove noise. By contrast, larger $h$ can obtain a more delicate segmentation for some unsmooth boundaries. 
On the other hand, $\lambda_i(i = 1,2), \epsilon_1$ and $\sigma$ are parameters in the LBF term.  $\lambda_i(i = 1,2)$ are corresponding to the coefficients in front of the evolving forces of internal and external regions. In most cases, we set $\lambda_1 = \lambda_2 = 1$. $\sigma$ is the standard deviation in the Gaussian kernel. As shown in \cite{li2008minimization}, a more accurate
segmentation result can be obtained for images with inhomogeneity if $\sigma$ is small. Furthermore, a reasonably large $\sigma$ is more suitable for many real-world images of which the intensity inhomogeneity is not so severe. It is also pinpointed in \cite{li2008minimization} that a large $\epsilon_1$ can fasten curve evolution due to the fast emergence of new contours at strong edges while a small $\epsilon_1$ will have higher accuracy in
the final contour location. In \cite{li2008minimization}, $\epsilon_1$ is always chosen as 1 while in this paper it is chosen as 0.5 or 1.}

\begin{figure}[t!]
	\centering
	\subfigure{
		\includegraphics[width=3cm]{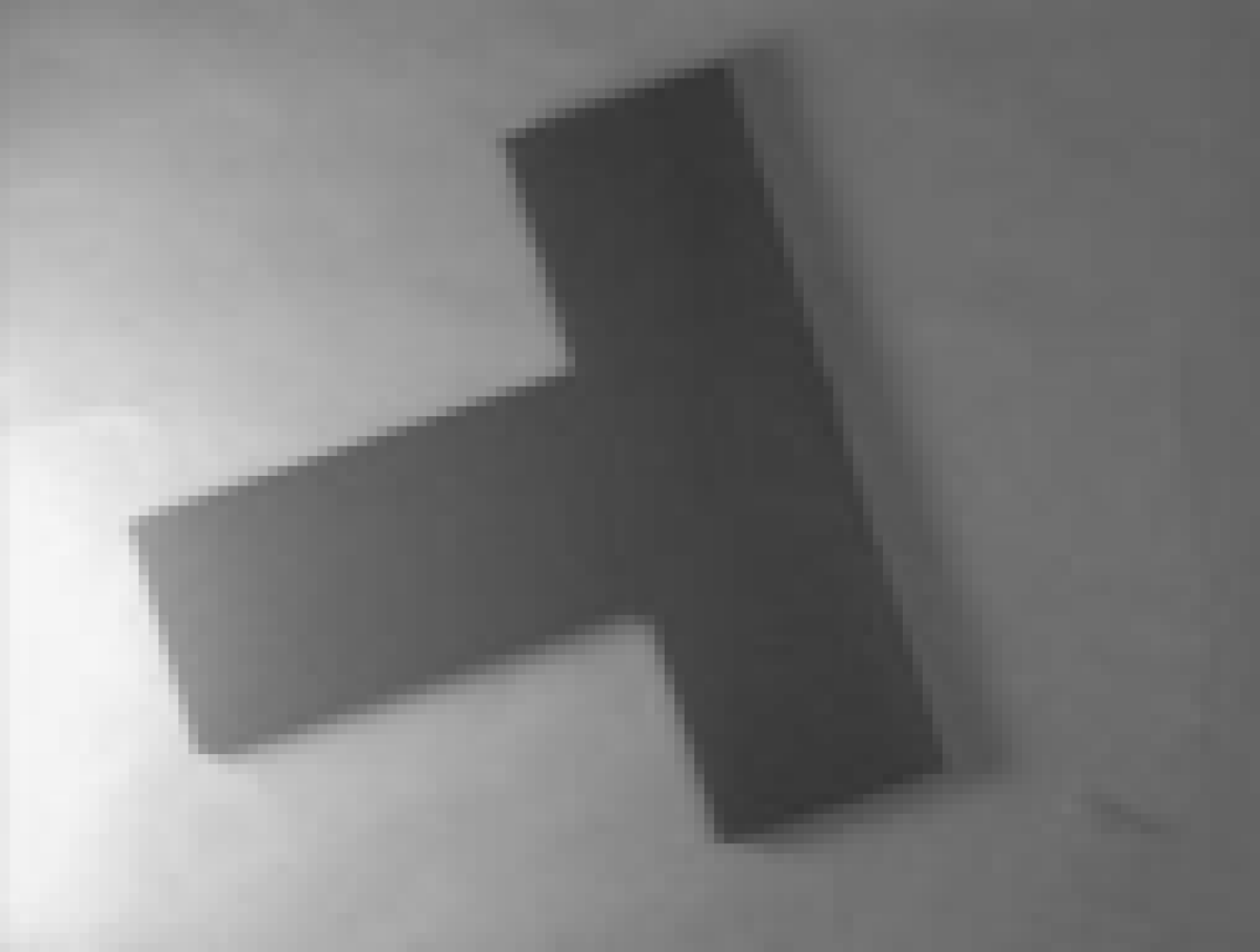}
		\includegraphics[width=3cm]{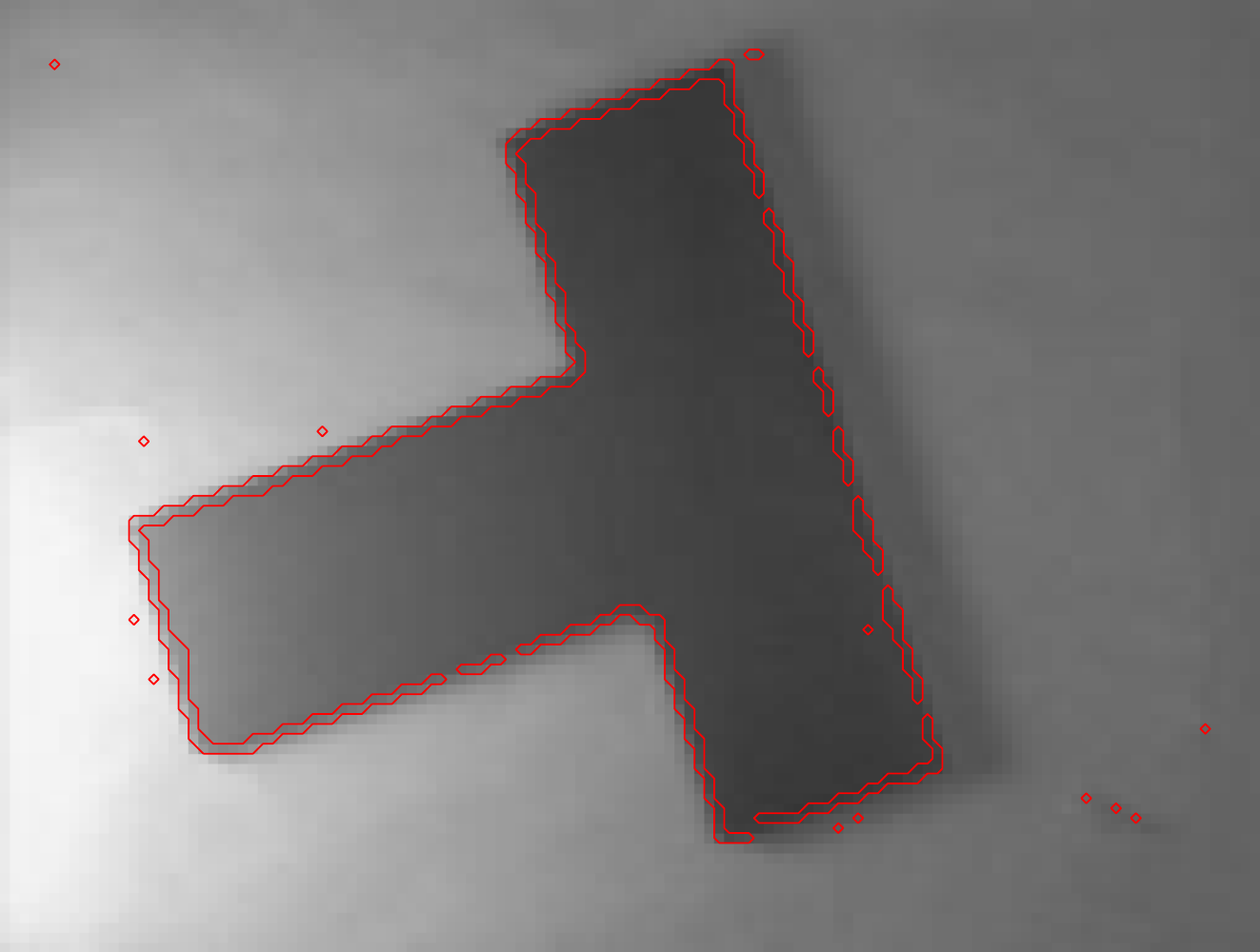}
		\includegraphics[width=3cm]{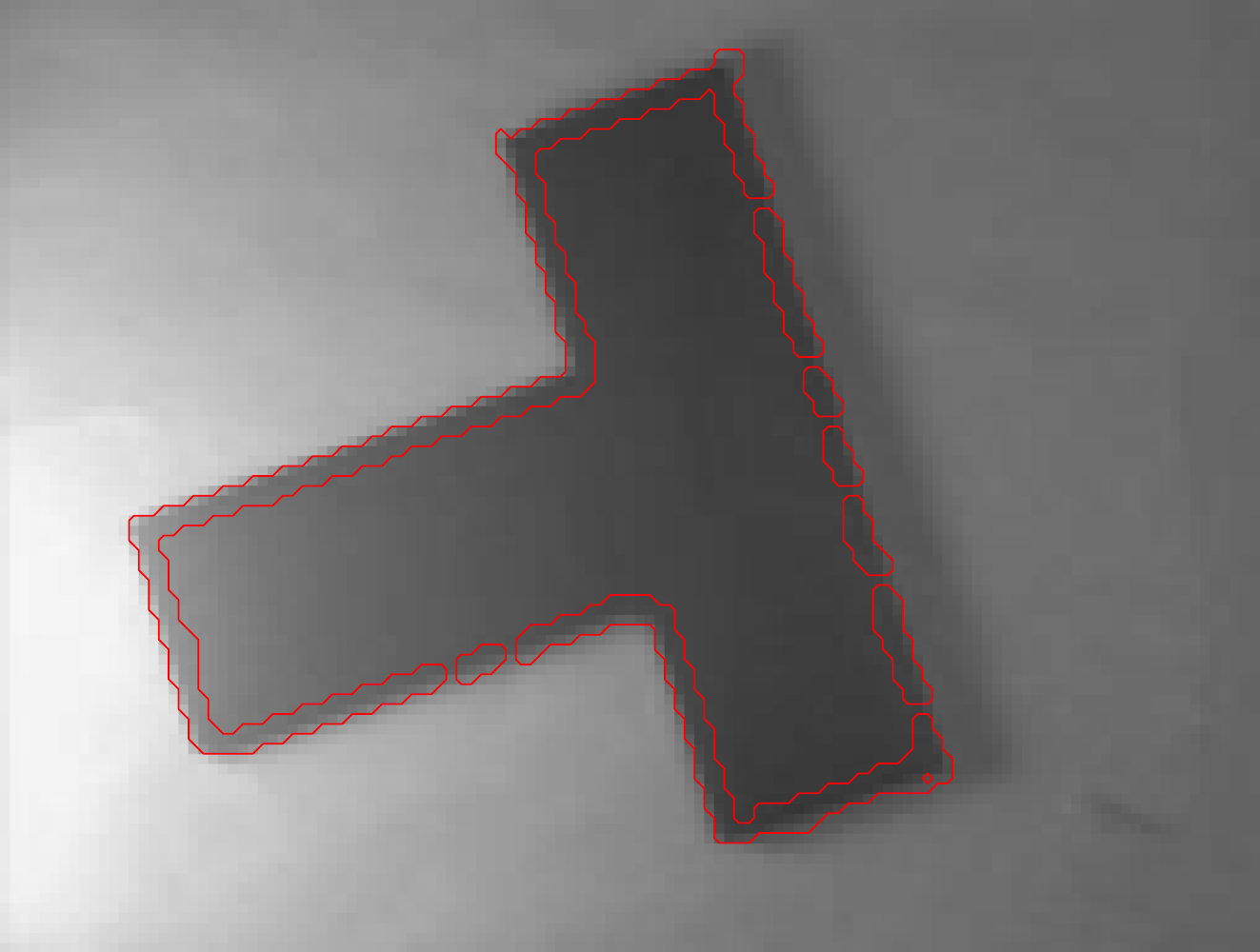}
		\includegraphics[width=3cm]{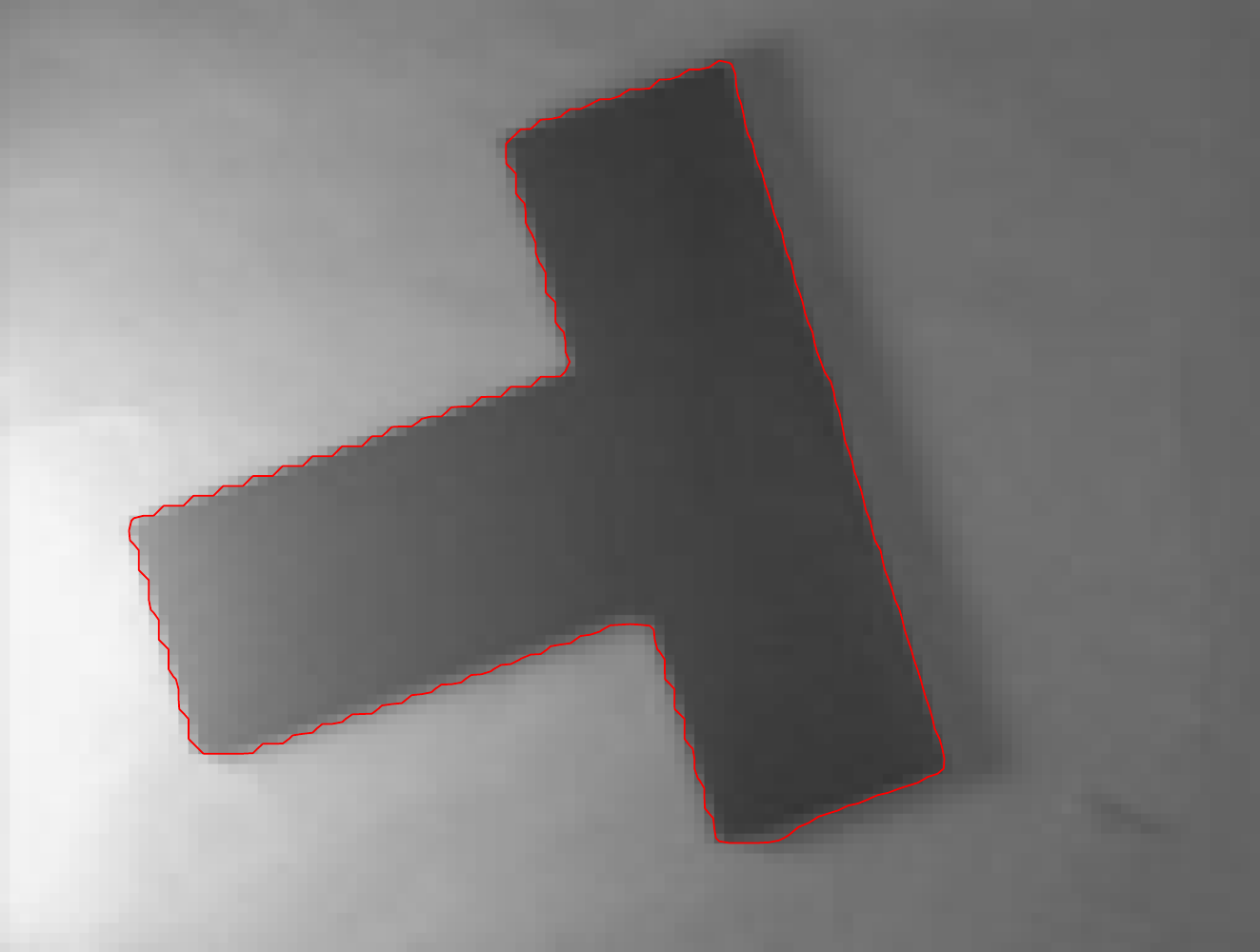}
	}
	\subfigure{
		\includegraphics[width=3cm]{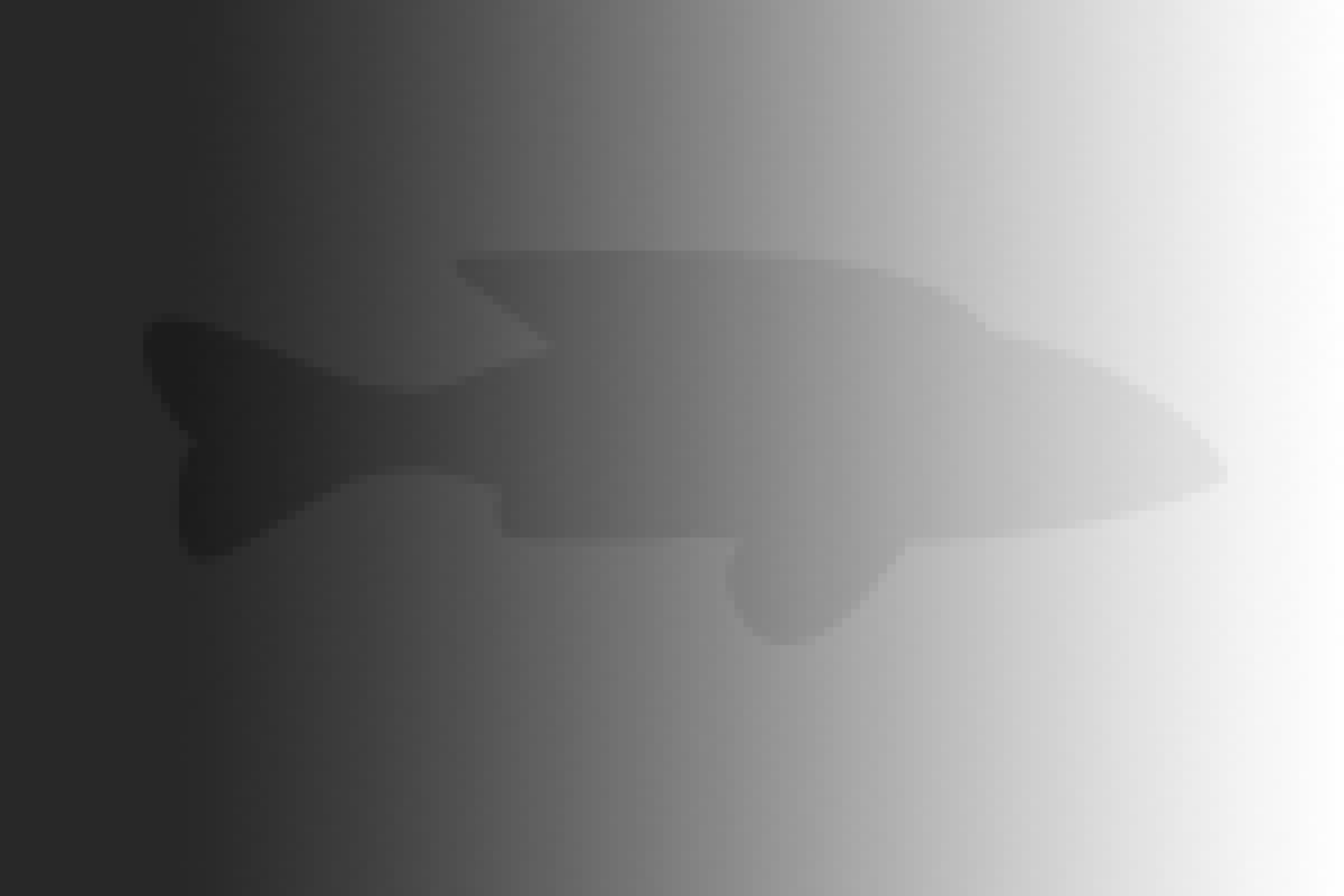}
		\includegraphics[width=3cm]{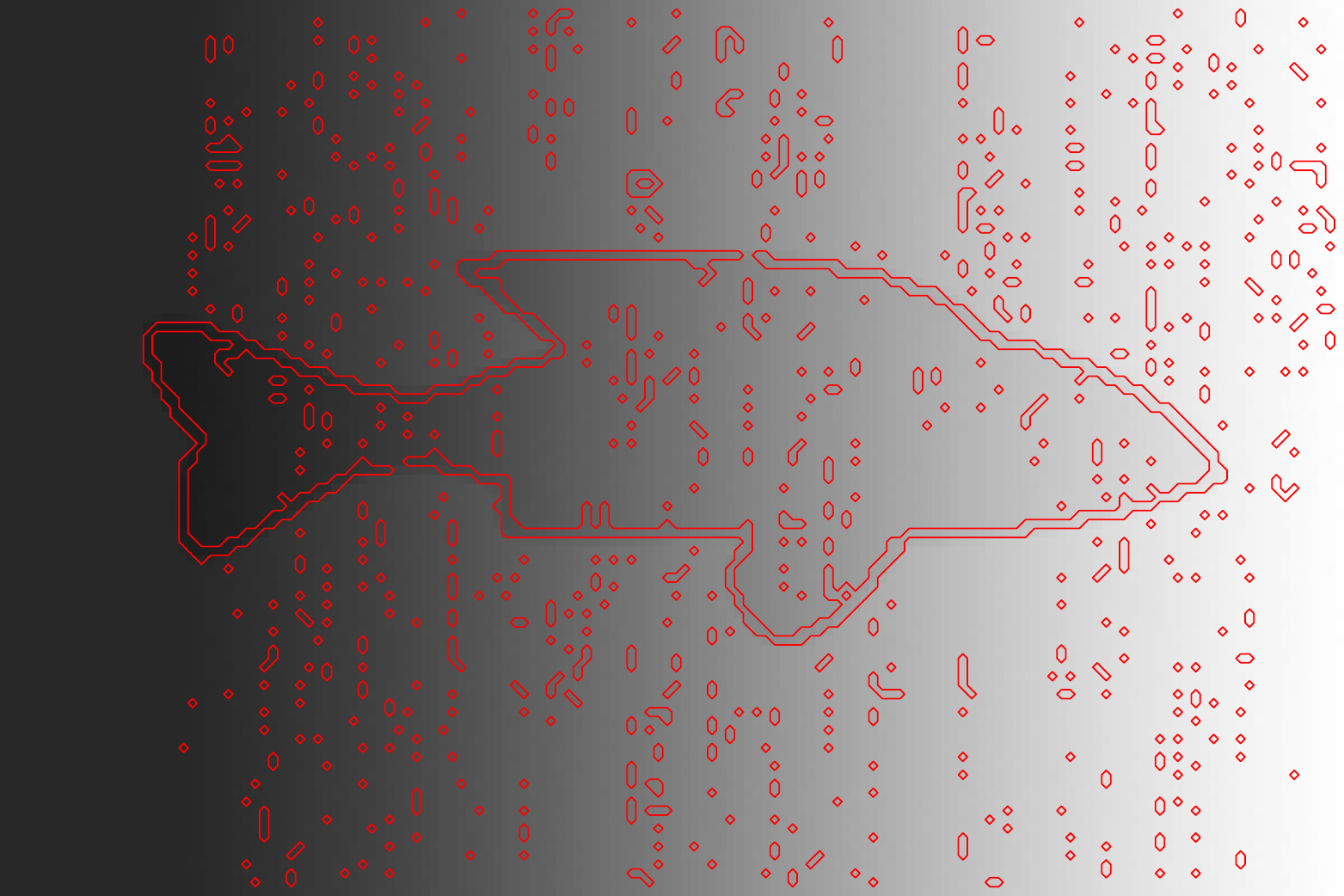}
		\includegraphics[width=3cm]{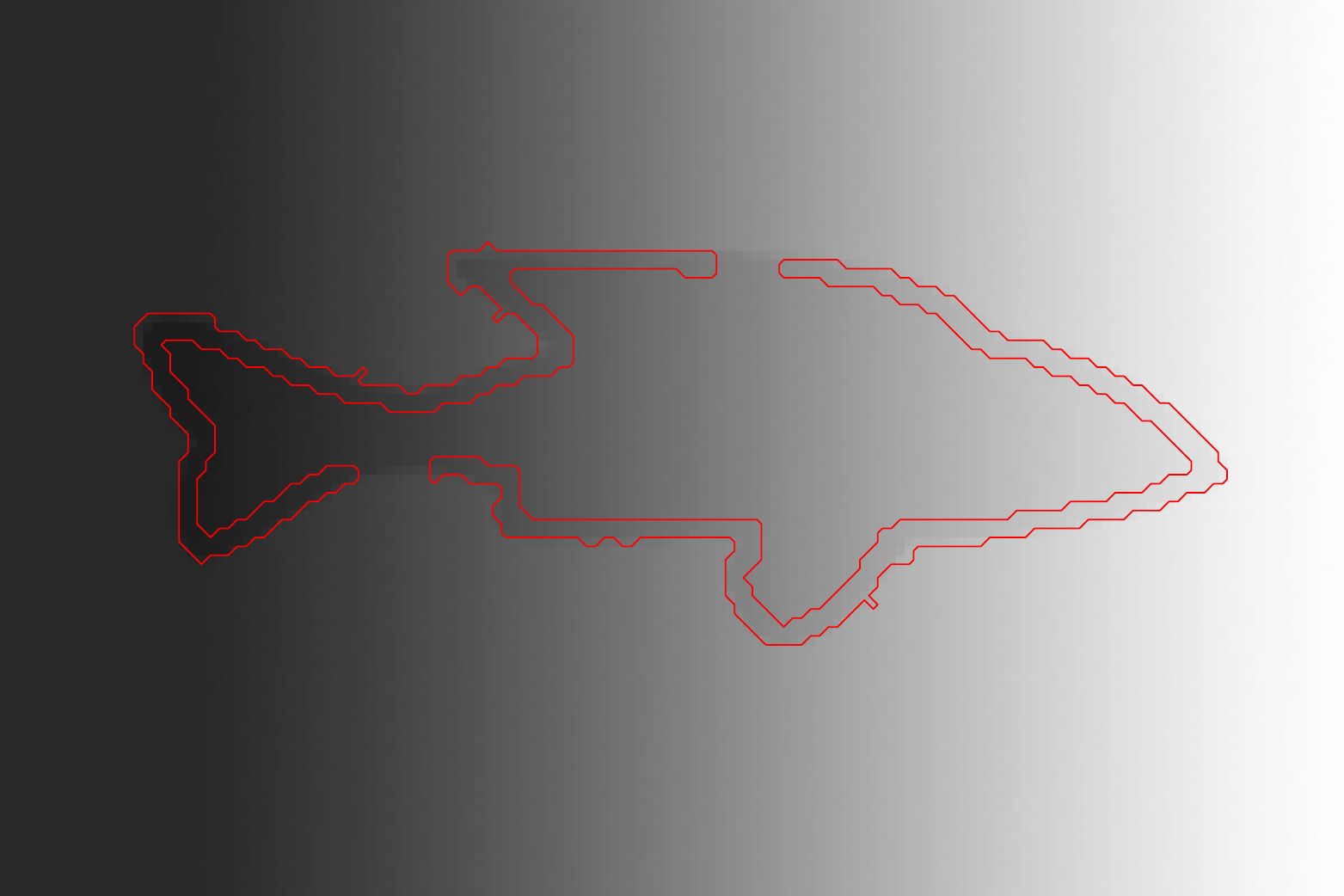}
		\includegraphics[width=3cm]{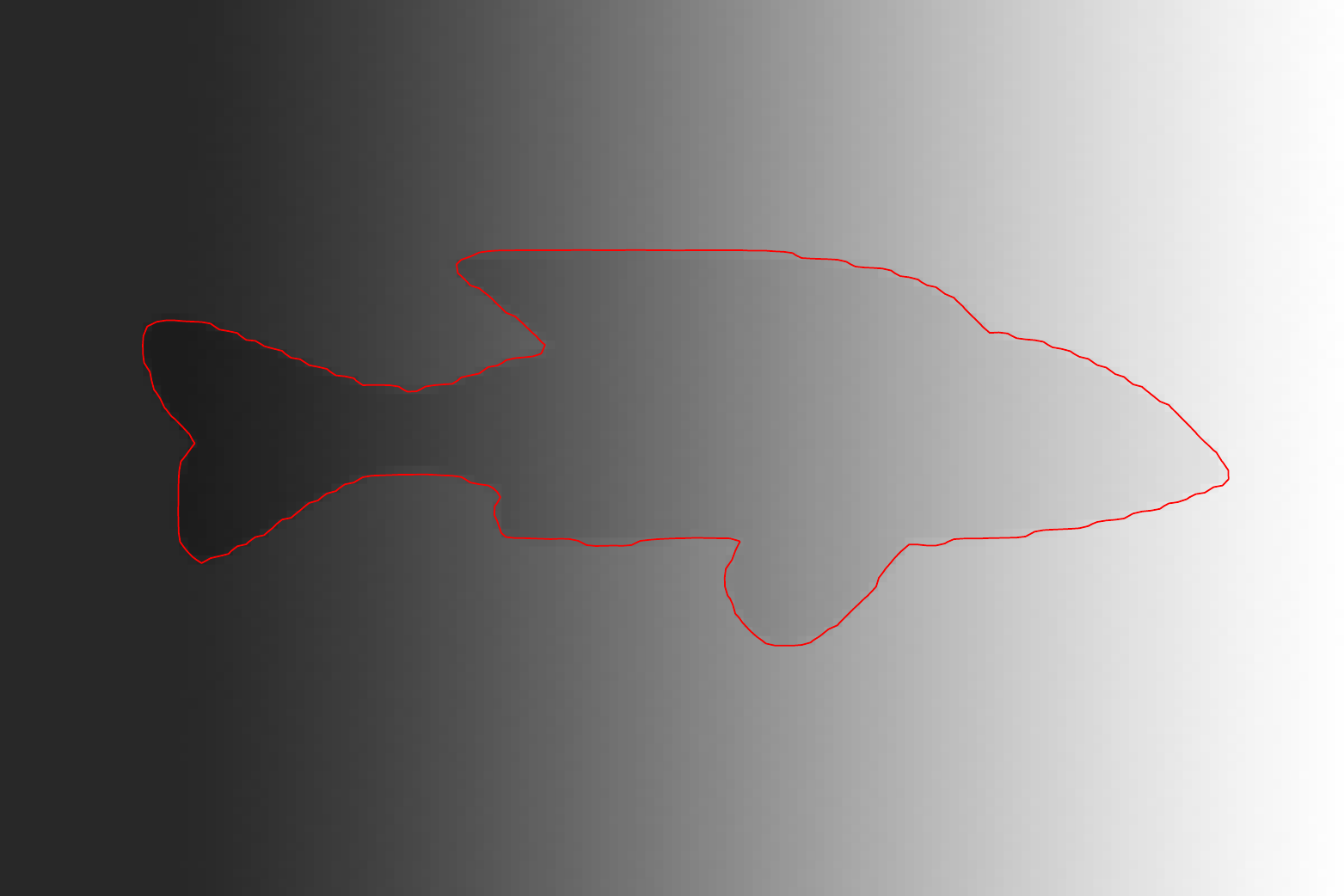}
	}
	\subfigure{
		\includegraphics[width=3cm]{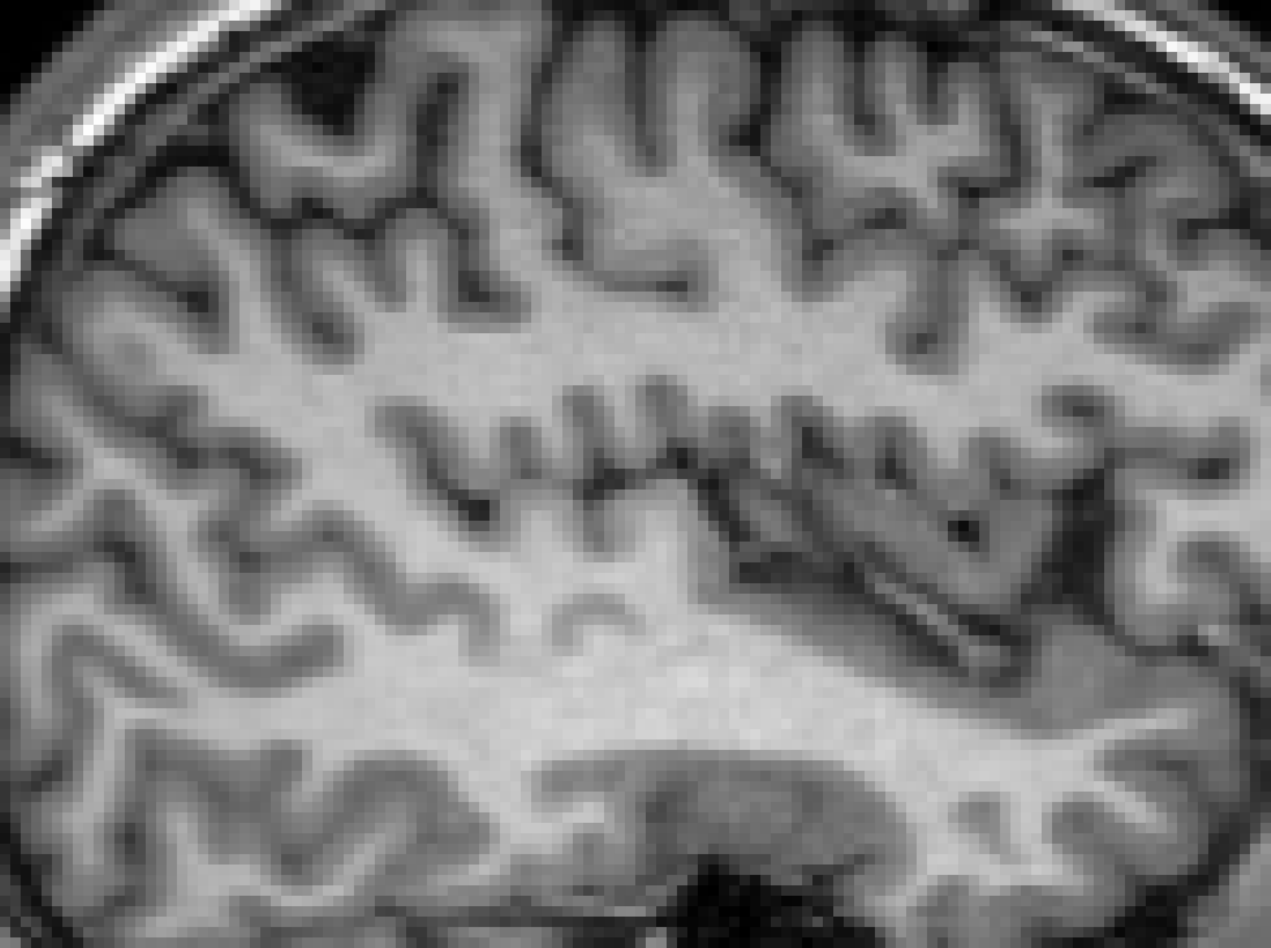}
		\includegraphics[width=3cm]{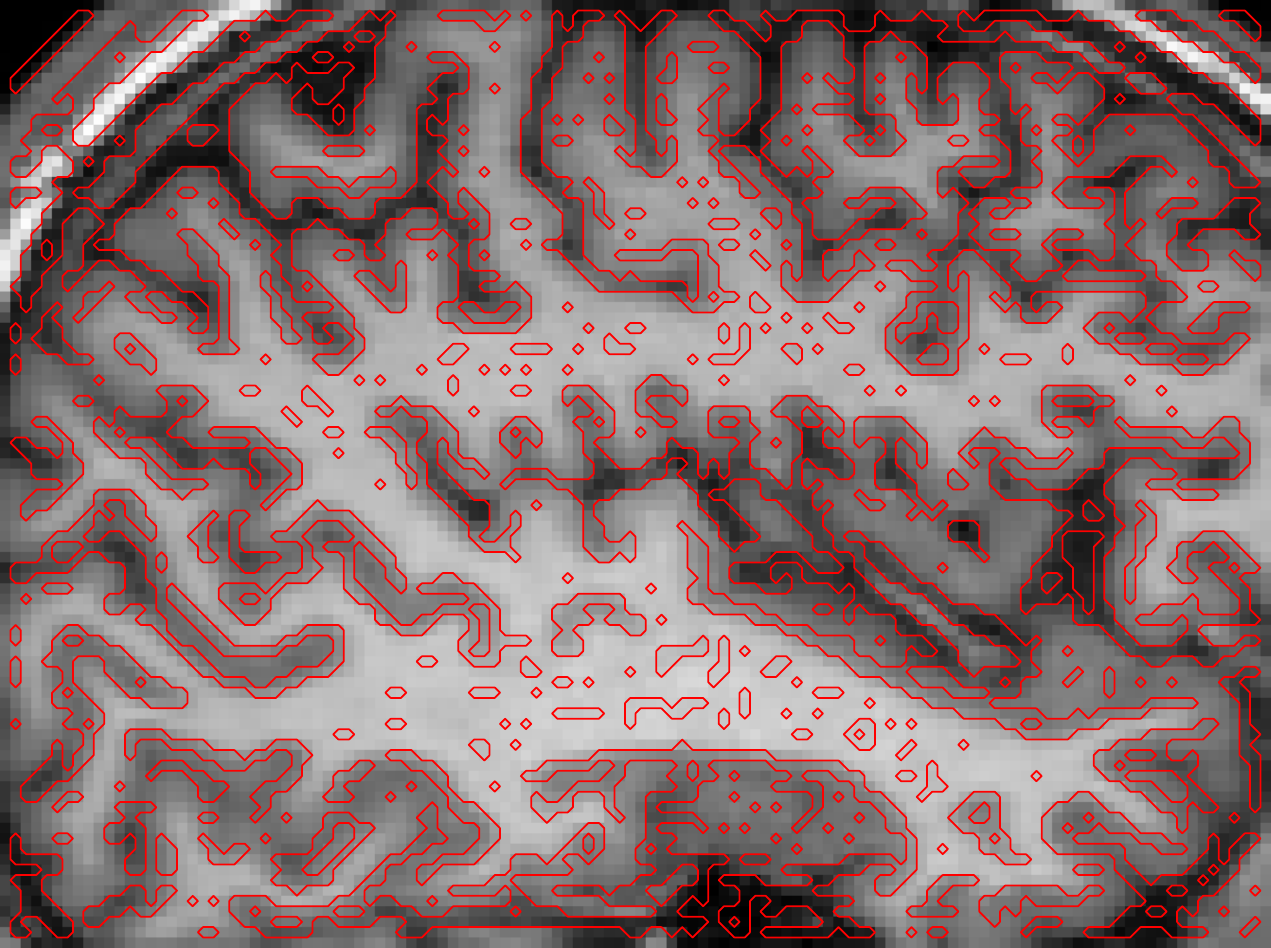}
		\includegraphics[width=3cm]{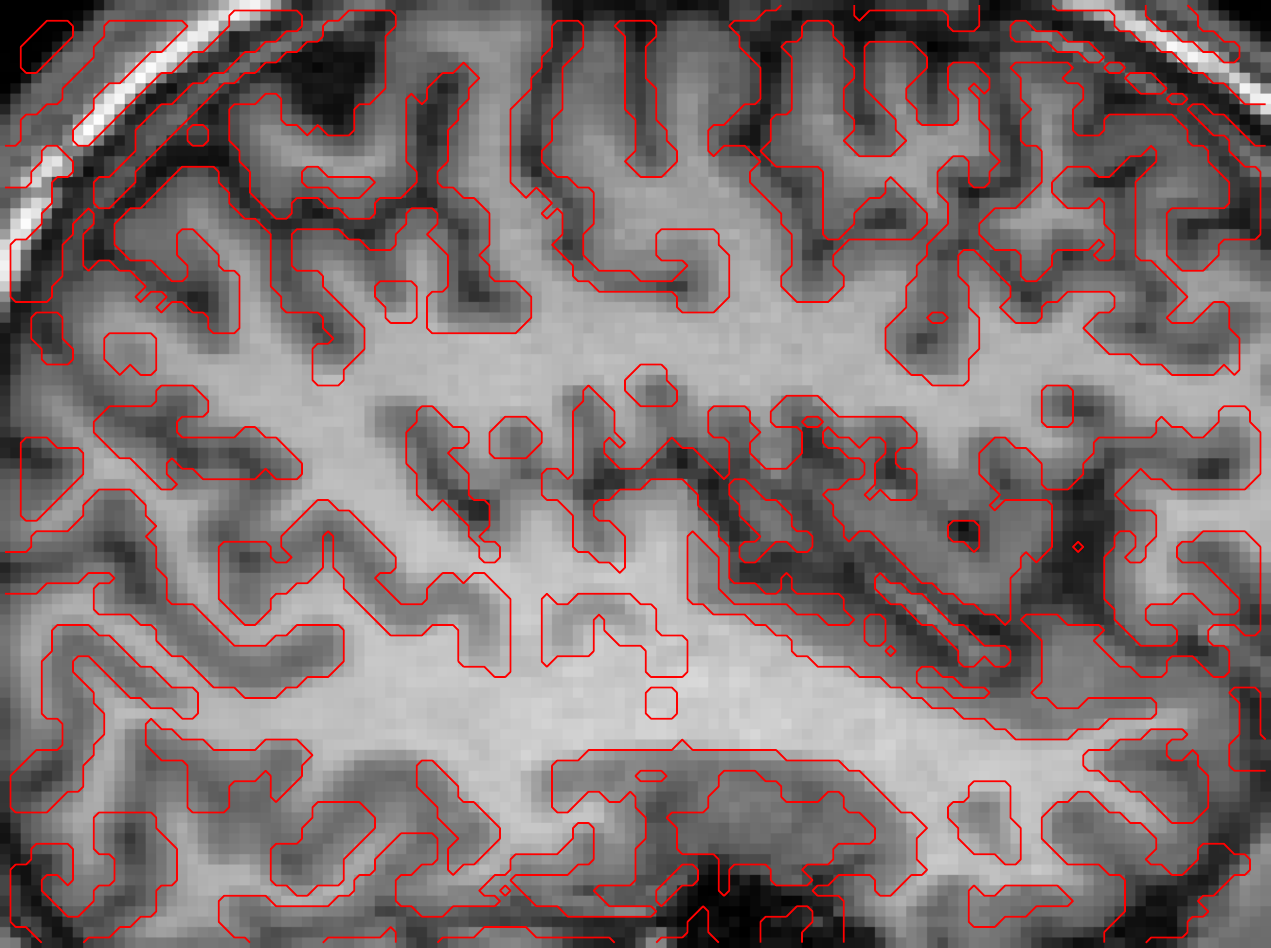}
		\includegraphics[width=3cm]{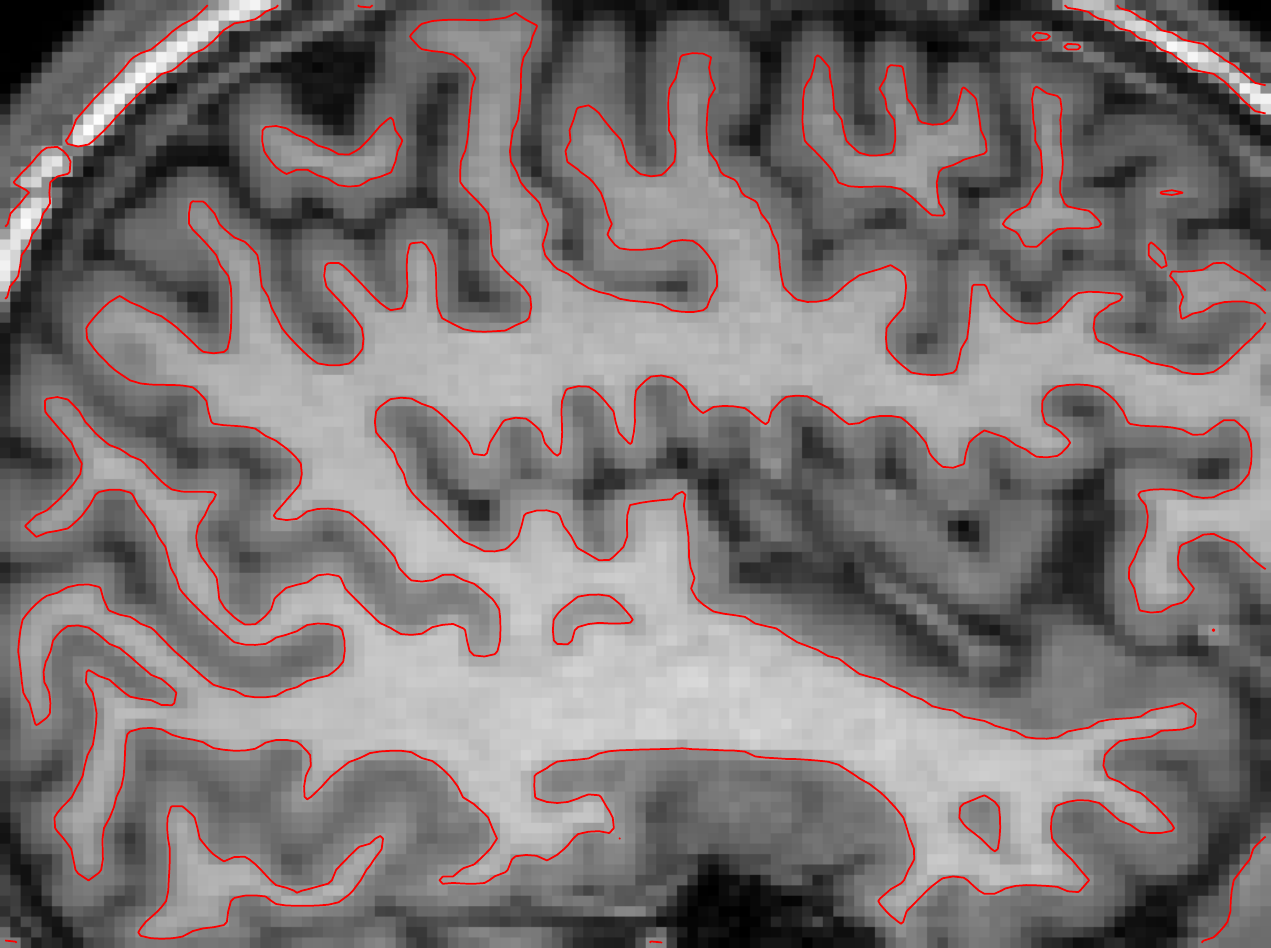}
	}
	\subfigure{
		\includegraphics[width=3cm]{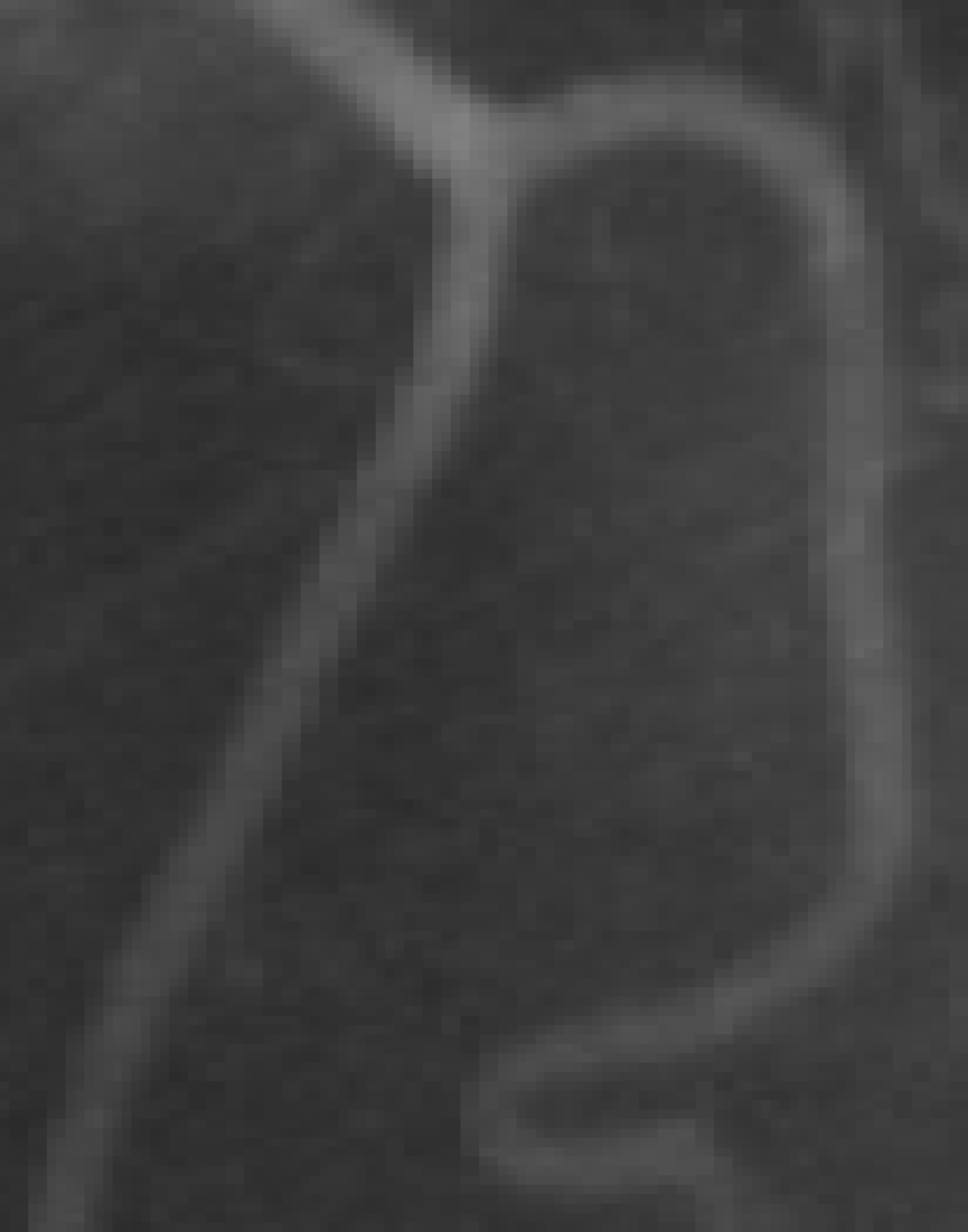}
		\includegraphics[width=3cm]{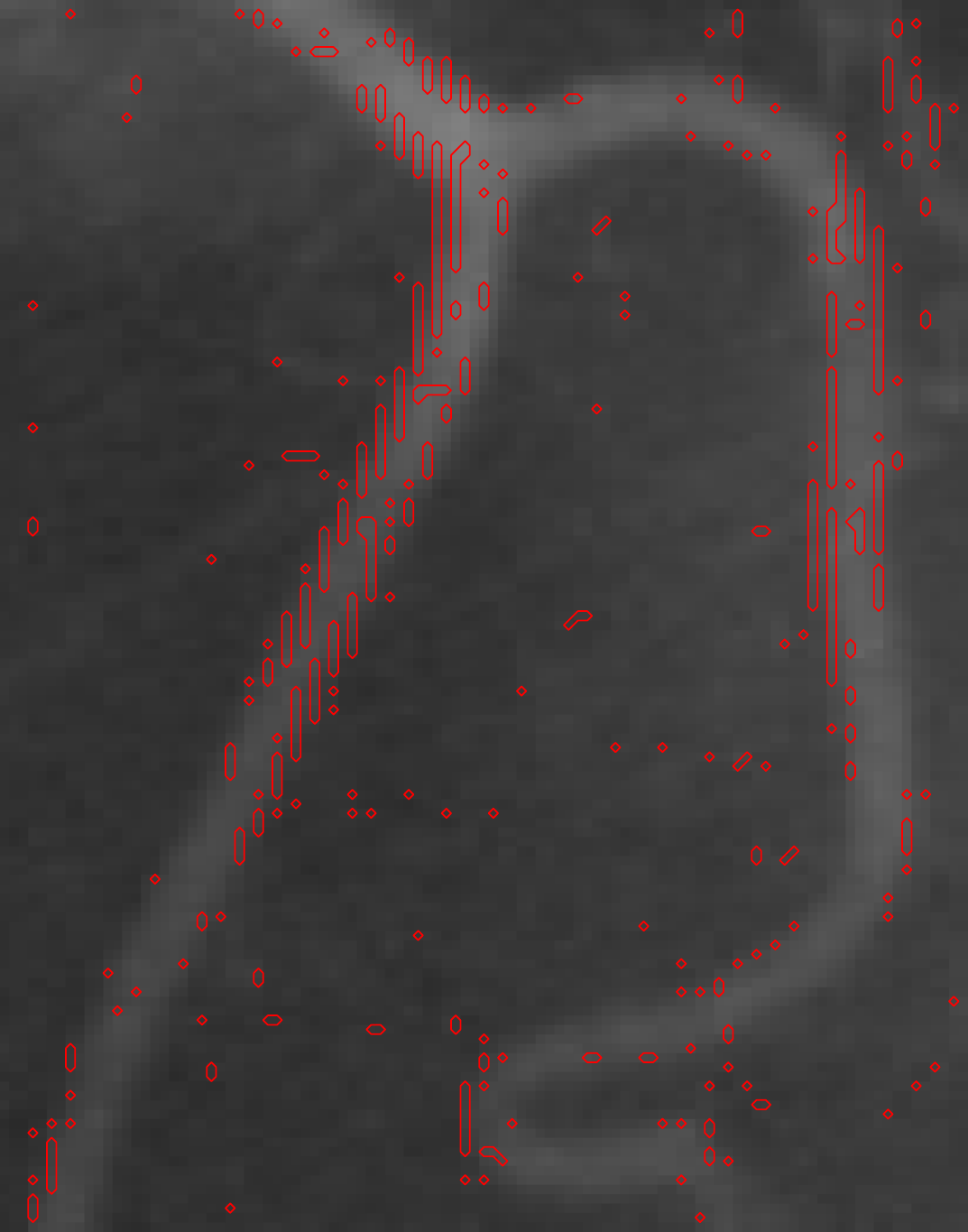}
		\includegraphics[width=3cm]{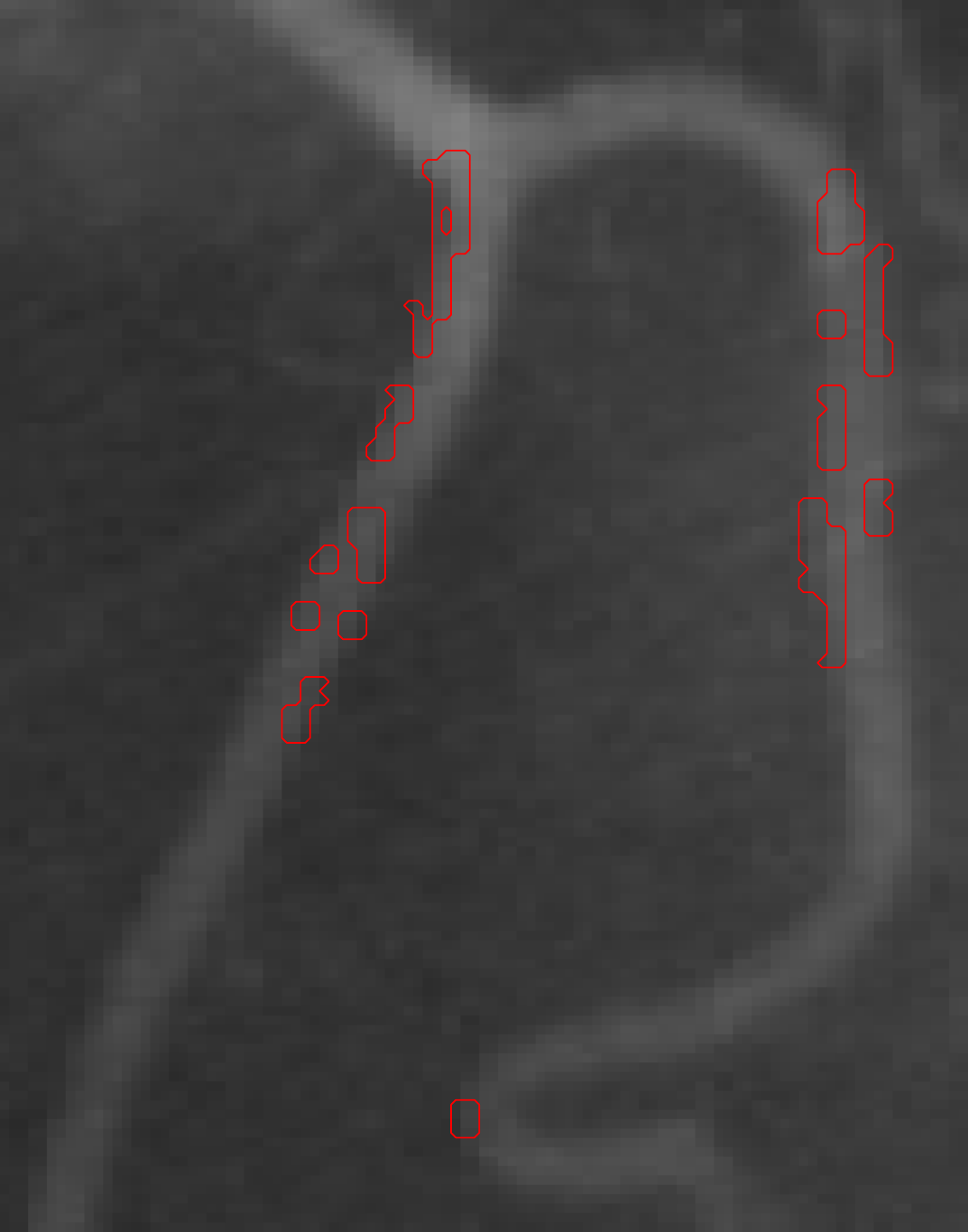}
		\includegraphics[width=3cm]{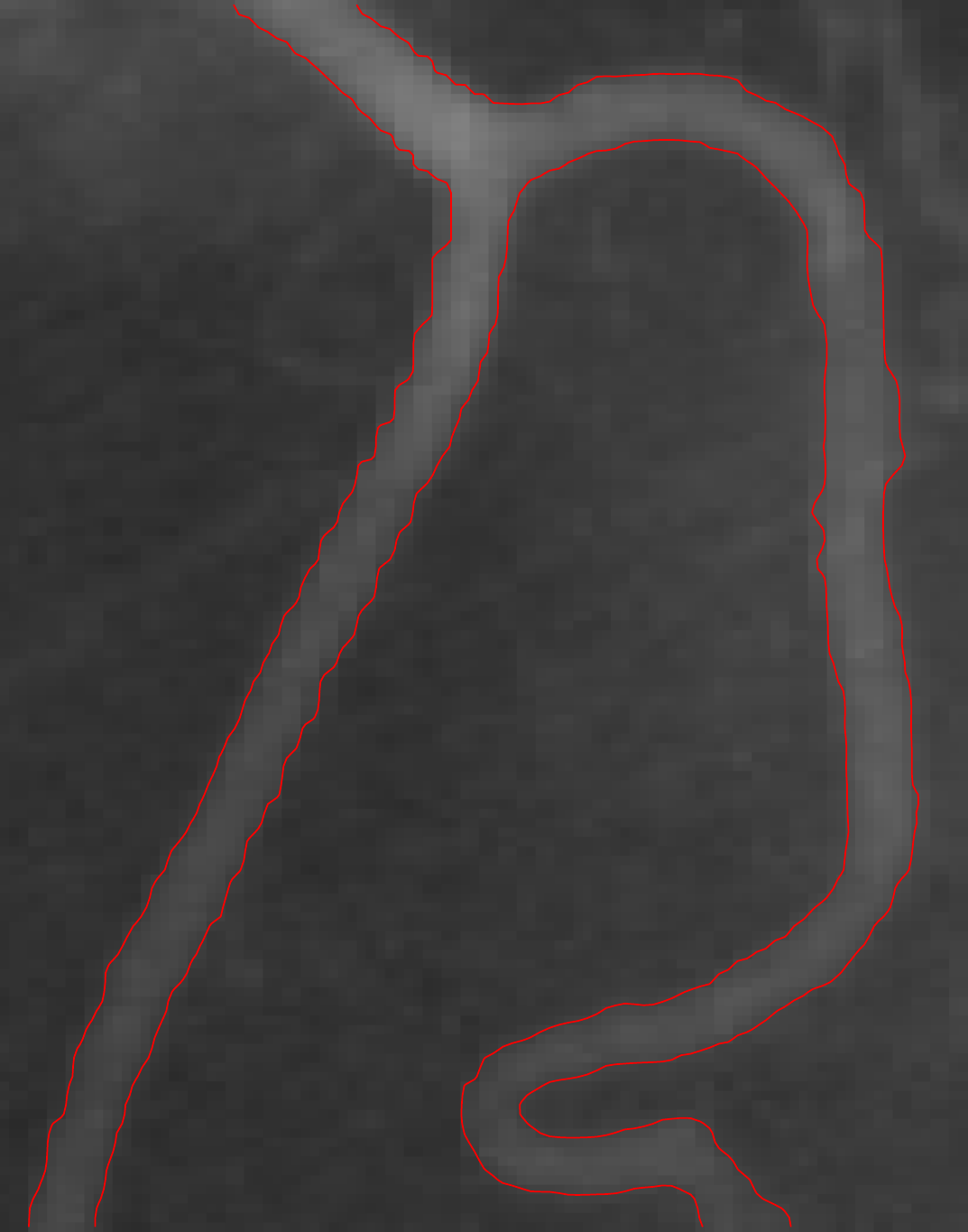}
	}
	\subfigure{
		\includegraphics[width=3cm]{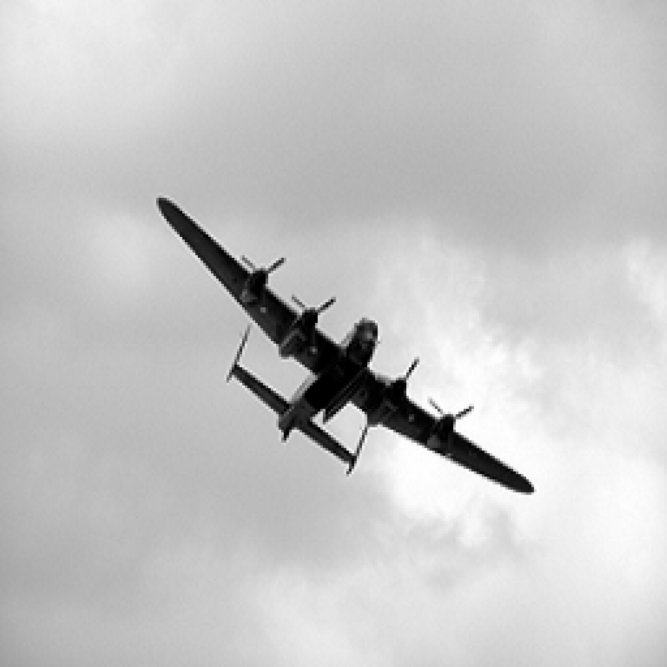}
		\includegraphics[width=3cm]{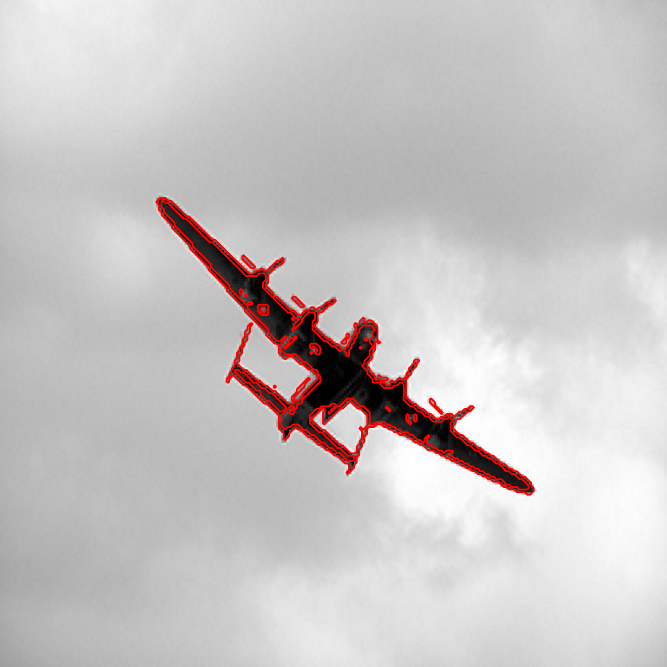}
		\includegraphics[width=3cm]{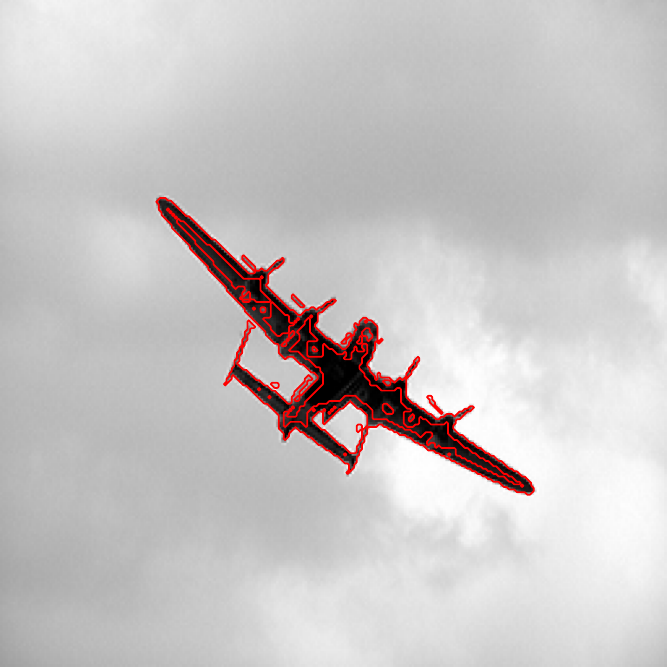}
		\includegraphics[width=3cm]{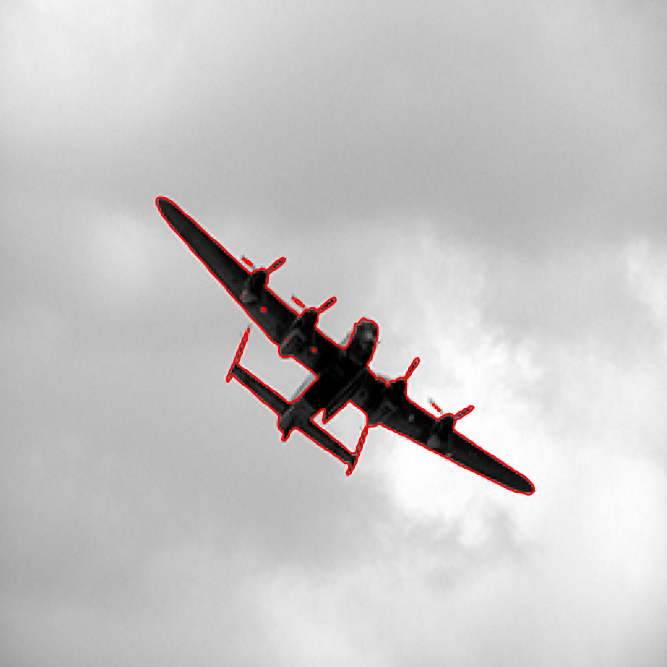}
	}
	\subfigure{
		\includegraphics[width=3cm]{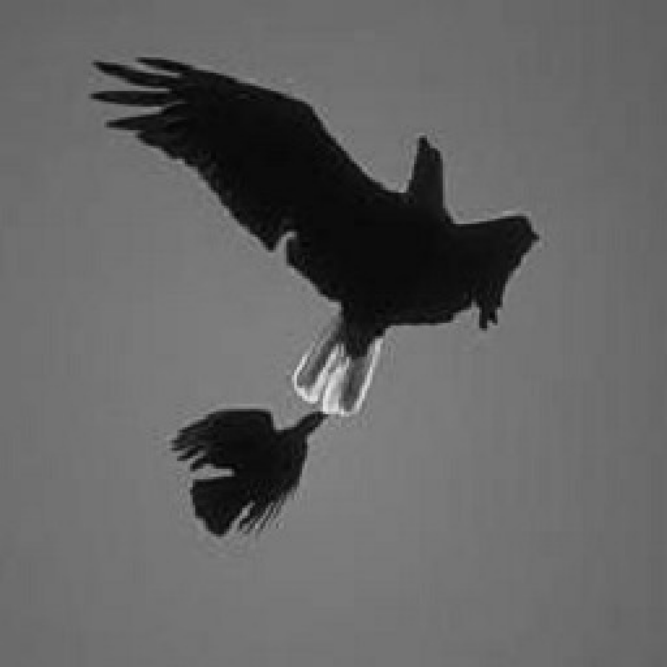}
		\includegraphics[width=3cm]{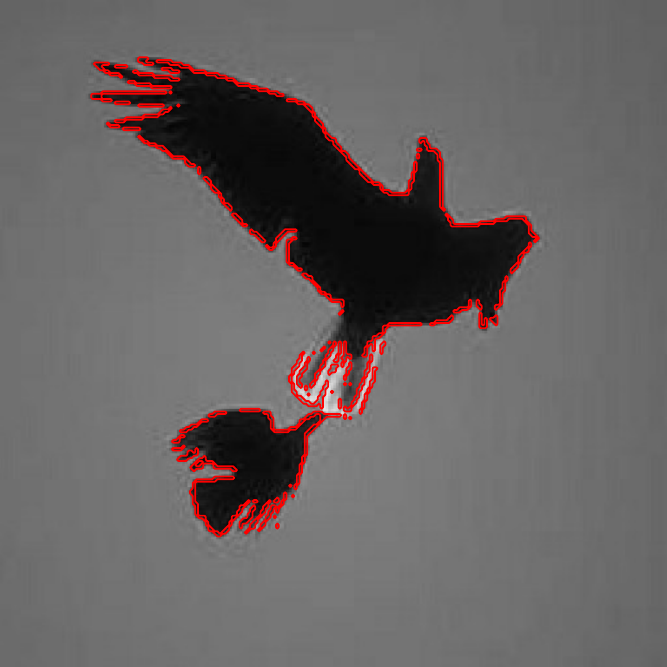}
		\includegraphics[width=3cm]{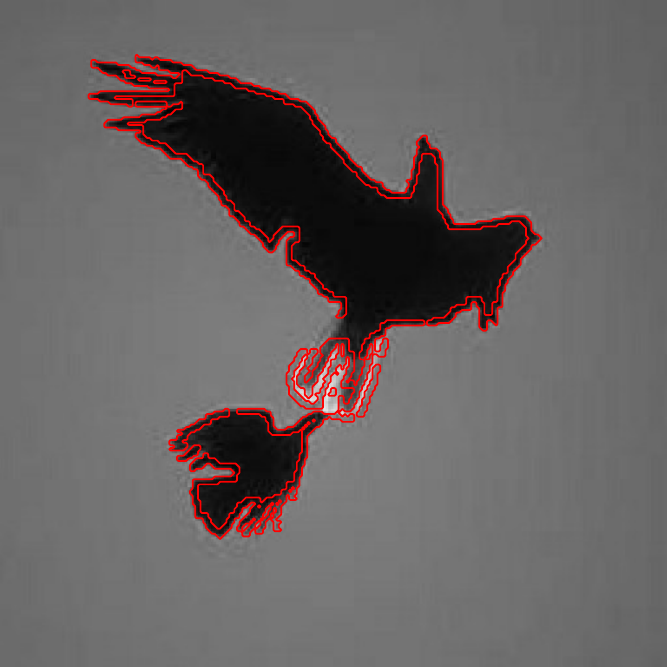}
		\includegraphics[width=3cm]{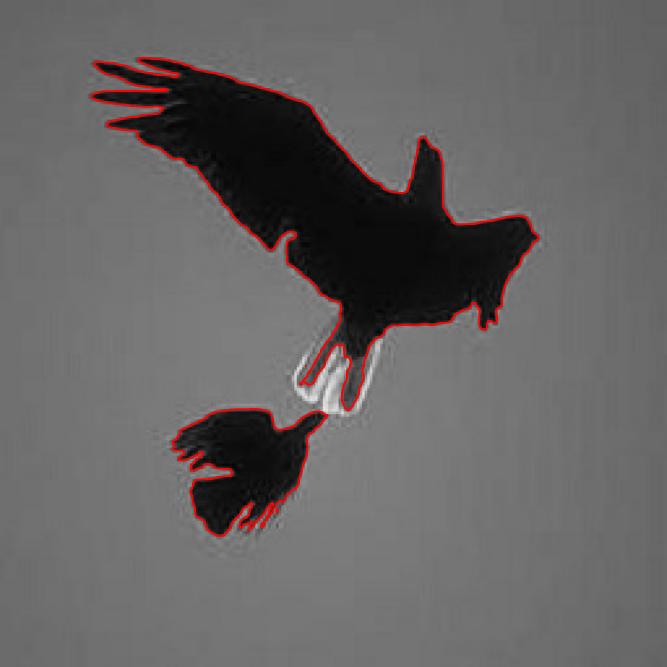}
	}
	 \caption{ First column: original images; Second column: rough initial contours; Third column: regions extended from initial contours after denoising; Fourth column: segmentation results. From top to bottom, the iteration number is 25, 42, 21, 35, 9 and 17, respectively.} \label{fig2}
\end{figure}

All the results are exhibited in Fig. \ref{fig2}. In each row, five figures are displayed to illustrate the whole segmentation process, where the original image is shown in the first figure, the second one is the rough initial contour derived from the inhomogeneous graph Laplacian operator, the third one indicates regions extended from initial contours after denoising, i.e., $S_p\cup R_p$ $(S_n\cup R_n)$, and the last one shows the final segmentation result of the ACLBF model. The contours are highlighted in red.

The experiment results indicate that initial contours given by IGLIM are reliable and almost precisely lie on the boundary of the object, which enables the ACLBF model to segment images effectively and efficiently, even for images with intensity inhomogeneity.

\begin{table}[t!]
	\resizebox{\textwidth}{!}{
	\begin{tabular}{llllp{1.2cm}llllllll}
		\hline
		\multirow{2}{*}{Images}&
		\multicolumn{4}{l}{IGLIM} & \multicolumn{7}{l}{ACLBF} \cr\cmidrule(lr){2-5}\cmidrule(lr){6-13}
		&$\lambda$ &$k_i$ &$M$ &initial contour &$\lambda_i$ &$\mu$ &$\sigma$ &$h$ &$\Delta t$ &$\varepsilon$ &$\varepsilon_1$ &$S$\cr
	\hline
		Row1 &50  &$k_1 = k_2 =0.01$  &1 &$S_p$ &$\lambda_1 =\lambda_2 =1$
		& 500 &1 &0.01 &0.1 &0.5 &1 &10*$\mu\varepsilon_1$ \cr
		Row2 &50  &$k_1 = k_2 =0.002$  &5 &$S_p$ &$\lambda_1 =\lambda_2 =1$
		& 500 &1 &0.01 &0.1 &0.5 &0.5 &10*$\mu\varepsilon_1$ \cr
		Row3 &50  &$k_1 = k_2 =0.01$  &5 &$S_p$ &$\lambda_1 =1,\lambda_2 =2.5$
		& 500&6 &0.01 &0.1 &0.5&0.5&150*$\mu\varepsilon_1$ \cr
		Row4 &50  &$k_1 = k_2 =0.01$  &3 &$S_n$ &$\lambda_1 =\lambda_2 =1$
		& 150 &3 &0.01 &0.1 &0.5 &0.5 &60*$\mu\varepsilon_1$ \cr
		Row5 &50  &$k_1 = k_2 =0.1$  &0 &$S_p$ &$\lambda_1 =\lambda_2 =1$
		& 500 &5 &0.01  &0.1 &0.5 &1 &500*$\mu\varepsilon_1$ \cr
		Row6 &50  &$k_1 = k_2 =0.1$  &0 &$S_p$ &$\lambda_1 =\lambda_2 =1$
		& 500 &10 &0.01  &0.1 &0.5 &1 &500*$\mu\varepsilon_1$ \cr\hline
	\end{tabular}
}
	\label{parameters}\caption{Parameters settings for images in Fig. \ref{fig2}.}
\end{table}

\subsection{Experiments on images with varying levels of noise}

\begin{figure}[t!]
\resizebox{\textwidth}{!}{
	\centering
	\subfigure{
	\begin{minipage}[b]{0.18\textwidth}
		\includegraphics[width=1.8cm]{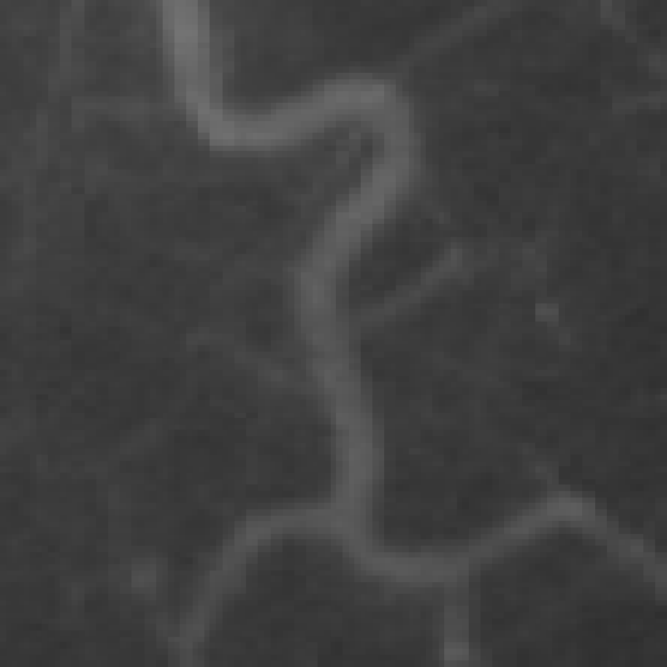}
		\includegraphics[width=1.8cm]{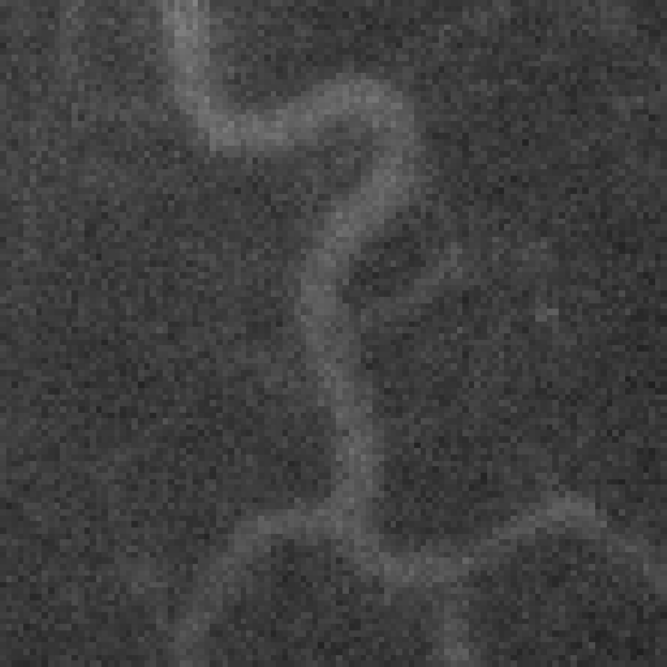}
		\includegraphics[width=1.8cm]{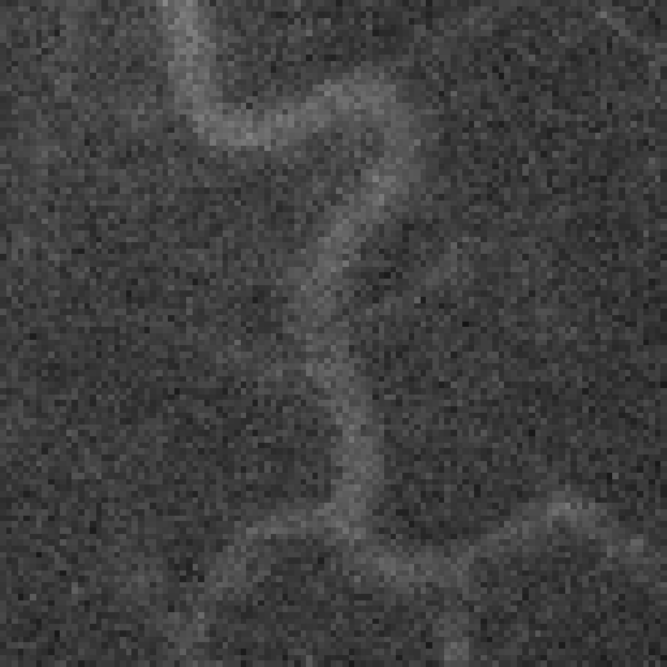}
		\includegraphics[width=1.8cm]{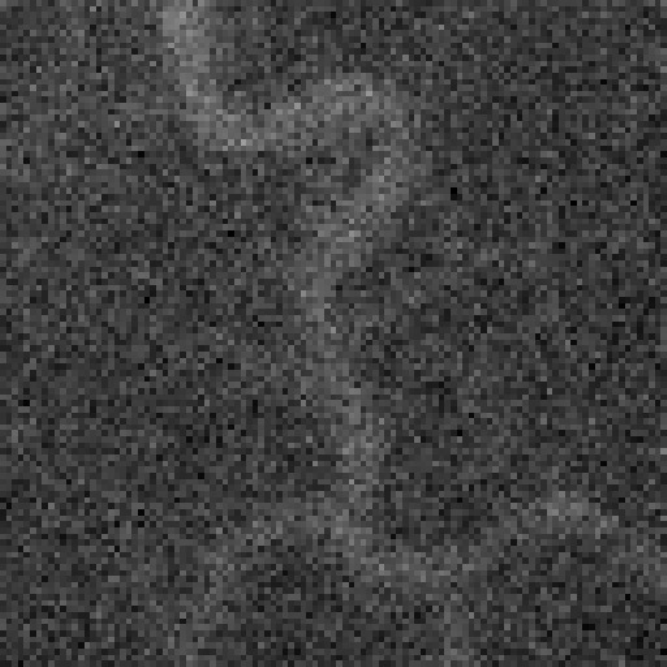}
		\centerline{(a)}
	\end{minipage}
	}
  \hspace{-9.8mm}
	\subfigure{
	\begin{minipage}[b]{0.18\textwidth}
		\includegraphics[width=1.8cm]{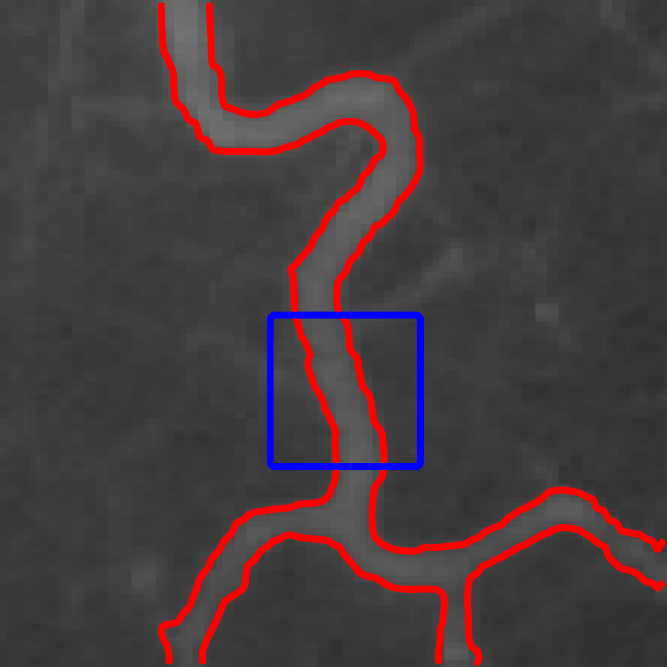}
		\includegraphics[width=1.8cm]{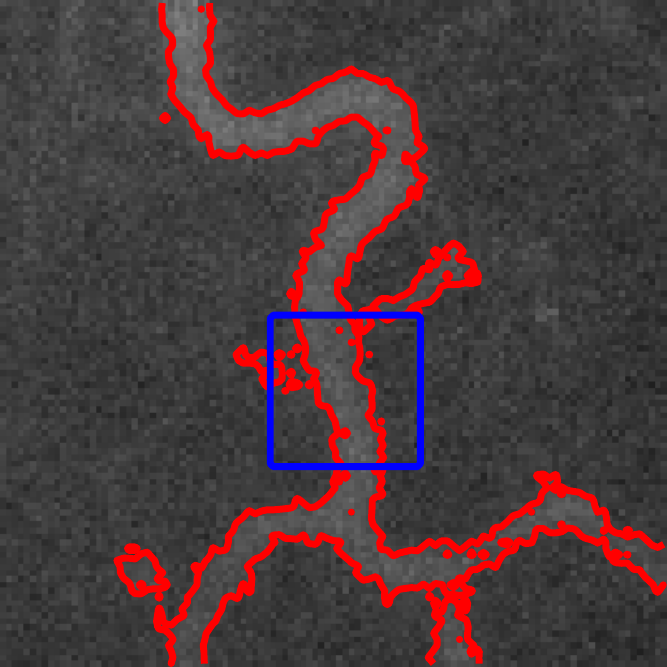}
		\includegraphics[width=1.8cm]{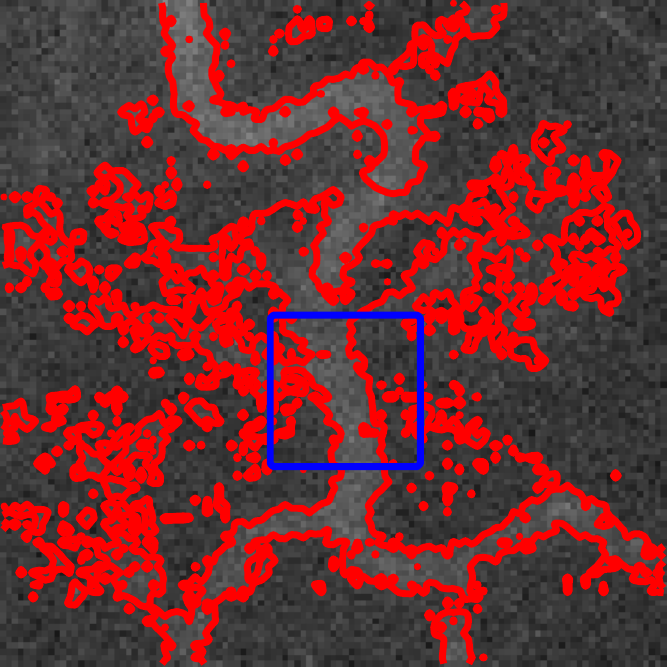}
		\includegraphics[width=1.8cm]{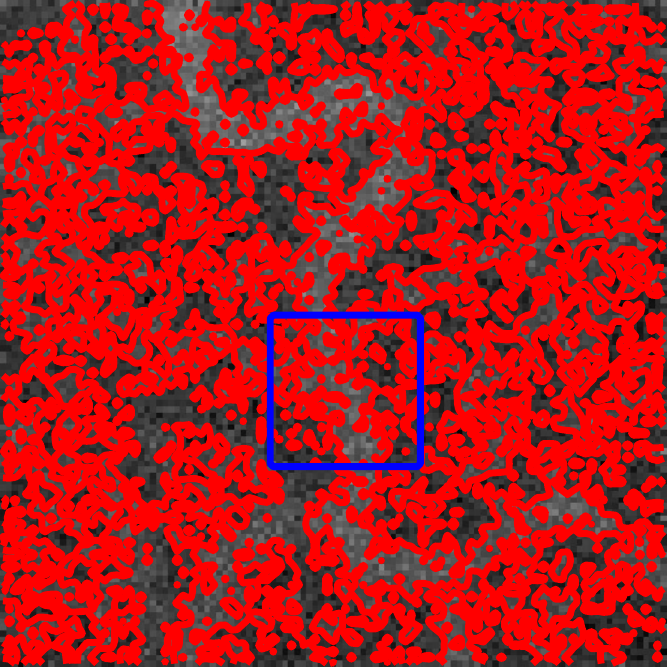}
		\centerline{(b)}
  \end{minipage}
	}
  \hspace{-9.8mm}
	\subfigure{
	\begin{minipage}[b]{0.18\textwidth}
	\includegraphics[width=1.8cm]{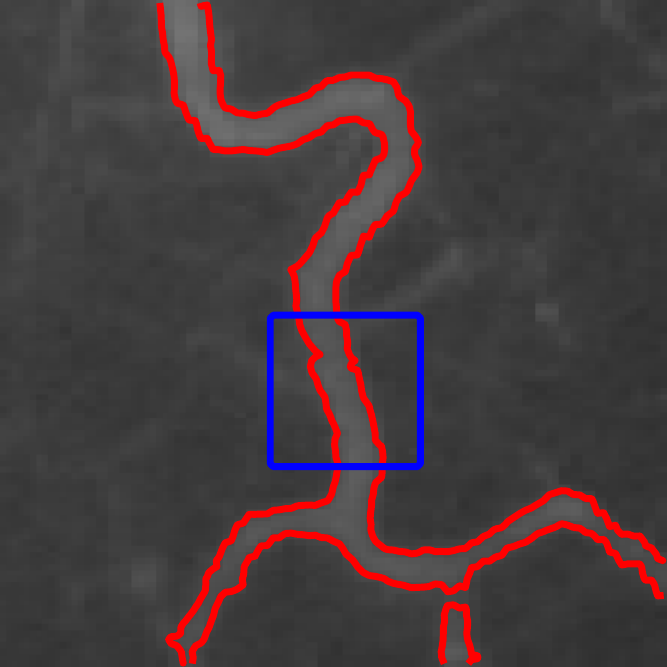}
	\includegraphics[width=1.8cm]{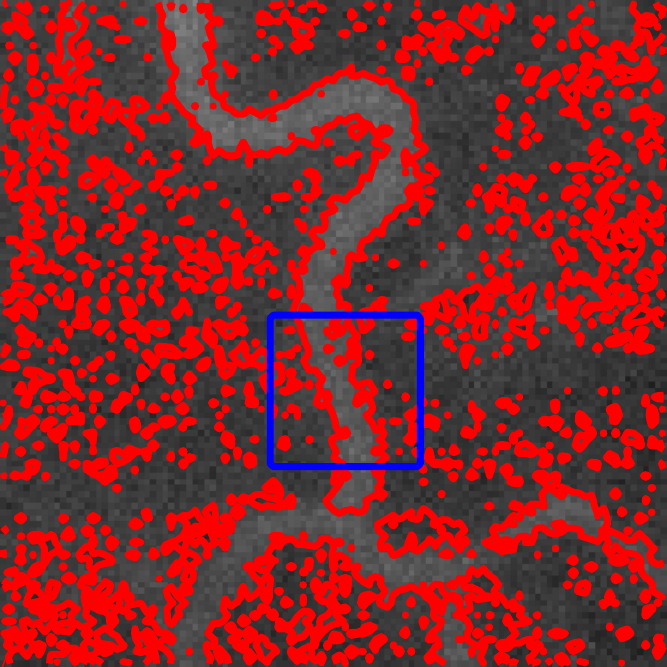}
	\includegraphics[width=1.8cm]{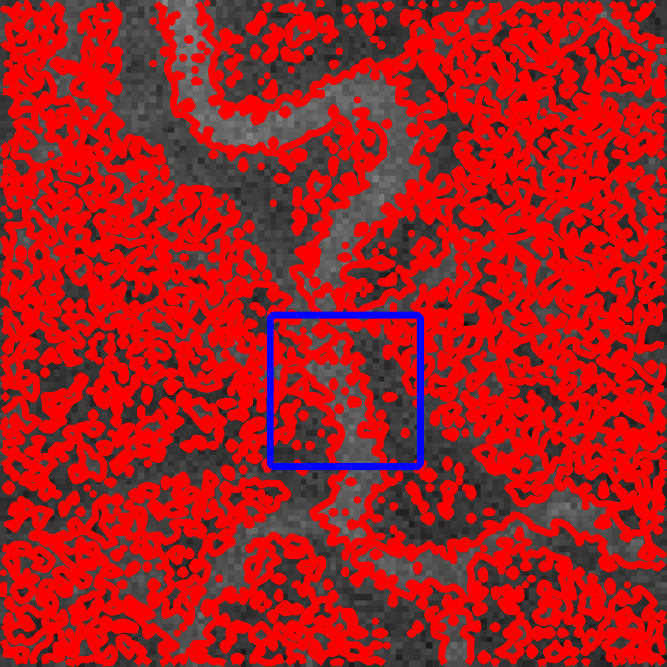}
	\includegraphics[width=1.8cm]{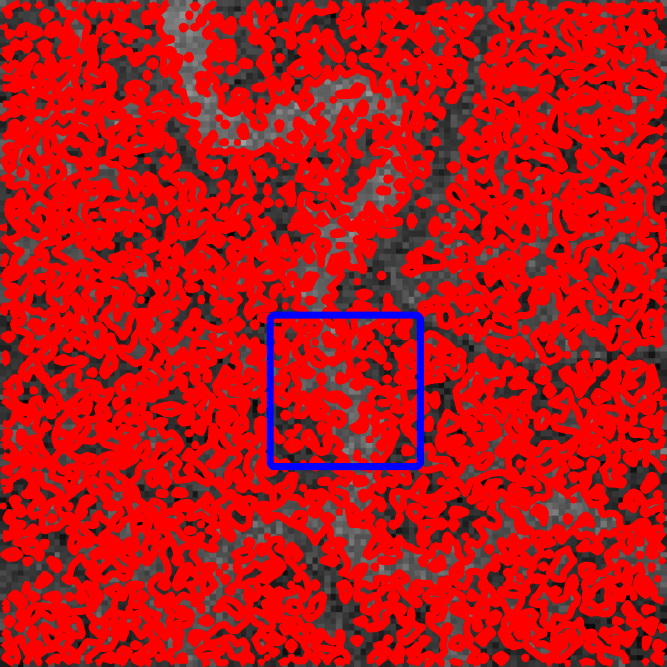}
	\centerline{(c)}
	\end{minipage}
	}
  \hspace{-9.8mm}
	\subfigure{
	\begin{minipage}[b]{0.18\textwidth}
	\includegraphics[width=1.8cm]{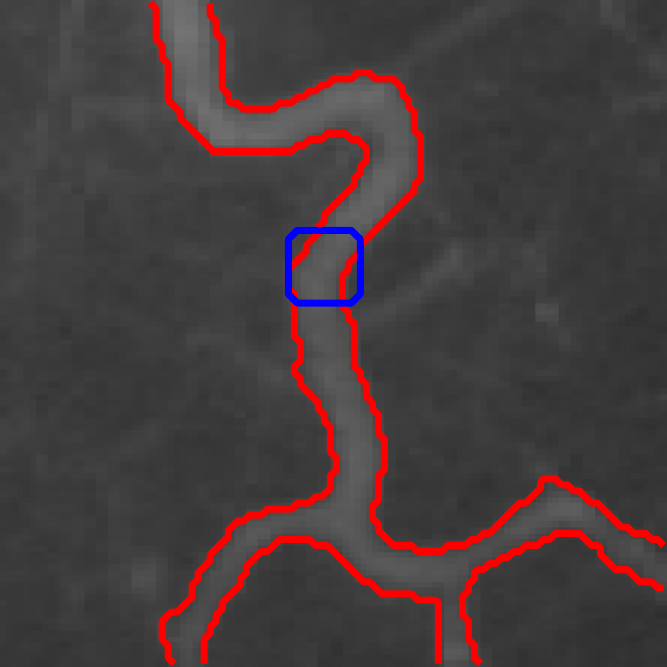}
	\includegraphics[width=1.8cm]{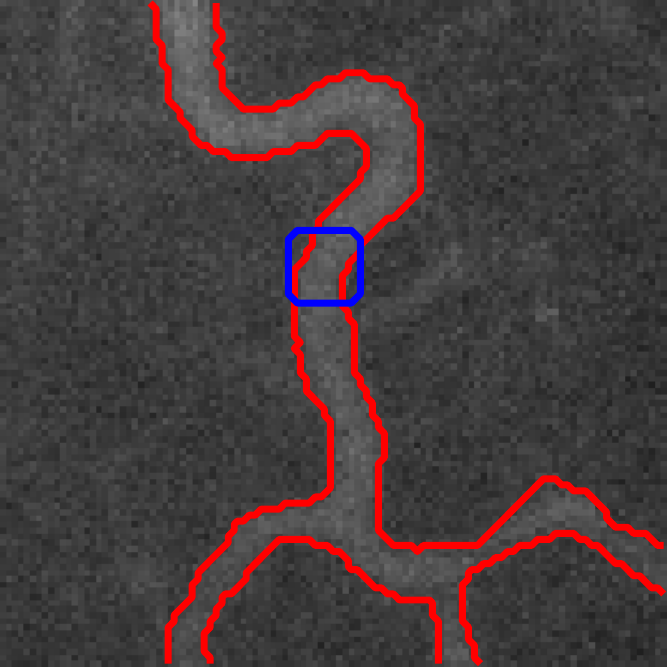}
	\includegraphics[width=1.8cm]{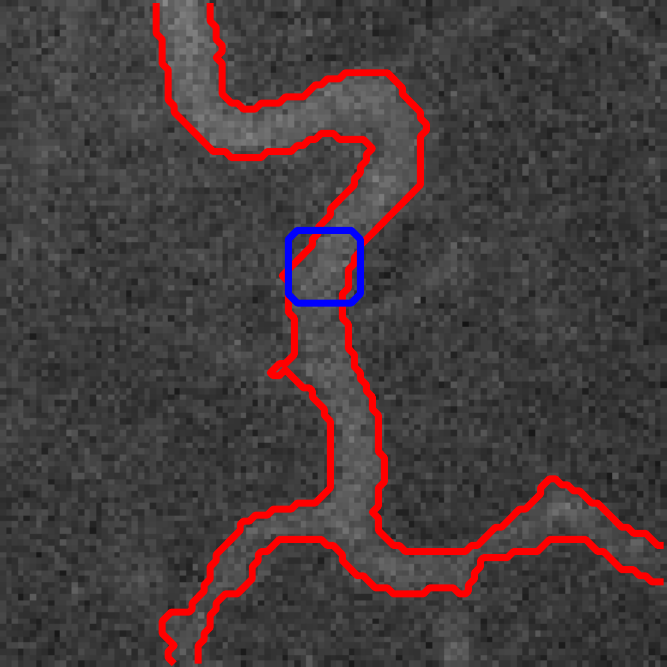}
	\includegraphics[width=1.8cm]{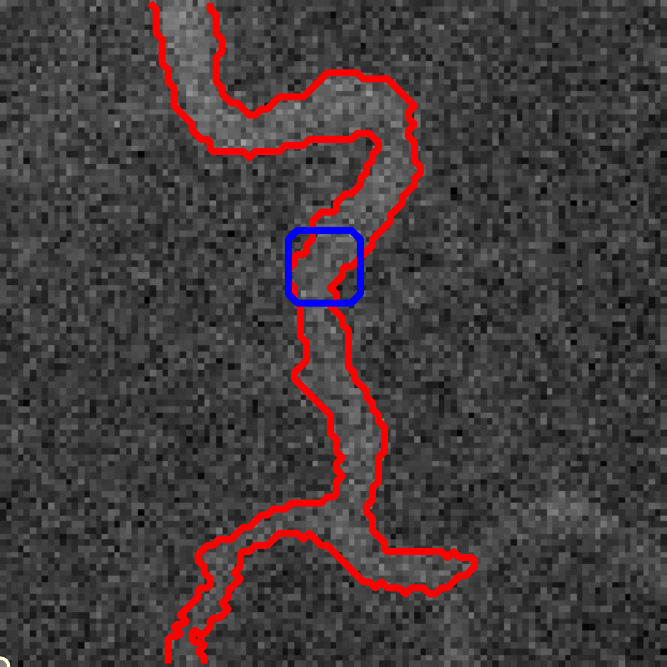}
	\centerline{(d)}
	\end{minipage}
	}
  \hspace{-9.8mm}
	\subfigure{
	\begin{minipage}[b]{0.18\textwidth}
	\includegraphics[width=1.8cm]{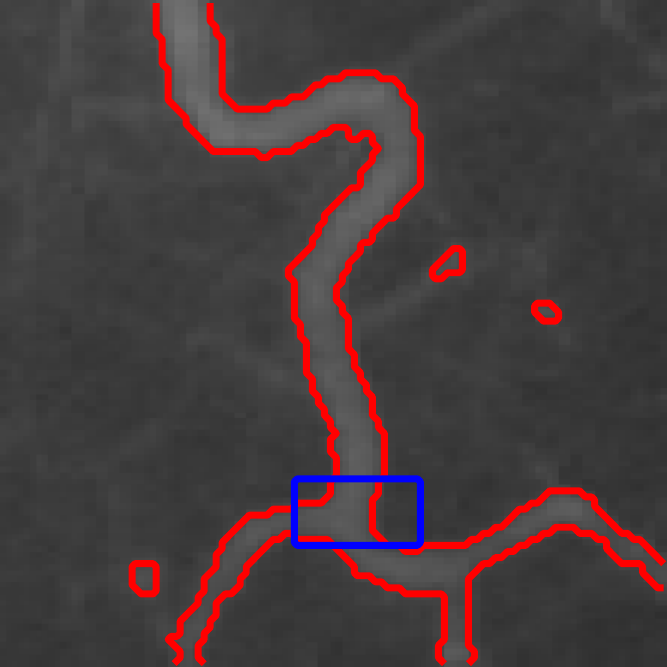}
	\includegraphics[width=1.8cm]{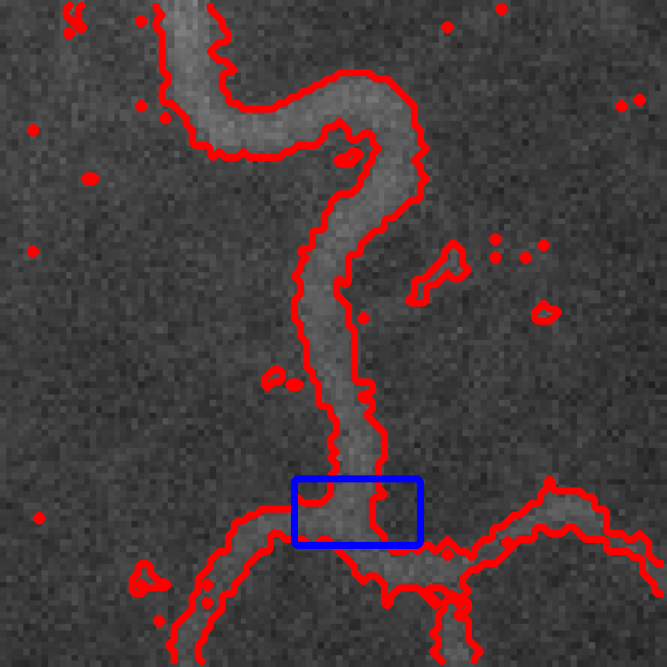}
	\includegraphics[width=1.8cm]{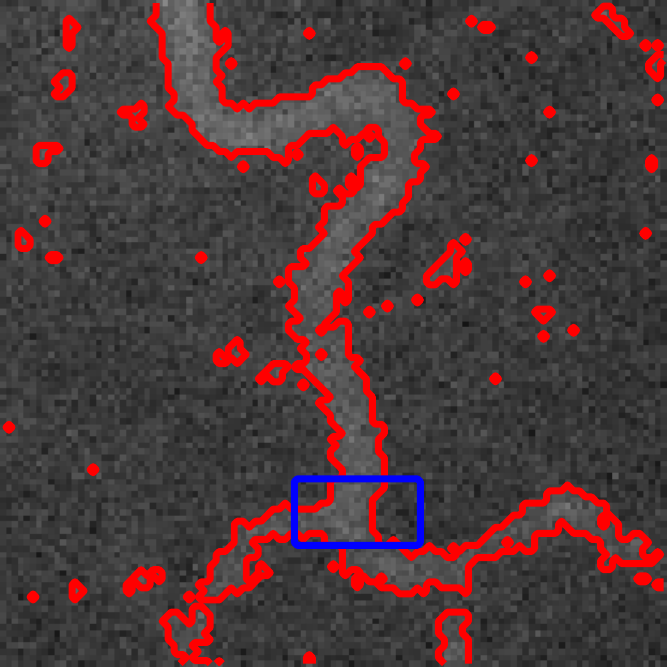}
	\includegraphics[width=1.8cm]{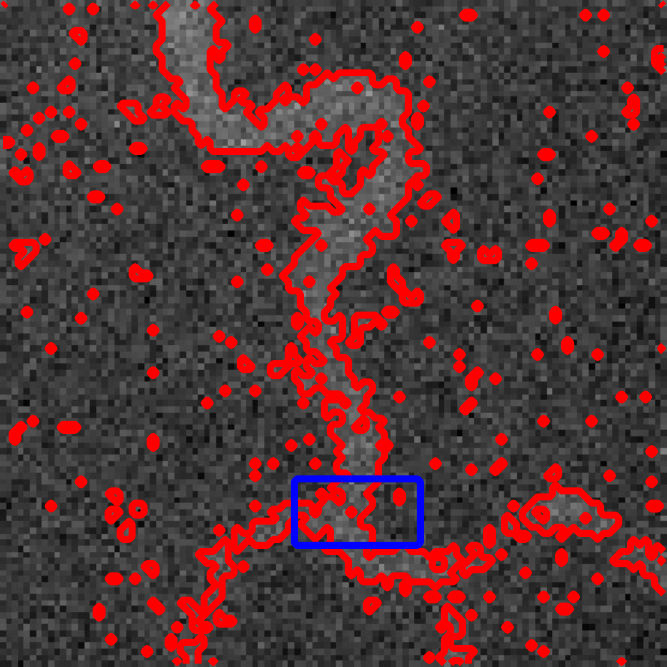}
	\centerline{(e)}
	\end{minipage}
	}
  \hspace{-9.8mm}
	\subfigure{
	\begin{minipage}[b]{0.18\textwidth}
	\includegraphics[width=1.8cm]{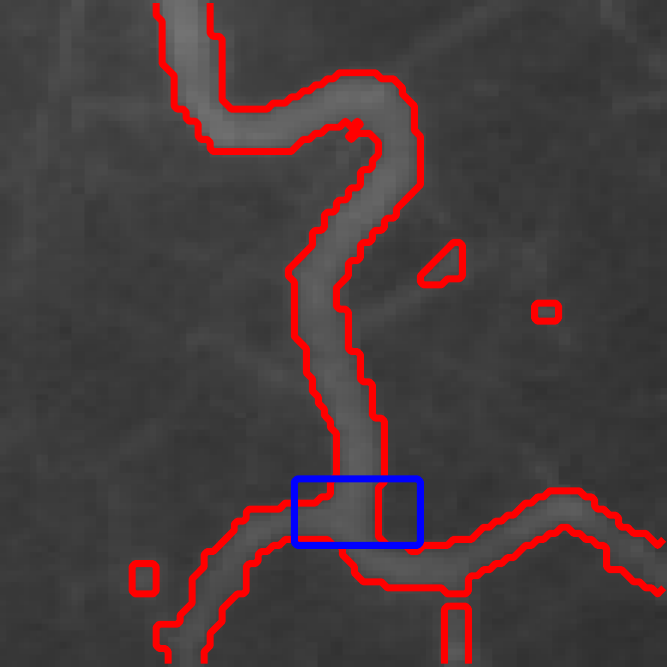}
	\includegraphics[width=1.8cm]{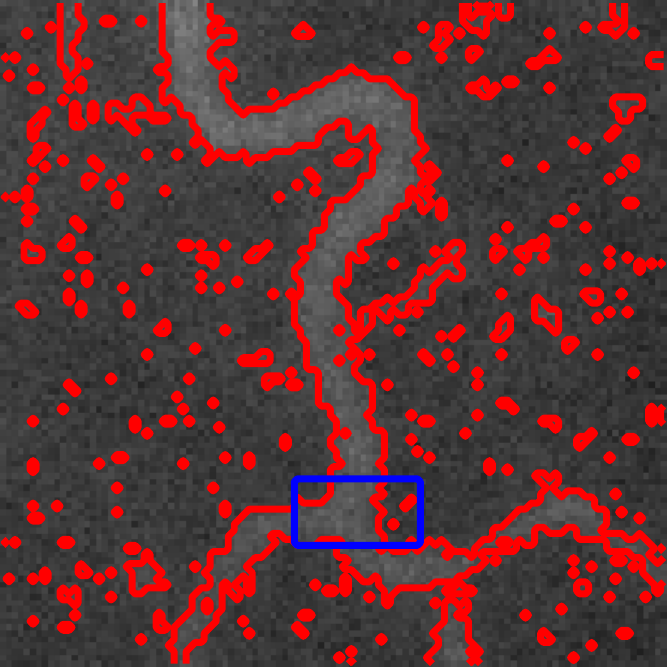}
	\includegraphics[width=1.8cm]{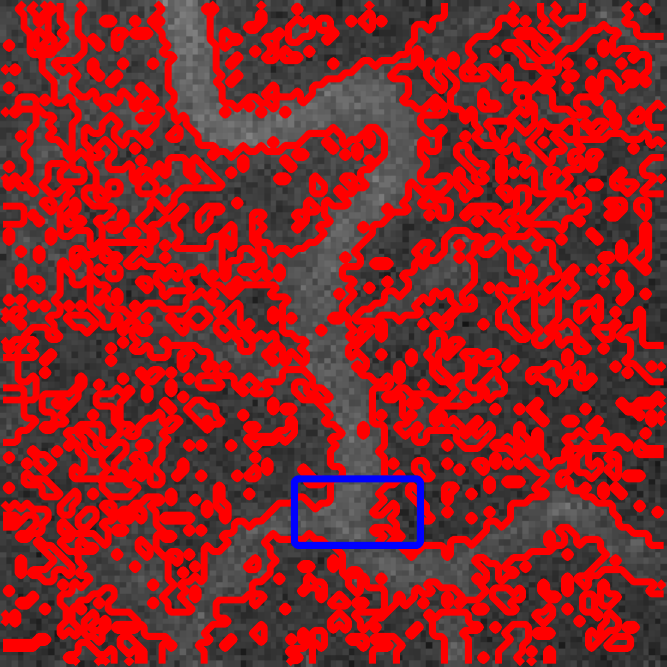}
	\includegraphics[width=1.8cm]{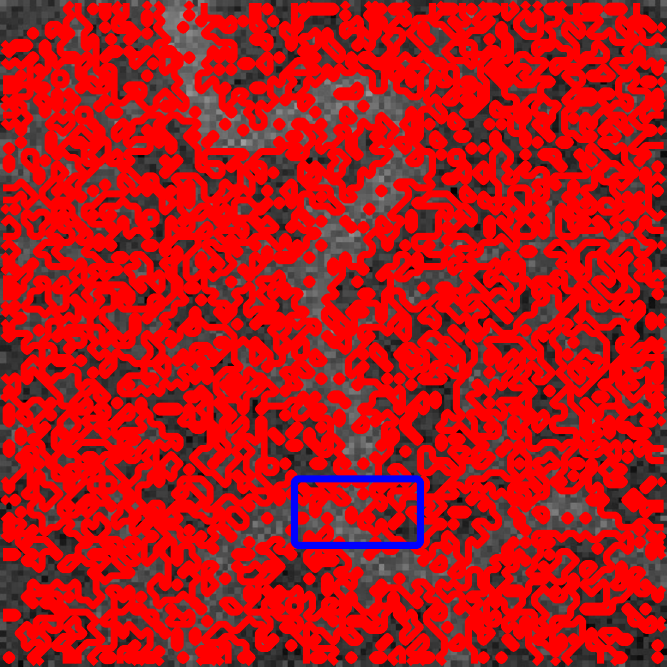}
	\centerline{(f)}
	\end{minipage}
	}
  \hspace{-9.8mm}
	\subfigure{
	\begin{minipage}[b]{0.18\textwidth}
	\includegraphics[width=1.8cm]{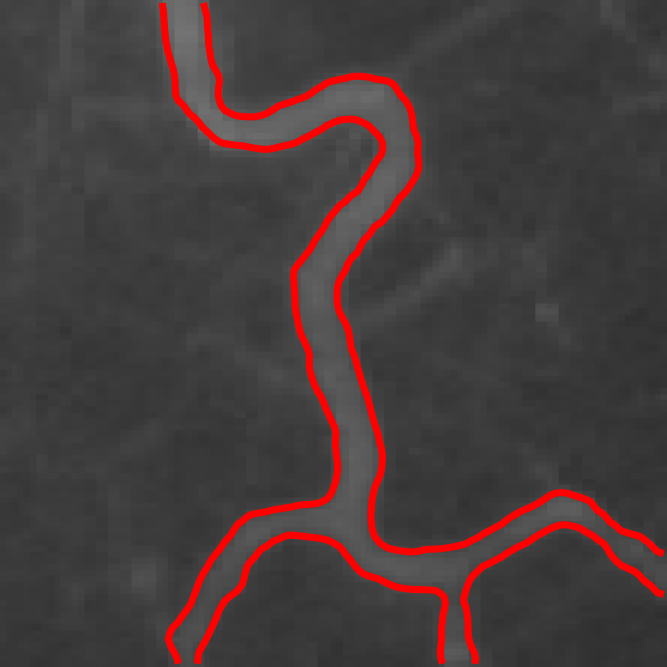}
	\includegraphics[width=1.8cm]{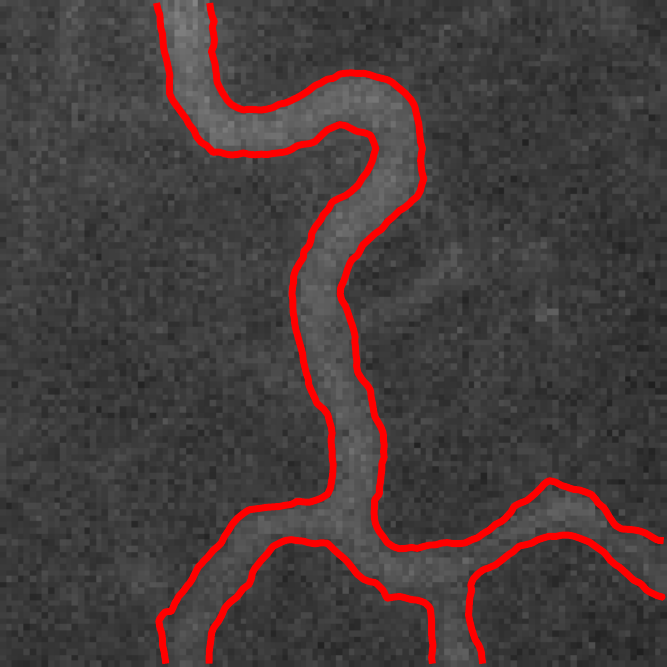}
	\includegraphics[width=1.8cm]{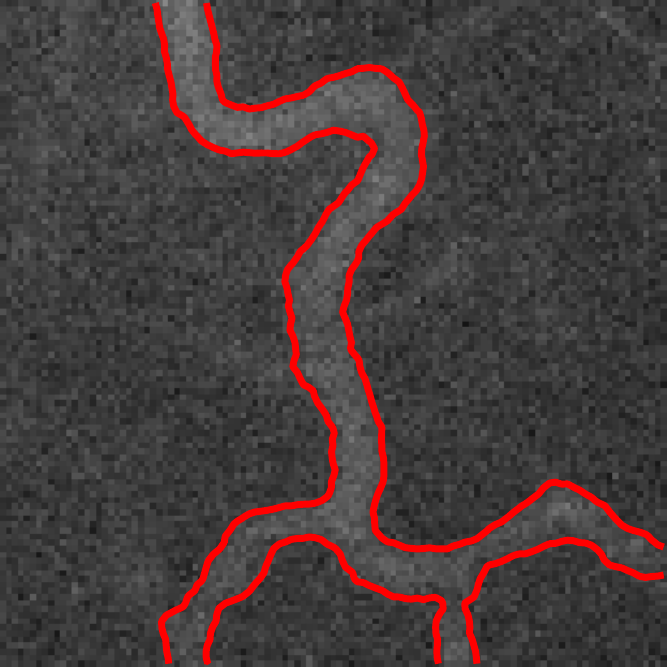}
	\includegraphics[width=1.8cm]{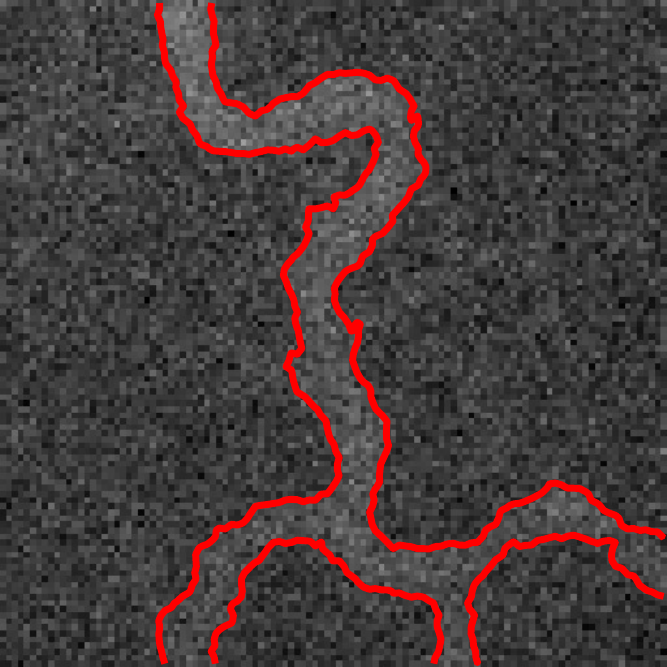}
	\centerline{(g)}
	\end{minipage}
	}
}
	 \caption{ Comparison with other models of segmentation results of images with different noise levels (The variance from the first row to the last row: 50, 100, 300). (a) Original image; segmentation results of (b) LBF model solved by level set method; (c) LBF model solved by ICTM; (d) RLSF model; (e) CV-XB model; (f) the model in \cite{min2021inhomogeneous}; (g) the ACLBF model.} \label{fig3}
\end{figure}

To demonstrate the robustness of our algorithm on segmentation of images with noise, we test the performance of the ACLBF model on five images corrupted by different levels of Gaussian noise. Fig. \ref{fig3} exhibits segmentation results of the ACLBF model and several celebrated models for these images. As the level of noise increases, the identification of the blood vessel becomes more challenging. Our results of the ACLBF model are still reliable for relatively strong noise.

\subsection{Experiments on images with both intensity inhomogeneity and noise}
\begin{figure}[htbp]
	\centering
	\subfigure{
	\begin{minipage}[b]{0.18\textwidth}
	\includegraphics[width=2cm]{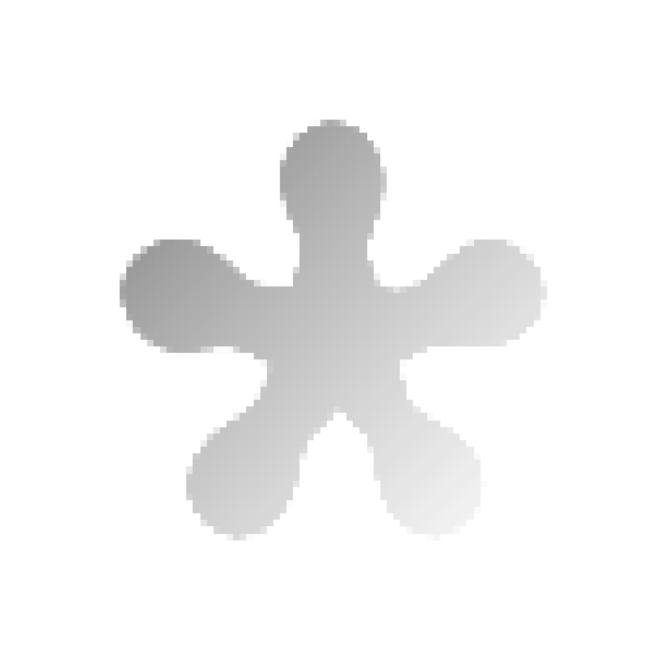}
	\includegraphics[width=2cm]{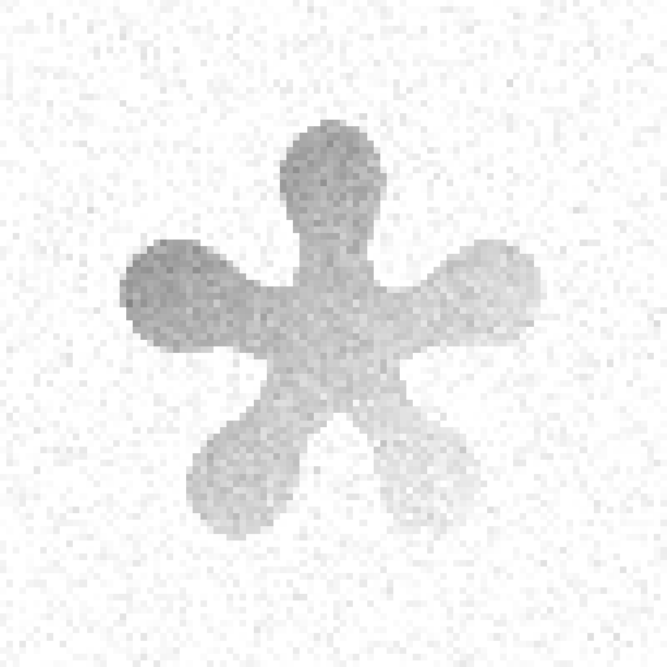}
	\includegraphics[width=2cm]{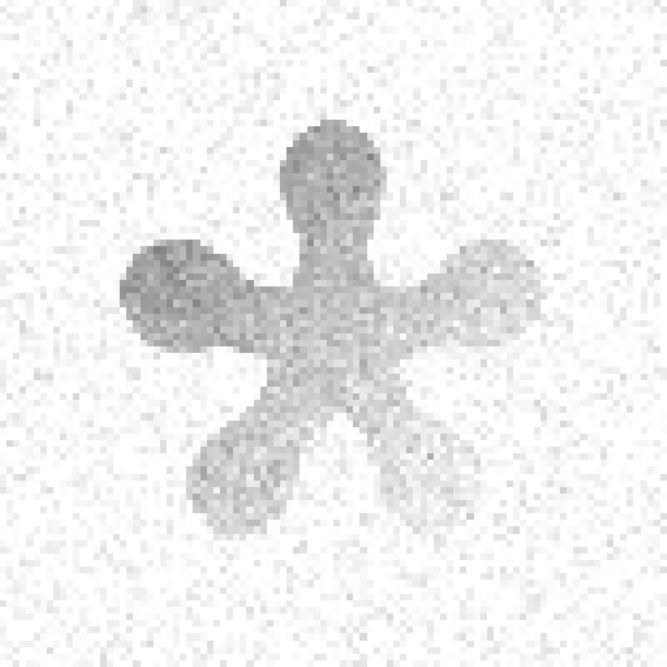}
	\includegraphics[width=2cm]{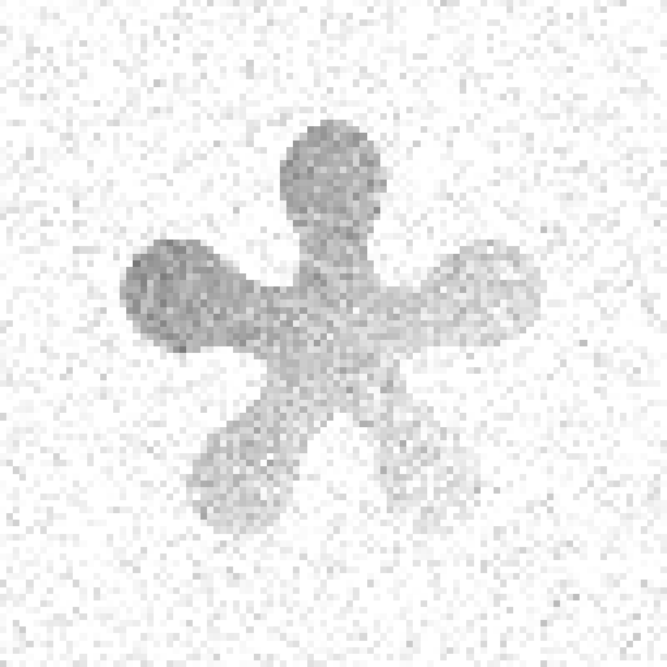}
	\includegraphics[width=2cm]{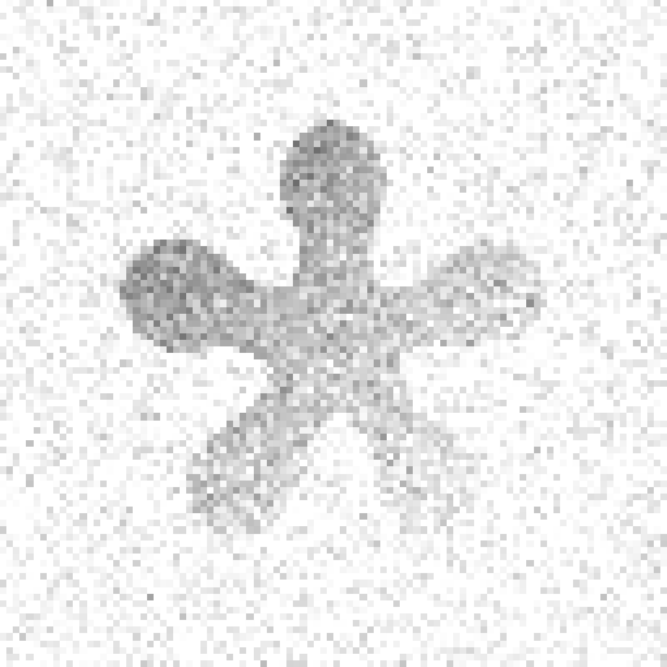}
	\centerline{(a)  }
	\end{minipage}
	}
	\hspace{-10.25mm}
	\subfigure{
	\begin{minipage}[b]{0.18\textwidth}
	\includegraphics[width=2cm]{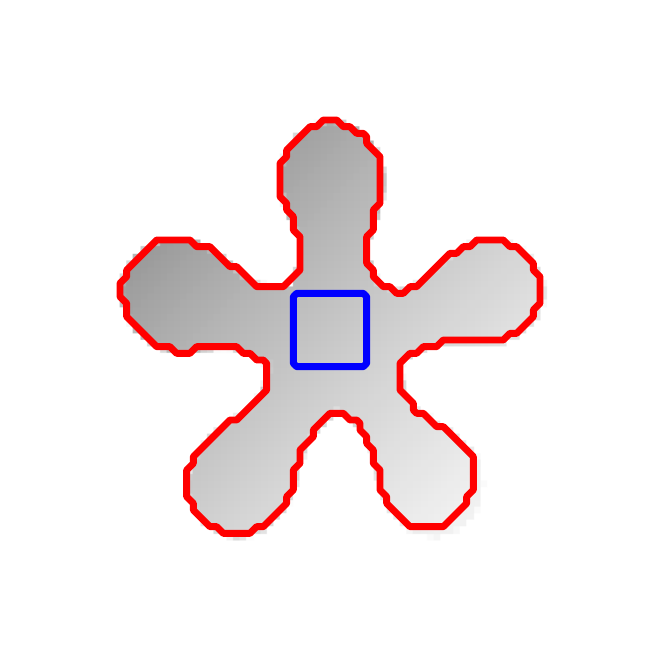}
	\includegraphics[width=2cm]{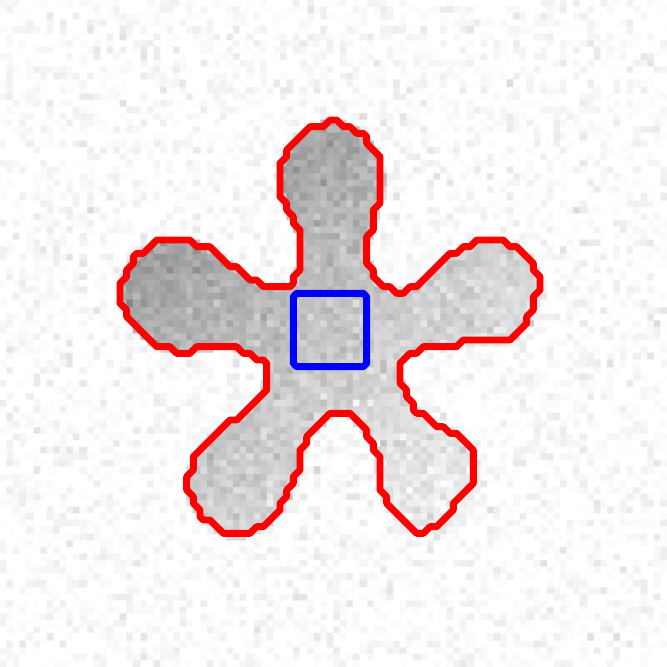}
	\includegraphics[width=2cm]{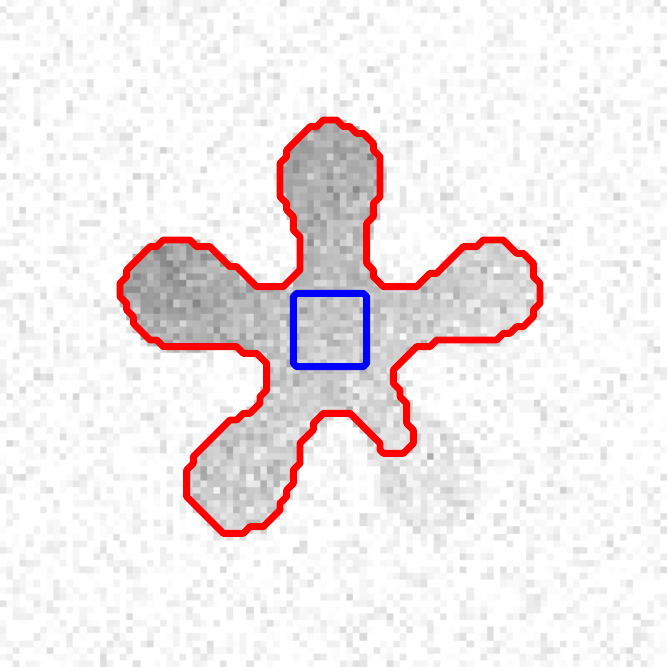}
	\includegraphics[width=2cm]{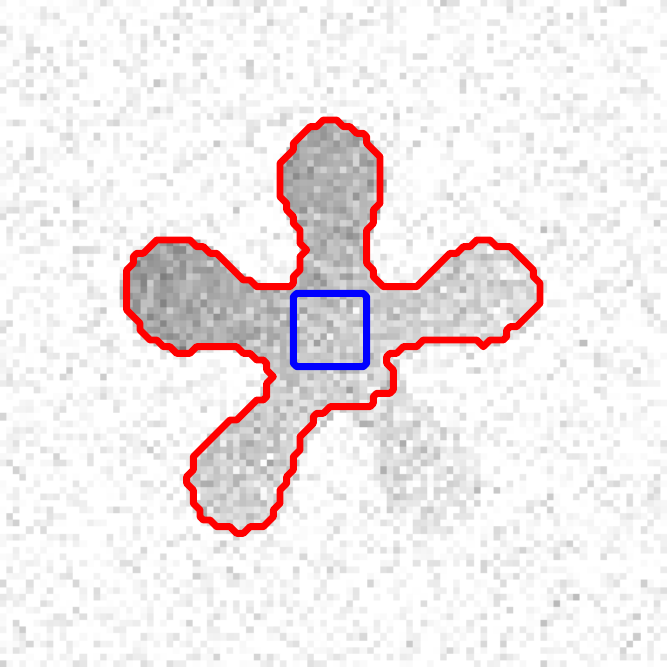}
	\includegraphics[width=2cm]{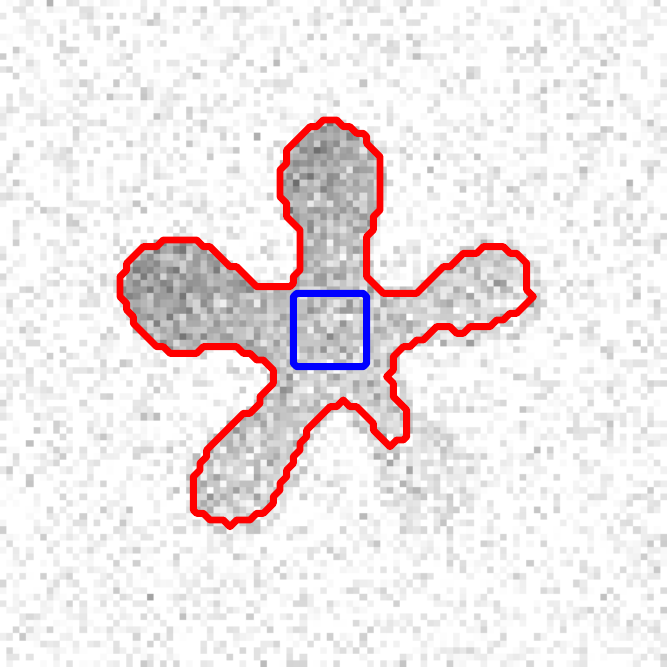}
	\centerline{(b)}
	\end{minipage}
	}
	\hspace{-10.25mm}
	\subfigure{
	\begin{minipage}[b]{0.18\textwidth}
	\includegraphics[width=2cm]{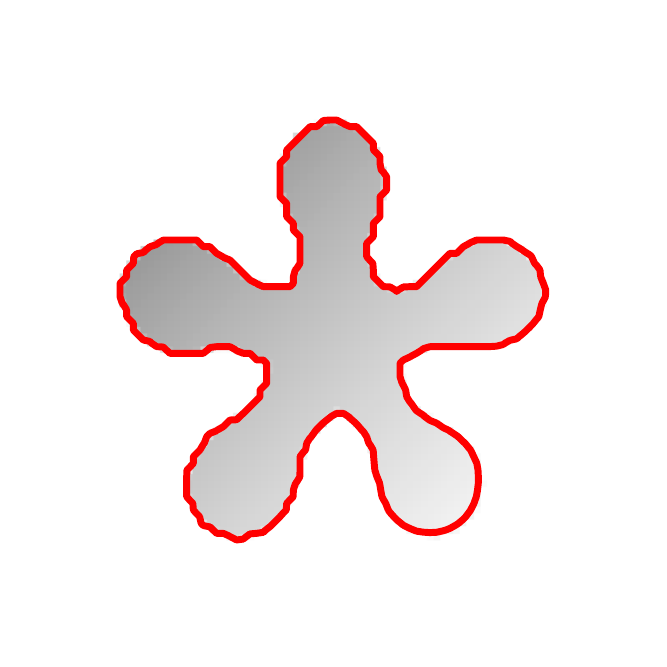}
	\includegraphics[width=2cm]{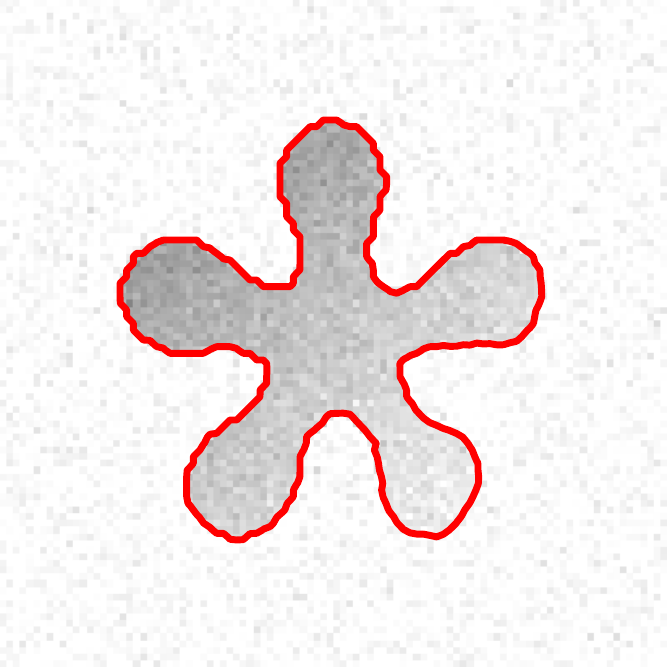}
	\includegraphics[width=2cm]{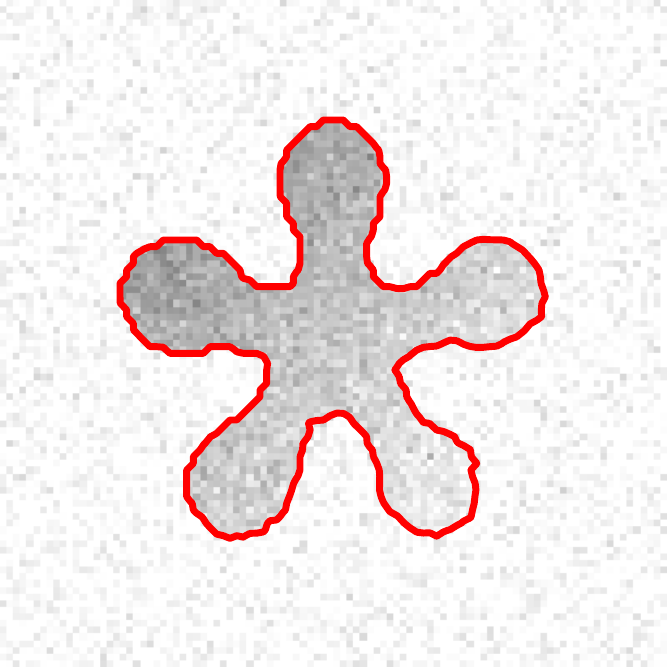}
	\includegraphics[width=2cm]{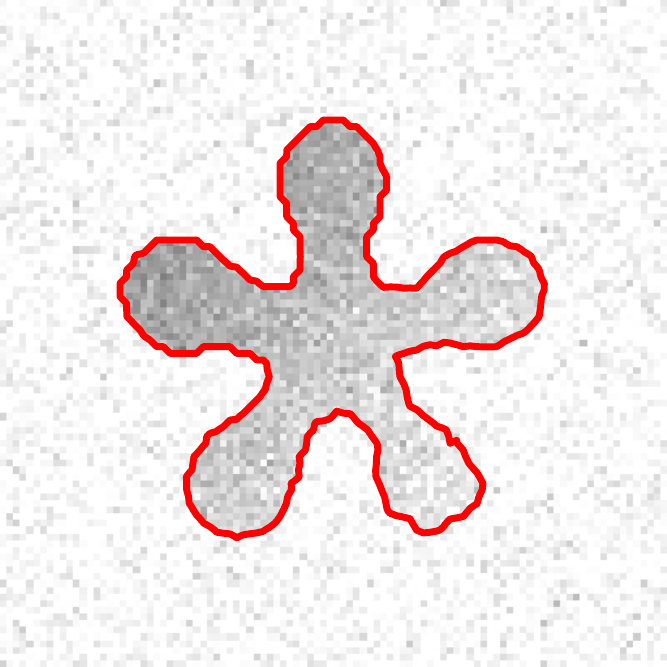}
	\includegraphics[width=2cm]{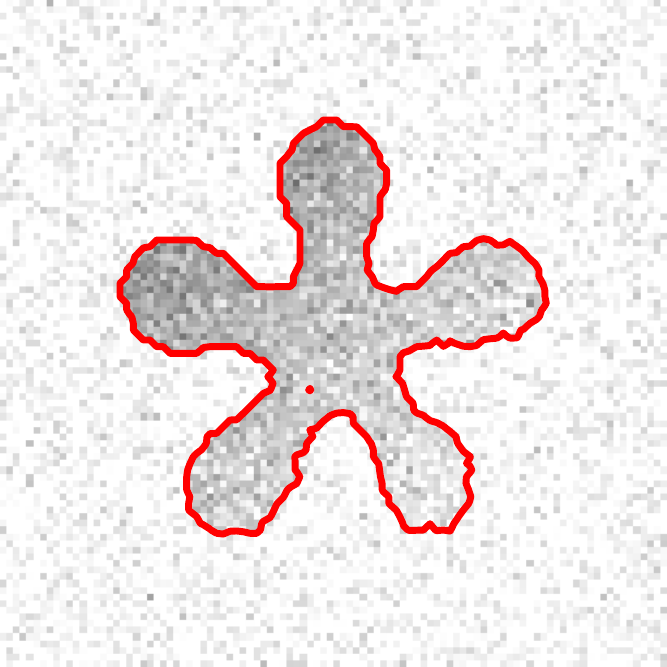}
	\centerline{(c)}
	\end{minipage}
	}
	\hspace{-10.25mm}
	\subfigure{
	\begin{minipage}[b]{0.18\textwidth}
	\includegraphics[width=2cm]{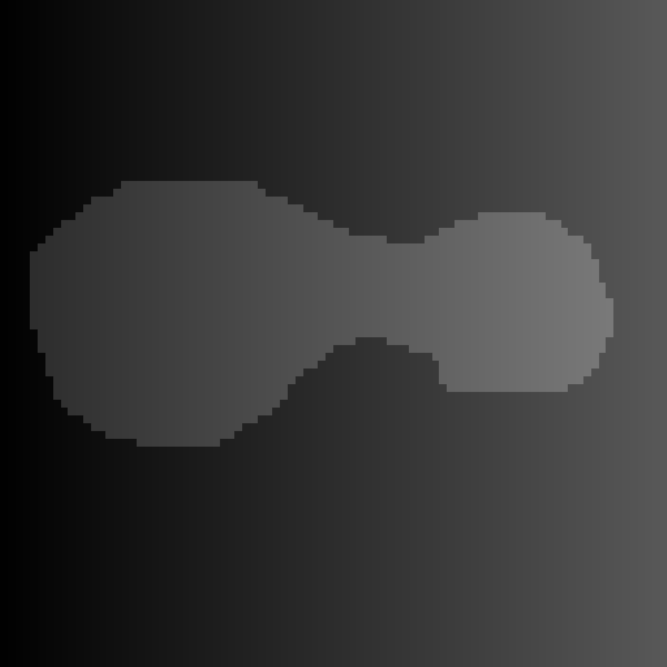}
	\includegraphics[width=2cm]{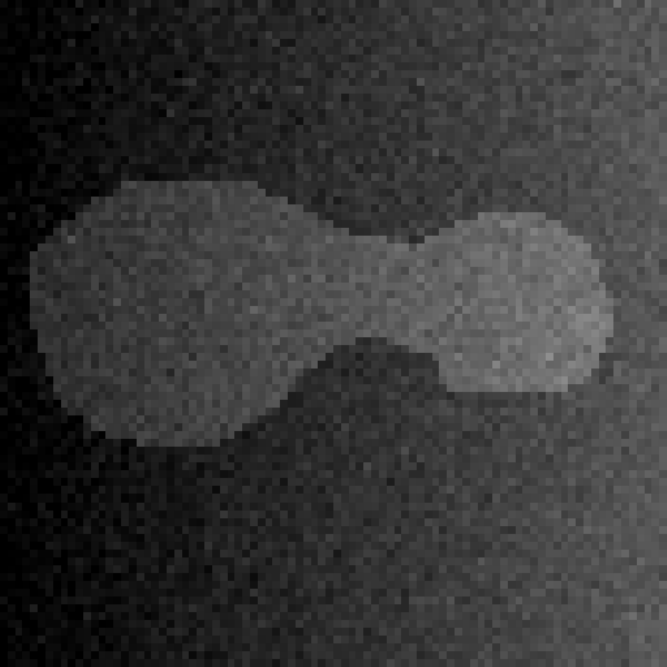}
	\includegraphics[width=2cm]{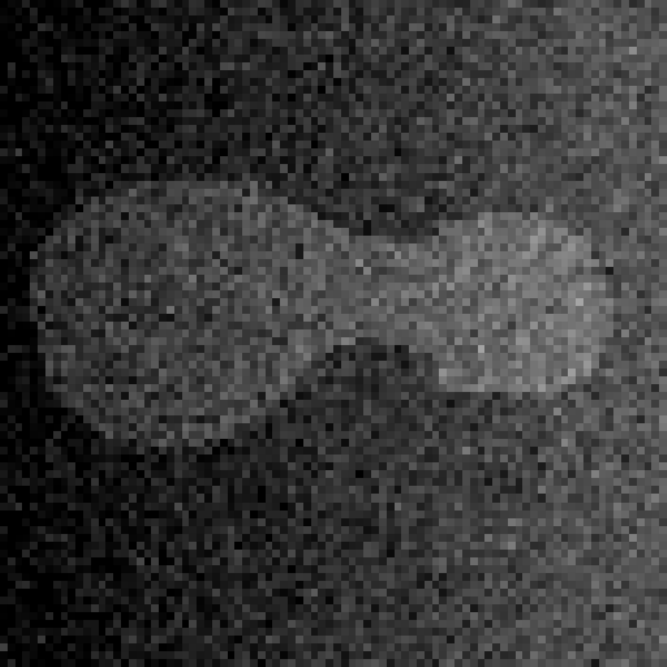}
	\includegraphics[width=2cm]{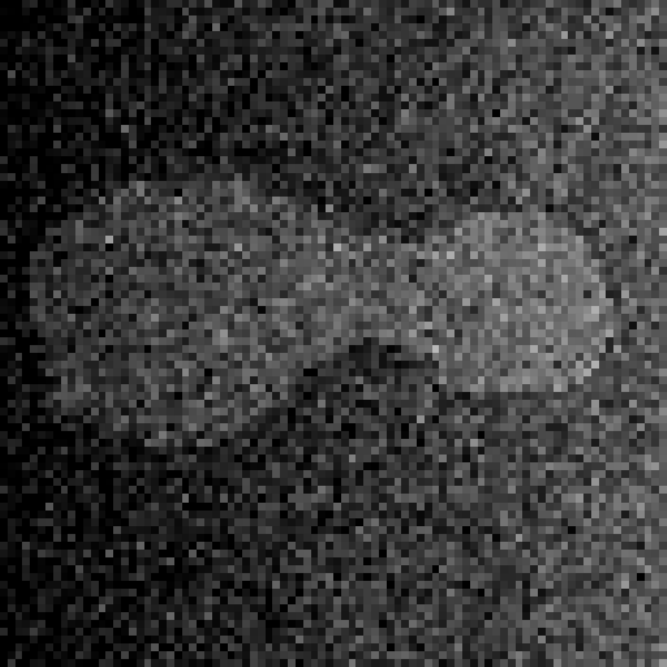}
	\includegraphics[width=2cm]{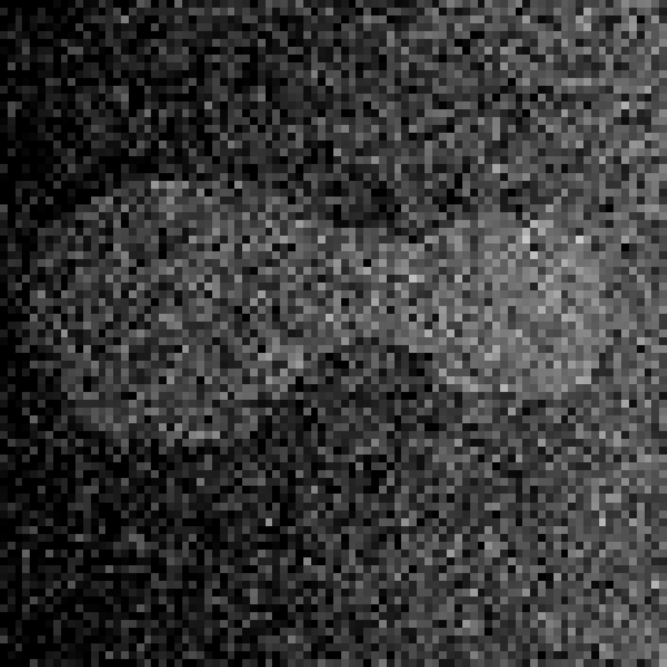}
	\centerline{(d)}
	\end{minipage}
	}
	\hspace{-10.25mm}
	\subfigure{
	\begin{minipage}[b]{0.18\textwidth}
	\includegraphics[width=2cm]{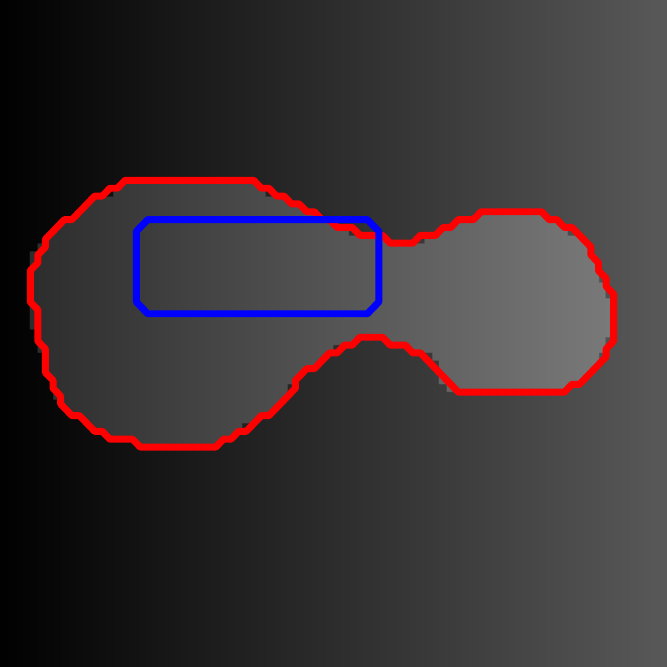}
	\includegraphics[width=2cm]{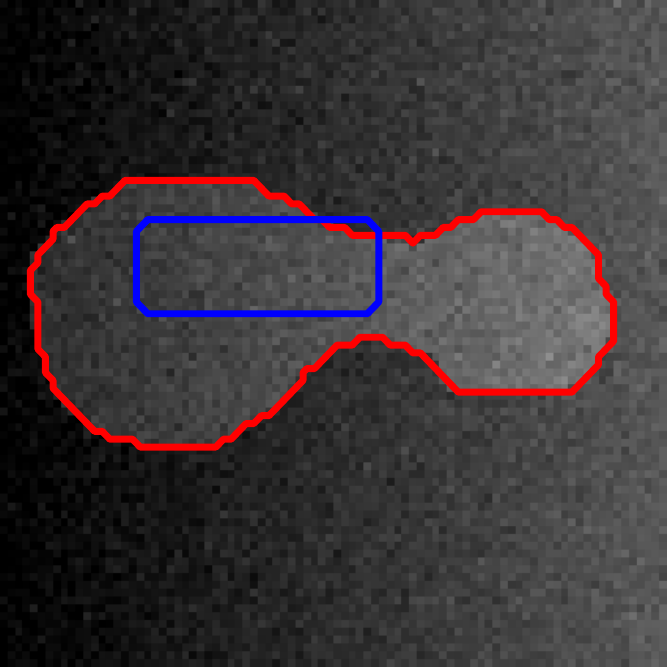}
	\includegraphics[width=2cm]{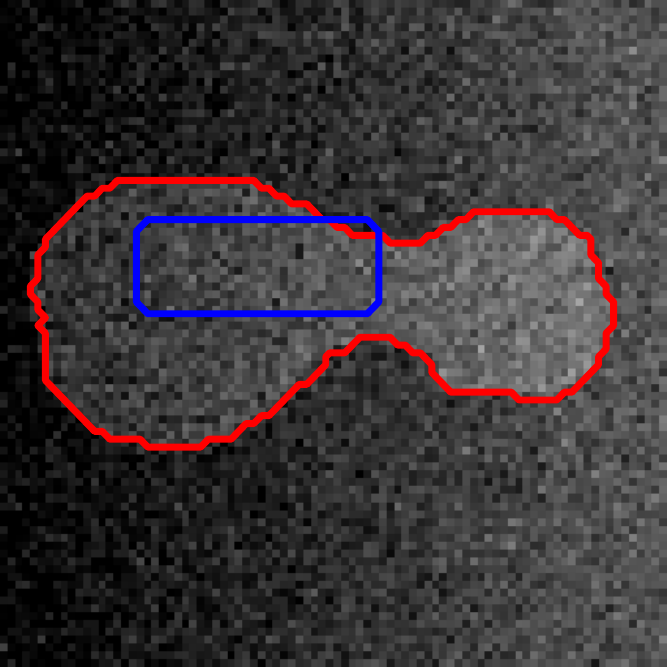}
	\includegraphics[width=2cm]{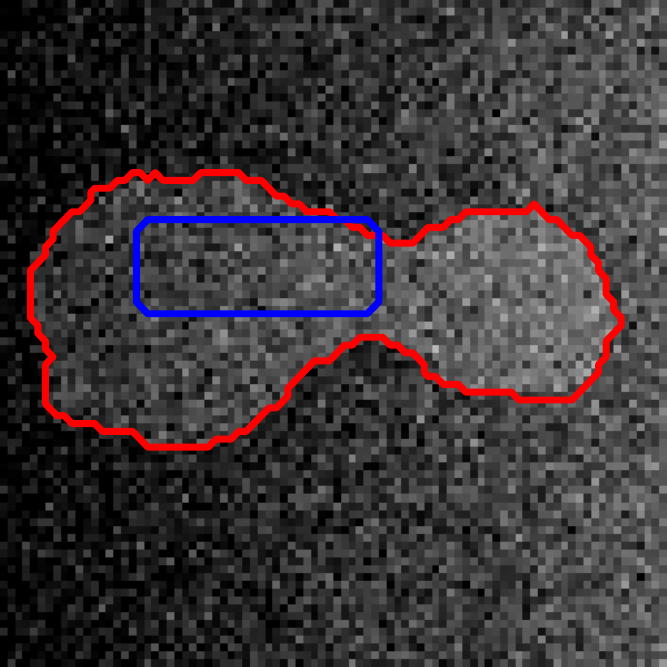}
	\includegraphics[width=2cm]{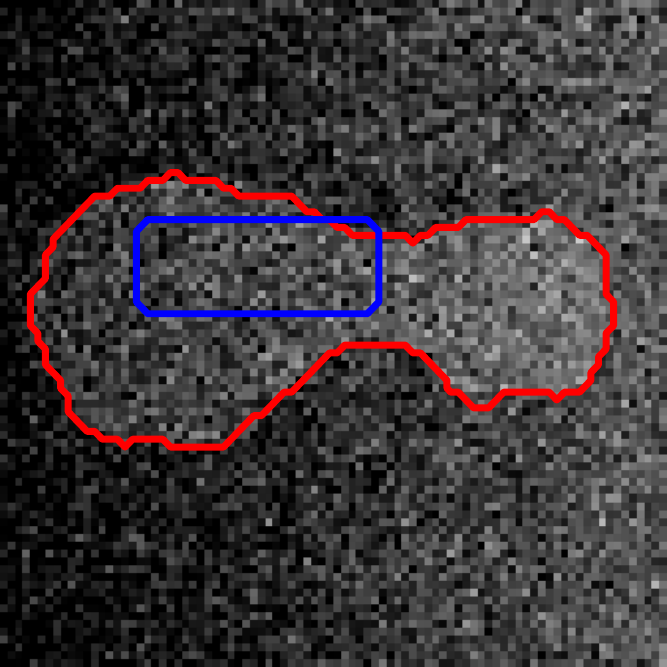}
	\centerline{(e)}
	\end{minipage}
	}
	\hspace{-10.25mm}
	\subfigure{
	\begin{minipage}[b]{0.18\textwidth}
	\includegraphics[width=2cm]{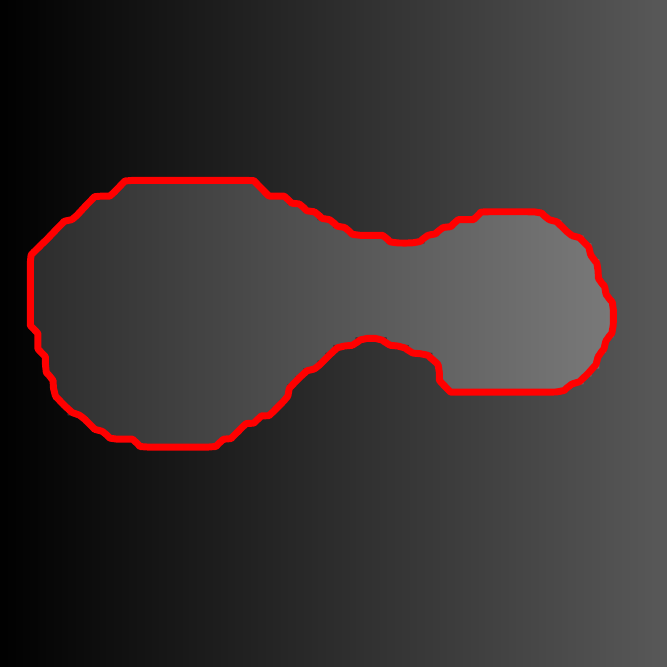}
	\includegraphics[width=2cm]{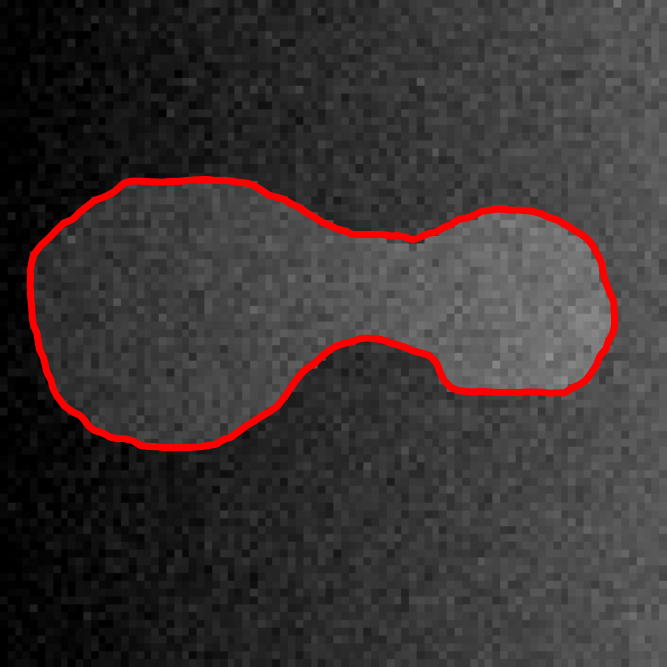}
	\includegraphics[width=2cm]{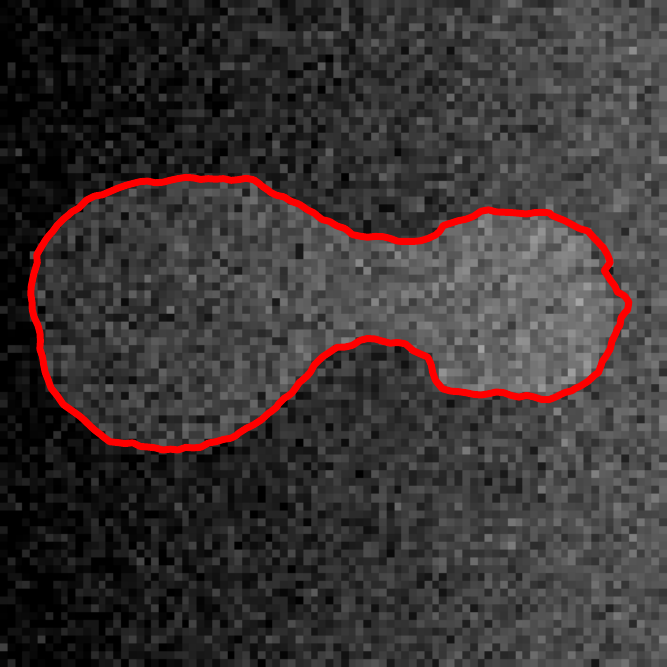}
	\includegraphics[width=2cm]{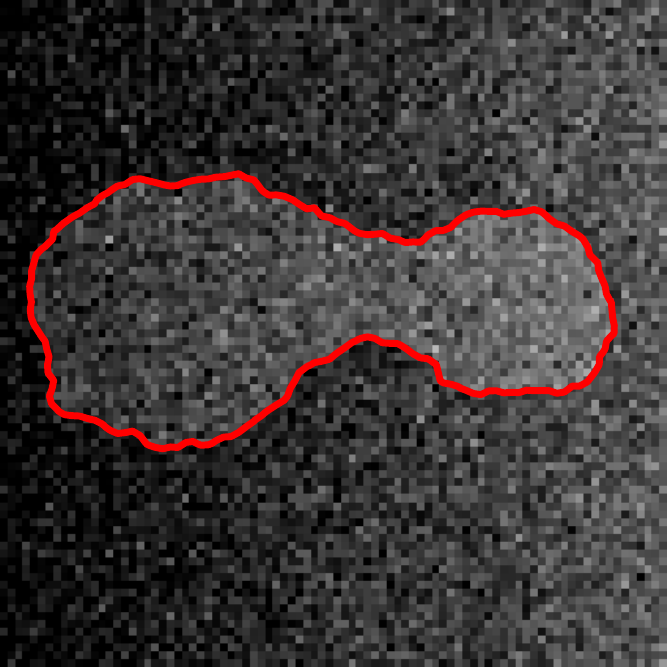}
	\includegraphics[width=2cm]{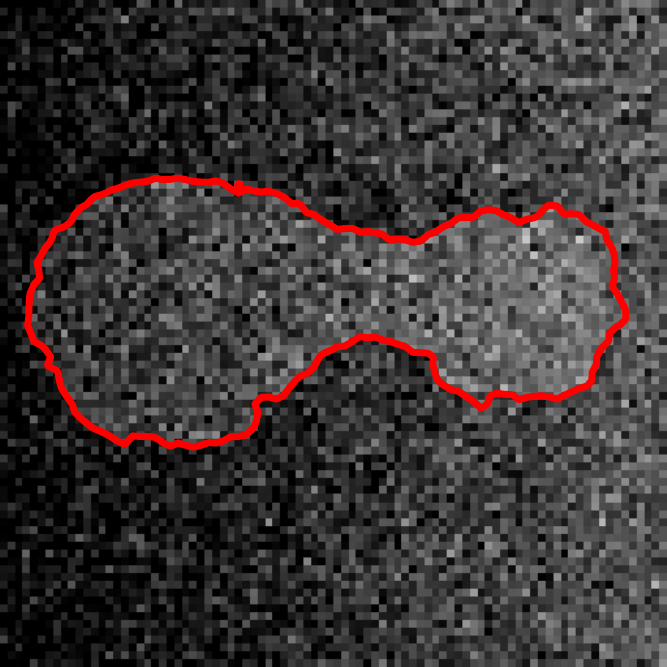}
	\centerline{(f)}
	\end{minipage}
	}
	\caption{
	(a) and (d) the original images with different Gaussian noise level (The variance from the first row to the last row: $100, 200, 300, 500$ for (a) and $0.001*255^2, 0.005*255^2, 0.01*255^2, 0.015*255^2$ for (d)); (b) and (e) corresponding segmentation results produced by the RLSF model; (c) and (f) corresponding segmentation results of the ACLBF model.
	}\label{fig4}
\end{figure}

Finally, we evaluate the ACLBF model on three images with both intensity inhomogeneity and strong Gaussian noise. Results of the RLSF model and our model are shown in Fig.\ref{fig4}. Table \ref{table2} records the iteration number and CPU time for the RLSF model and ACLBF model. 
Experiment results reveal that the ACLBF model, which combines the advantages of the LBF model and phase-field term, has a strong capability of segmenting images with intensity inhomogeneity and noise.
{ One can see that the iteration number of the ETDRK2 scheme is generally less than that of the ETD1 scheme for each segmentation result. In addition, the ETDRK2 scheme also performs better in terms of CPU time.}

\begin{table}[htbp]
\centering
  \begin{tabular}{llllllll}
		\hline
		\multicolumn{2}{l}{\multirow{2}{*}{Images}} &
		\multicolumn{2}{l}{RLSF} & \multicolumn{2}{l}{ {ACLBF-ETD1}}& \multicolumn{2}{l}{ {ACLBF-ETDRK2}} \cr\cmidrule(lr){3-4}
		\cmidrule(lr){5-6}\cmidrule(lr){7-8}
		\multicolumn{2}{l}{} &Ite.&Time(s)&{ Ite.}&{ Time(s)}&Ite.&Time(s)\cr
		\hline
		\multirow{5}{*}{(a)}
		&Row1 &550 &5.016371  &{ 18}&{ 0.168486} &13&0.137160\cr
		&Row2 &800 &7.218534  &{ 19}&{ 0.158388} &15&0.156374\cr
		&Row3 &460 &4.192593  &{ 18}&{ 0.171649} &12&0.136397\cr
		&Row4 &360 &3.265598  &{ 13}&{ 0.122694} &12&0.126147 \cr
		&Row5 &700 &7.312978  &{ 21}&{ 0.162694} &14&0.153064 \cr\hline
		\multirow{5}{*}{(d)}
		&Row1 &290 &2.060543  &{ 14}&{ 0.124064} &8&0.115751\cr
		&Row2 &240 &1.654115  &{ 21}&{ 0.147375} &12&0.138367\cr
		&Row3 &240 &1.700579  &{ 36}&{ 0.198531} &22&0.153823\cr
		&Row4 &270 &1.980383  &{ 68}&{ 0.262130} &34&0.184207 \cr
		&Row5 &270 &1.899098  &{ 58}&{ 0.232933} &39&0.205767 \cr\hline
	\end{tabular}
	\label{table2}\caption{Iteration number and CPU time for experiments in Fig.\ref{fig4}.}
\end{table}

\section{Conclusion}
\label{sec5}
In this paper, we proposed a novel IGLIM for the initial edge detection. Wherein the inhomogeneous graph Laplacian operator can be regarded as an anisotropic Laplacian operator that can recognize most edges of images with intensity inhomogeneity. And the noise-removal method is applied to remove part of irrelevant noise in the initial contour. Then to achieve a better segmentation of images with noise, we adopt a phase-field approach to the LBF model. Based on IGLIM, our proposed method avoids artificial selection of the initial value and obtains a satisfactory segmentation result by solving the derived Allen-Cahn equation. Besides, the ETD schemes we adopt have energy stability.  The ETDRK2 method usually gives better results with less CPU time in the simulation. Numerical experiments exhibit that this model with IGLIM can handle various images effectively and efficiently. Comparison made to other models shows the necessity of the IGLIM and the strong ability of our phase-field approach for segmenting images with noise. In fact, our IGLIM and phase-field approach can be applied to many other models that have issues with robustness on the initialization and noise. Currently, our method can only solve two-phase image segmentation problems because it is difficult to judge which phase the edges obtained by IGLIM belong to when dealing with multi-phase images. In the future, we will employ some classification techniques to generalize our IGLIM to multi-phase image segmentation problems.

\section*{Acknowledgments}
 {We thank the anonymous referees very much for their very insightful remarks and beneficial suggestions.} We thank Prof. Sijie Niu at University of Jinan for providing the MATLAB codes of \cite{niu2017robust} and thank Prof. Tieyong Zeng at The Chinese University of Hong Kong for providing the MATLAB codes of \cite{min2021inhomogeneous}. We also thank Dr. Tingting Wu at Nanjing University of Posts and Telecommunications for comments on the manuscript.

\bibliographystyle{plain}
\bibliography{aclif}

\end{document}